%% file: main.tex
\documentclass{article}

\usepackage[utf8]{inputenc}
\usepackage{amsmath, amsfonts,amsthm}
\allowdisplaybreaks
\usepackage[left=1in,top=1in,right=1in,bottom=1in,letterpaper]{geometry}
\usepackage{multirow,booktabs,makecell}
\usepackage{graphicx,epstopdf,subcaption}
\usepackage{tikz}
\usetikzlibrary{positioning}
\usepackage{hyperref}
\usepackage{algorithm,algorithmic}
\usepackage{color}
\usepackage{natbib}

\input{notations}

\title{Convergence Analysis and Acceleration of Fictitious Play for General Mean-Field Games via the Best Response}
\author{
Jiajia Yu\thanks{The Department of Mathematics, Duke University, Durham, NC.} \and
Xiuyuan Cheng\footnotemark[1] \and
Jian-Guo Liu\footnotemark[1] \thanks{The Department of Physics, Duke University, Durham, NC.} \and
Hongkai Zhao\footnotemark[1]
}

\date{}
\begin{document}

\maketitle
\begin{abstract}

A mean-field game (MFG) seeks the Nash Equilibrium of a game involving a continuum of players, where the Nash Equilibrium corresponds to a fixed point of the best-response mapping. However, simple fixed-point iterations do not always guarantee convergence. Fictitious play is an iterative algorithm that leverages a best-response mapping combined with a weighted average.
Through a thorough study of the best-response mapping, this paper develops a simple and unified convergence analysis, providing the first explicit convergence rate for the fictitious play algorithm in MFGs of general types, especially non-potential MFGs. We demonstrate that the convergence and rate can be controlled through the weighting parameter in the algorithm, with linear convergence achievable under a general assumption.
Building on this analysis, we propose two strategies to accelerate fictitious play. The first uses a backtracking line search to optimize the weighting parameter, while the second employs a hierarchical grid strategy to enhance stability and computational efficiency. We demonstrate the effectiveness of these acceleration techniques and validate our convergence rate analysis with various numerical examples.


\end{abstract}

\section{Introduction}
\label{sec: intro}

The mean-field game (MFG) studies the Nash Equilibrium of a non-cooperative game involving a continuum of players. 
Since Larsy and Lions \cite{lasry2007mfg} and Huang, Malham\'{e} and Caines \cite{huang2006mfg} introduced the model independently, researchers from various fields studied the application of MFGs in, for example, economics \cite{achdou2022income}, epidemics \cite{lee2021controlling}, crowd motion \cite{achdou2019mean,degond2014large}, data science \cite{weinan2019mean,guo2019learning,carmona2023model}, etc.
Recently, people in the machine learning community have been interested in connections between MFG and generative models \cite{zhang2023meanfieldgenerative,Lin2021apacnet,huang2023nfmfg}, and between MFG and reinforcement learning \cite{yang2018mean,guo2019learning,lauriere2022scalable}. 

In an MFG, each player seeks to optimize a strategy that minimizes their individual cost, in response to a given state distribution of the entire player population. 
The individual best strategies collectively shape a new state distribution of the player population. 
The equilibrium of an MFG is the fixed point of this iterative process. 
More precisely, at the Nash Equilibrium, the strategy that optimizes each player’s cost also generates the equilibrium state distribution. 
This equilibrium can be characterized by a set of partial differential equations, comprising a Hamilton-Jacobi-Bellman (HJB) equation for determining each player's optimal control, alongside a Fokker-Planck (FP) equation that depicts the players' state distribution.
Specifically, the solution of the MFG system is the pair of $(\rho,\phi)$ that solves
\begin{subequations}
\begin{align}
    -\partial_t\phi-\nu\Delta\phi+H(x,\nabla\phi(x,t))=f(x,\rho(t)),&
    \quad \phi(x,T)=f_T(x,\rho(T)),   \\
    \partial_t\rho-\nu\Delta\rho-\nabla\cdot(\rho D_pH(x,\nabla\phi))=0,& 
    \quad \rho(x,0)=\rho_0(x).   
\end{align}
\end{subequations}
Despite its diverse applications, there are many difficulties in numerically solving the MFG system.
Besides the nonlinearity in the HJB equation, the main challenge is the intertwined structure of the system.
Specifically, the value function, $\phi$,  of the HJB equation is contingent on the FP solution, $\rho$, and the gradient of the HJB solution drives the FP equation. Furthermore, the HJB equation is evolving backward in time while the FP equation is evolving forward in time, essentially a two-point boundary value problem, which prevents the direct implementation of a time-marching solver.

This paper focuses on an approach known as fictitious play \cite{cardaliaguet2017ficplay}.
This strategy decouples the system and facilitates us to solve the HJB and FP alternatively by simple time-marching methods.
It has a natural interpretation in game theory and economics. 
In each iteration, the players seek the best strategy in response to the current state distribution. They then update the state distribution by taking a weighted average of the current distribution and the distribution resulting from the best strategy, which guides the next strategy search.

An important study of fictitious play is its convergence analysis. 
For potential mean-field games (MFGs), the convergence has been established with the help of the potential functional, as shown in \cite{cardaliaguet2017ficplay}, and the sublinear and linear convergence rates for different $\delta_k$ weights were recently proven in \cite{Lavigne2023gcg-ficplay}.
However, for non-potential MFGs, since a variational structure is absent, the convergence analysis is more challenging and less understood. 

In this paper, we establish the convergence rate of fictitious play for general mean-field games (MFGs) through an analysis of the best-response mapping. 
First, we reveal an interesting geometric relation between the current state distribution and best response distribution it leads to: they lie on the "opposite side" with respect to the equilibrium state, assuming the monotonicity of the interaction cost (see Proposition \ref{lem: stability}). 
This observation explains why direct fixed-point iteration may fail and why fictitious play yields a more stable and accurate approximation to the equilibrium.
Next, we study the stability of the best response.
Although the best-response objective is typically convex but not uniformly strongly convex, we show that it exhibits directional strong convexity between any two distinct function pairs in the feasible set (see Section \ref{subsec: aspt, cvx}). This observation ensures the stability of the best response mapping.
We then define the \textit{gain of the best response} as the difference between the cost of the given distribution and the minimal cost achieved by the best response. 
This is also known as exploitability in \cite{hadikhanloo2017anonymous}.
By definition, the gain of the best response is always non-negative and vanishes only when the given distribution corresponds to a Nash Equilibrium (see Lemma \ref{lem: NE - dual gap=0}).
Our main Theorem \ref{thm: non-pot cvg rate} shows that the gain of the best response acts as a Lyapunov function and decays to zero under the directional convexity of the dynamic cost and monotonicity of the interaction cost. 
It therefore indicates the convergence and provides a convergence rate of fictitious play for general MFGs.
To the best of our knowledge, this is the first convergence rate established for fictitious play without assuming the existence of a potential.

To accelerate convergence in fictitious play, we introduce a backtracking line search method to dynamically select the largest feasible weight $\delta_k$ at each iteration. 
Our convergence analysis shows that the convergence rate is explicitly influenced by the choice of $\delta_k$.
While traditional fictitious play algorithms often use a predefined decay rate, such as $\delta_k = \frac{1}{k}$, this approach results in slower convergence as $k$ increases.
In contrast, our backtracking line search approach adapts $\delta_k$ to maximize efficiency at each step. 
While previous work has applied similar line search techniques for potential MFGs using potential costs \cite{Lavigne2023gcg-ficplay}, our method is independent of potential structures and is therefore suitable for general MFGs.

In addition, we propose a hierarchical grid approach to accelerate the algorithm. 
Recognizing that the MFG system is a two-boundary value problem, we interpret the fictitious play as a bi-directional shooting method. 
To expedite the propagation of boundary conditions, we initiate fictitious play on a coarse grid and progressively refine it.
Specifically, once fictitious play on a coarse grid converges to a criterion proportional to the grid size, the resulting solution serves as an initial guess for the next finer grid. 
This process continues iteratively, using the solution on each grid as a starting point for the next until the finest grid is reached for a more accurate solution. 
Additionally, by incorporating a numerical viscosity term proportional to the grid size in the Hamiltonian discretization, this hierarchical grid approach automatically generates a vanishing viscosity solution, which also stabilizes the fictitious play when physical viscosity is absent.

We summarize our contributions as follows:
\begin{enumerate}
    \item We connect the Nash Equilibrium and fictitious play with the best response mapping. The best response yields a Lyapunov function, and the stability analysis of the best response provides an explicit convergence rate for the fictitious play algorithm for a general class of mean-field games (MFGs). To the best of our knowledge, this is the first convergence rate established for non-potential MFGs.
    \item We show that the convergence rate essentially depends on the averaging weight and propose an acceleration technique that optimizes the weight via a backtracking line search.
    \item We introduce a hierarchical grid method to enhance the stability and efficiency of fictitious play for general MFGs.
\end{enumerate}

\paragraph{Related work on decoupling approaches:}
Different ways to solve for the optimal control and to update the density distribution lead to different implementations of fictitious play. 
In this paper, we focus on alternatively solving HJB and FP equations. 
Among the existing literature, the work presented in \cite{inoue2023ficplaylsmfg} closely aligns with our investigation. This study utilizes the Hopf-Cole transformation technique to linearize the MFG system, which was initially introduced in \cite{gueant2012mean}.
It then employs fictitious play combined with the finite difference method to tackle the linearized system. While \cite{gueant2012mean} focuses on systems with quadratic Hamiltonians, \cite{inoue2023ficplaylsmfg} extends the methodology to MFGs with generalized quadratic Hamiltonians.
In contrast, our paper adopts a direct discretization of the original MFG system, iterating between solving the HJB and FP equations. This approach enables us to work with general Hamiltonians.

When the state space is finite \cite{gomes2013continuous,Gao2024finitemfg}, many works adopt reinforcement learning strategies to obtain the optimal control instead of solving the HJB equation. 
While this reinforcement learning can be inaccurate, \cite{elie2020convergence} studies the error propagation in fictitious play and shows that fictitious play is still effective even when the model is not known and there are accumulated errors in the learning process.
For more studies using reinforcement learning in fictitious play, we refer readers to the review paper \cite{laurière2024learningreview}.

Besides fictitious play, when the Hamiltonian is quadratic, \cite{gueant2012mean} applies the Hopf-Cole transformation to linearize the system and then implements fixed-point iteration. They also prove that this reformulation of fixed-point iteration converges. 
For the general Hamiltonian, another strategy adapts the policy iteration in control theory to the MFG setting. It solves the HJB equation in two phases: initially freezing the action to linearize the HJB for solving the value function, followed by action adjustment based on this value function. 
\cite{cacase2021policyitercvg,camilli2022policyitercvg,camilli2022policyitercvgrate} establish the convergence rate of the algorithm for polynomial growth Hamiltonian and \cite{Laurière2023nonseparableHam} for non-separable Hamiltonian.
\cite{tang2024smpolicyiter} proposes and proves the convergence of smoothed policy iterations, which use an average of the policy to update the density in policy iteration.

\paragraph{Numerical methods for MFGs involving multi-grid vs. hierarchical grid:}
Multi-grid is an effective preconditioner that can accelerate an iterative solver for elliptic equations. The key idea is that the iterative solver is effective in reducing errors in high frequencies (relative to the grid size). So one can use multiple grid sizes to deal with all frequencies effectively. 
Some, but not extensive, literature has explored the application of the multi-grid idea to solve an MFG problem or its subproblem.
\cite{achdou2012multigrid} explores accelerating Newton's methods for the coupled global nonlinear discretized system in space and time with multigrid and preconditioning. 
Additionally, \cite{briceno2018multigridadmm,andreev2017multigridadmm} apply the multigrid methods to solve the (bi-)Laplacian equation that arises in the subproblem of optimization algorithms like primal-dual and ADMM. 
All the above use of multi-grid is aimed as a preconditioner to deal with problems of elliptic nature.
However, as noted in \cite{achdou2012multigrid,AchdouLauriere2020review}, for second-order MFGs, the elliptic operator is second-order in time and fourth-order in space, leading to anisotropy that renders a full-coarsening multigrid approach ineffective.
Therefore \cite{achdou2012multigrid,briceno2018multigridadmm,andreev2017multigridadmm} only apply multigrid on the spatial domain for second-order MFGs.
Our hierarchical grid strategy is totally different. The key idea of using a hierarchical grid is to accelerate information propagation from the boundary for an MFG as a two-point boundary value problem. Hence our strategy applies to both first-order and second-order MFGs.
In addition, while multigrid deals with frequencies from high to low and hence goes from fine to coarse grid, our strategy always goes from coarse to fine, using the approximate solution on a coarse grid, which can achieve a coarse balance between boundary conditions at both ends more quickly and stably, as an initial guess to a finer grid to find a more accurate solution. Moreover, our hierarchical strategy automatically incorporates a numerical viscosity that is proportional to the grid size in our finite difference scheme. It allows a natural computation of the vanishing viscosity solution for the best response for first-order MFGs which enhances the numerical stability. 
In this sense, the spirit of the multilevel approach in \cite{yu2024pgdeucmfg} is the closest to ours, where the interpolation of the solution on the coarse grid serves as a good initial guess on the fine grid. However, their method is an optimization approach for the global (in space and time) and coupled discrete system based on the variational formulation. It can not be applied to non-potential MFGs. Also, due to the use of straightforward numerical discretization for the HJB equation, the stability can be a problem for first-order MFGs or MFGs with small diffusion.  
To the best of our knowledge, our work is the first one employing the full-coarsening hierarchical grid and combining it with the decoupling approach. 

\paragraph{Other numerical methods for MFGs:} 
Many existing numerical methods for MFG view the system for both PDEs of $\phi$ and $\rho$ in space and time as a whole system and aim to solve it altogether. 
For example, \cite{achdou2010mean} applies finite difference to the MFG system and solves the entire system with Newton solver. 
A more common approach focuses on potential MFG with a variational formulation and applies an optimization solver to it, see \cite{benamou2015augmented,benamou2017variational,liu2021nonlocal,yu2024pgdeucmfg} for low-dimensional cases and \cite{Lin2021apacnet,Ruthotto2020mfgmachinelearning,huang2023nfmfg} for high-dimensional cases. 
While this coupled approach directly handles the global structure of the PDE system, the full discretization in space and time results in a large nonlinear system to solve. Moreover, as the grid size becomes finer, the nonlinear system not only becomes larger but also worse conditioned due to the differential operation.
In general, a direct solver ignoring the physics and structures of the underlying PDE system can be inefficient. 
Besides computational concerns, all optimization strategies rely on the equivalent variational formulation of the problem and are not easy to adapt to non-potential MFG problems.

\paragraph{Organization:} 
We summarize the notations in Table \ref{tab: notation}.
The remainder of the paper is organized as follows.
Section \ref{sec: review mfg} introduces the problem setting and reviews the formulation of mean-field games (MFGs) along with the fictitious play algorithm.
In Section \ref{sec: ficplay cvg}, we reformulate the MFG and fictitious play dynamics using the concept of the \textbf{best response}, analyze the stability of the best response map, and provide a detailed convergence analysis for both potential and non-potential MFGs.
Section \ref{sec: imple} describes our discretization and implementation details.
Section \ref{sec: num result} verifies the convergence theory and demonstrates the effectiveness of the proposed methods through numerical experiments.
Finally, Section \ref{sec: conclu} concludes the paper and discusses directions for future work.

Our main contributions lie in the convergence analysis developed in Section \ref{sec: ficplay cvg}, with a focus on Sections \ref{subsec: best response} and \ref{subsec: non pot mfg}, as well as in the proposed acceleration strategies presented in Sections \ref{subsec: acc btls} and \ref{subsec: multigrid}. The analysis is built upon three key components: the Lyapunov function defined in \eqref{eq: def exploit}, the qualitative property established in Proposition \ref{lem: stability}, and the notion of directional strong convexity in Section \ref{subsec: aspt, cvx}.

For simplicity, we focus on second-order MFGs ($\nu>0$) and periodic boundary conditions in Sections \ref{sec: review mfg} and \ref{sec: ficplay cvg}. 
We believe our results naturally extend to problems posed on $\mathbb{R}^d$ or on a bounded domain with suitable boundary conditions. Although first-order MFGs ($\nu = 0$) involve technical subtleties, we expect our results to hold when the solution is sufficiently smooth and $\rho > 0$ throughout. 
To illustrate this, in Section \ref{sec: imple}, we detail a discretization on a bounded domain with homogeneous Neumann conditions for the HJB equation and no-flux conditions for the Fokker–Planck equation, which are used in the experiments in Section \ref{sec: num result}.

\begin{table}[htbp]
\caption{Summary of Notations}
\label{tab: notation}
\centering
\renewcommand{\arraystretch}{1.2}  
\begin{tabular}{m{14cm}}  
\toprule
$\bbT^d=\bbR^d/\mathbb{Z}^d$ is the space domain, with $L^2$ inner product $\langle\cdot,\cdot\rangle_{\bbT^d}$ and $L^2$ norm $\|\cdot\|_{\bbT^d}$. \\
$T>0$ is the given terminal time. \\
$Q_T:=\bbT^d\times(0,T)$ is the joint domain, with $L^2$ inner product $\langle\cdot,\cdot\rangle_{Q_T}$ and $L^2$ norm $\|\cdot\|_{Q_T}$. \\
$\calP(\bbT^d)$ is the set of Borel probability measures on $\bbT^d$. \\ \hline
$\rho=(\rho_t)_{t\in[0,T]}, \rho_t\in\calP(\bbT^d)$ is the distribution/density evolution in most cases. Depending on context, $\rho\in[0,+\infty)$ denotes the density at one point, or $\rho\in\calP(\bbT^d)$ denotes a probability distribution or its density. \\
$v:\bbT^d\times[0,T]\to\bbR^d$ is the velocity field in most cases. Depending on context, $v\in\bbR^d$ is a vector at one point, or $v:\bbT^d\to\bbR^d$ is the velocity field at some time $t$. \\
$m:\bbT^d\times[0,T]\to\bbR^d,m=\rho v$ is the momentum in most cases. Depending on context, $m\in\bbR^d$ is a vector at some $(x,t)$, or $m:\bbT^d\to\bbR^d$ is the velocity field at some time $t$. \\ \hline
$H:\bbT^d\times\bbR^d \to \bbR,(x,p)\mapsto H(x,p)$ is the Hamiltonian. \\
$L:\bbT^d\times\bbR^d,(x,v)\mapsto L(x,v)$ is the convex conjugate of $H$ defined in \eqref{eq: lag}. \\
$f,f_T$ are interaction and terminal costs, either nonlocal $\bbT^d\times\calP(\bbT^d)\to\bbR$ or local $\bbT^d\times[0,+\infty)\to\bbR$. \\
$C$ is the constraint set of $(\rho,m)$, defined in \eqref{eq: cstr set m}. \\ 
$R:C\to\bbR$ is the dynamic cost, defined in \eqref{eq: def R}. \\
$J(\cdot,\cdot;\rho): C\to\bbR$ is the optimal control cost under given density evolution $\rho$, defined in \eqref{eq: def J(,;rho)}. \\
$g:C\to\bbR$ is the gain of the best response, defined in \eqref{eq: def exploit}. \\ \hline
$F,F_T:\calP(\bbT^d)\to\bbR$ are the potential of $f,f_T$ if $f,f_T$ derive from potential. \\
$J: C\to\bbR$ is the potential cost of a potential MFG, defined in \eqref{eq: obj m}. \\
$h:C\to\bbR$ is the primal gap in potential MFG, defined in \eqref{eq: def h}. \\
\bottomrule
\end{tabular}
\end{table}

\section{A Brief Review of MFG and Fictitious Play}
\label{sec: review mfg}

Let the time domain be $[0,T]$ ($T > 0$), the spatial domain be $\mathbb{T}^d$, and define the space-time domain as $Q_T := \bbT^d \times (0,T)$. We denote the $L^2$ inner product and norm on $\bbT^d$ by $\langle\cdot,\cdot\rangle_{\bbT^d}$ and $\|\cdot\|_{\bbT^d}$, and those on $Q_T$ by $\langle\cdot,\cdot\rangle_{Q_T}$ and $\|\cdot\|_{Q_T}$.
$\mathcal{P}(\mathbb{T}^d)$ is the set of Borel probability measures on $\mathbb{T}^d$.
The Hamiltonian $H$, interaction cost $f$, and terminal cost $g$ are periodic in space.
The initial condition $\rho_0 \in \mathcal{P}(\mathbb{T}^d)$ is periodic in space.
When $\rho\in\calP(\bbT^d)$ is absolutely continuous, we slightly abuse notation and use $\rho$ to refer to either a probability distribution or its density function, depending on the context. When $\rho$ refers to the density function, $\rho>0$ implies $\rho(x)>0$ for any $x\in\bbT^d$.

\subsection{Second-order MFGs and monotonicity}

A second-order Mean Field Game (MFG) is the following system of $(\rho,\phi)$ with $\nu>0$:
\begin{equation}\left\{
\begin{aligned}
    -\partial_t\phi(x,t)-\nu\Delta\phi(x,t)+H(x,\nabla\phi(x,t))=f(x,\rho(\cdot,t)),
    \quad &(x,t)\in Q_T,  \\
    \partial_t\rho(x,t)-\nu\Delta\rho(x,t)-\nabla\cdot(\rho D_pH(x,\nabla\phi))(x,t)=0,
    \quad &(x,t)\in Q_T,  \\
    \phi(x,T)=f_T(x,\rho(\cdot,T)),\quad \rho(x,0)=\rho_0(x), 
    \quad &x\in\bbT^d.  \\
\end{aligned}\right.
\label{eq: mfg}
\end{equation}
Here $\rho(x, t)$ represents the population density and $\phi(x, t)$ denotes the value function of a representative player, both at position $x$ and time $t$. 
The value function $\phi$ solves a HJB equation involving $\rho$, indicating the influence of the population on the cost for a single player. Meanwhile, the population density $\rho$ solves a FP equation, where the drift is determined by $\phi$.

The two equations are coupled through the interaction cost $f$ and terminal cost $f_T$.
We consider two common types of couplings $f$ and $f_T$ appear in the MFG literature.
A nonlocal cost $f:\bbT^d \times \calP(\bbT^d)\to\bbR$ is used in the cases where players can influence each other even when they are not close.
For example, $f(x,\rho) := \int_{\bbT^d} K(x,y) \dd\rho(y)$, where $K:\bbT^d \times \bbT^d \to \bbR$ is a given kernel.
Alternatively, one may consider local interaction costs, where $f:\bbT^d \times [0,+\infty) \to \bbR$ takes the form $f(x,\rho) = f(x, \rho(x))$ for absolutely continuous distributions $\rho$. A typical example is $f(x,\rho)=\log(\rho(x))$.

The notion of a weak solution in mean-field games (MFGs) varies depending on the setting, such as whether the cost is local or nonlocal, or whether the system is first-order or second-order. Correspondingly, the assumptions required for the existence of weak solutions also differ across these cases. For a comprehensive discussion, we refer the reader to \cite{Achdou2021MFGbook,cardaliaguet2010notes,Gomes2014survey} and the references therein.
For simplicity, we assume that the system under consideration admits at least one classical solution.

To ensure the uniqueness of the solution, monotonicity of interaction and terminal cost is necessary in general (\cite{Achdou2021MFGbook} Remark 1.4). It also plays an important role in our convergence analysis.

\begin{definition}[Monotonicity \cite{lasry2007mfg}]
\label{def: monotone}
    A nonlocal cost $f:\bbT^d\times\calP(\bbT^d)$ is monotone if
    for any $\rho,\rho'\in\calP(\bbT^d)$,
    $$ \int_{\bbT^d}(f(x,\rho)-f(x,\rho'))\dd (\rho-\rho')(x) \geq 0. $$
    $f$ is strictly monotone if the strict inequality holds for any $\rho\neq\rho'$.
    For a local cost $f:\bbT^d\times[0,+\infty)$, we say $f$ is monotone if $f(x,\cdot)$ is nondecreasing and $f$ is strictly monotone if $f(x,\cdot)$ is increasing for any $x\in\bbT^d$.

    If both $f$ and $f_T$ are monotone, we call the system \eqref{eq: mfg} a monotone MFG. 

\end{definition}

If the solution to a monotone MFG exists and $H$ is strictly convex on $p$, then the solution is unique 
(\cite{Achdou2021MFGbook} Theorem 1.4, Theorem 1.11 and Theorem 1.12).

\subsection{Fictitious Play}

Directly solving \eqref{eq: mfg} is difficult because $\rho$ and $\phi$ are intertwined and the HJB is backward in time while the FP is forward in time.
A natural idea to resolve this issue is to solve the HJB and FP equations alternatively. This leads to a fixed-point iteration.
Precisely, during the $k$-th iteration,  one calculates $\phi^{(k+1)}$ by solving the HJB equation with $\rho=\rho^{(k)}$, and then utilize $-D_pH(x,\nabla\phi^{(k+1)})$ to update $\rho^{(k+1)}$ through the FP equation. 
However, Proposition \ref{lem: stability} and the numerical experiment in Section \ref{subsec: num fixpt} indicate that fixed-point iteration can be unstable, potentially oscillating between two states and failing to converge.
Instead of using the latest population density as guidance, the fictitious play \cite{cardaliaguet2017ficplay} uses a weighted average of \textbf{history} population densities.
Specifically, when solving for the value function $\hat{\phi}^{(k+1)}$, an weighted average of history and current density, $\rho^{(k)}=\left(1-\delta_k\right)\rho^{(k-1)}+\delta_k\hat{\rho}^{(k)} ~ (0\leq\delta_k\leq 1)$, is utilized, rather than the current density $\hat{\rho}^{(k)}$. 
We summarize the procedure in Algorithm \ref{alg: ficplay}.

\begin{algorithm}[htbp]
\caption{Fictitious Play}
\begin{algorithmic}
    \STATE{Parameters} $\rho_0,f,f_T,\nu,0 \leq \delta_k \leq 1$
    \STATE{Initialization} $\rho^{(0)}$
    \FOR{$k=1,2,\cdots,K$}
    
        \STATE{Solve the HJB \eqref{eq: ficplay hjb} for $\hat{\phi}^{(k)}$}
        
        \begin{equation}
        \left\{
        \begin{aligned}
            &-\partial_t\hat{\phi}^{(k)}-\nu\Delta\hat{\phi}^{(k)}+H(x,\nabla\hat{\phi}^{(k)}(x,t))=f(x,\rho^{(k-1)}(t)),\quad (x,t)\in Q_T, \\
            &\hat{\phi}^{(k)}(x,T)=f_T(x,\rho^{(k-1)}(T)),\quad x\in\bbT^d.
        \end{aligned}
        \right.
        \label{eq: ficplay hjb}
        \end{equation}
        \STATE{Solve the FP \eqref{eq: ficplay fp} for $\hat{\rho}^{(k)}$}
        \begin{equation}
        \left\{
        \begin{aligned}
            &\partial_t\hat{\rho}^{(k)}-\nu\Delta\hat{\rho}^{(k)}-\nabla\cdot(\hat{\rho}^{(k)} D_pH(x,\nabla\hat{\phi}^{(k)}))=0, \quad (x,t)\in Q_T,\\
            & \hat{\rho}^{(k)}(x,0)=\rho_0(x), \quad x\in\bbT^d.
        \end{aligned}
        \right.
        \label{eq: ficplay fp}
        \end{equation}

        \STATE{Execute density average to obtain $\rho^{(k)}$}
        \begin{equation}
            \rho^{(k)}=\left(1-\delta_k\right)\rho^{(k-1)}+\delta_k\hat{\rho}^{(k)}.
        \label{eq: ficplay avg}
        \end{equation}
    \ENDFOR
    \STATE{Output} $\rho^{(K)},\hat{\phi}^{(K)}.$
\end{algorithmic}
\label{alg: ficplay}
\end{algorithm}

\begin{remark}[Connection to the Folk Theorem \cite{wikipedia_folk_theorem}]
The use of a weighted average of historical population densities in fictitious play, rather than relying solely on the current density, emphasizes the role of history in shaping strategic decisions. This mirrors the insight of the Folk Theorem in repeated games, where long-term outcomes are sustained not only by current actions but also by the accumulated history of play. In both settings, past behavior influences future responses, enabling the emergence of equilibria that may not be supported in one-shot interactions. This connection highlights how memory and learning mechanisms can foster stable outcomes, even in decentralized or adaptive settings.
\end{remark}

The following assumptions are enforced throughout the paper.
\begin{assumption} For the given initial distribution, Hamiltonian $H$, interaction cost $f$ and terminal cost $f_T$, we make the following assumptions: 
    \begin{enumerate}
        \item The initial probability distribution $\rho_0$ is absolutely continuous with respect to the Lebesgue measure, and has a $C^2$ positive density, also denoted as $\rho_0$.

        \item $H$ is $C^2$ in $x$, smooth and strictly convex in $p$, and satisfies $$\lim_{|p|\to\infty}\frac{H(x,p)}{|p|}=+\infty,\quad\text{for any } x\in\bbT^d.$$

        \item $f$ satisfies one of the following assumptions:
        \begin{enumerate}
            \item $f$ is nonlocal, and $f(\cdot,\rho)\in C^2(\bbT^d)$ for any $\rho\in\calP(\bbT^d)$.
            \item $f$ is local, and $f\in C^2(\bbT^d\times[0,+\infty))$.
        \end{enumerate}

        \item $f_T$ satisfies one of the following assumptions:
        \begin{enumerate}
            \item $f_T$ is nonlocal, and $f_T(\cdot,\rho)\in C^3(\bbT^d)$ for any $\rho\in\calP(\bbT^d)$.
            \item $f_T$ is local, and $f_T\in C^3(\bbT^d\times[0,+\infty))$.
        \end{enumerate}
    \end{enumerate}
\label{asp: throughout}
\end{assumption}
Assumptions \ref{asp: throughout} do not guarantee the existence of a classical solution to the MFG system \eqref{eq: mfg}. However, standard results on semilinear parabolic equations ensure that they are sufficient to guarantee a unique classical solution at each step of the fictitious play algorithm, leading to the following lemma.
\begin{lemma}
\label{lem: ficplay welldefined}
    Under Assumption \ref{asp: throughout}, if $\rho^{(0)}\in C^2(\overline{Q_T})$, $\rho^{(0)}(\cdot,t)\in\calP(\bbT^d)$, $\rho^{(0)}(\cdot,0) = \rho_0$, and $\rho^{(0)}>0$, then Algorithm \ref{alg: ficplay} satisfies that for $k=1,2,\cdots$,
    \begin{enumerate}
        \item \eqref{eq: ficplay hjb} has a unique classical solution $\hat{\phi}^{(k)}$.
        \item \eqref{eq: ficplay fp} has a unique classical solution $\hat{\rho}^{(k)}$ and $\hat{\rho}^{(k)}>0$.
        \item $\rho^{(k)}>0$.
    \end{enumerate}
\end{lemma}

\subsection{Potential MFGs}

In the literature, the assumption that the MFG has potential is key to proving the convergence of fictitious play. 

\begin{definition}[Potential]
    A nonlocal cost $f:\bbT^d\times \calP(\bbT^d)\to\bbR$ or local cost $f:\bbT^d\times[0,+\infty)\to\bbR$ derive from potential if there exists $F:\calP(\bbT^d)\to\bbR$ such that 
    for any $\rho,\rho'\in\calP(\bbT^d)$,
    $$ \lim_{s \to 0} \frac{F(\rho +s(\rho'-\rho))-F(\rho)}{s} = \int_{\bbT^d} f(x,\rho) \dd(\rho'-\rho)(x). $$
    A MFG is a potential MFG if the interaction and terminal costs $f,f_T$ derive from potentials.
\end{definition}
When the game is potential, the monotonicity of $f$ and $f_T$ is linked to the convexity of the potentials $F$ and $F_T$.
Along with the convexity of $H$, \cite{lasry2007mfg} points out that the MFG system \eqref{eq: mfg} can be connected to the optimality condition of an optimization problem.
The associated objective function then serves as a Lyapunov function in the convergence analysis.

To be precise, we first introduce a change of variable $(\rho,m)=(\rho,\rho v)$ and define the convex constraint set $C$
\begin{equation}
    C:=\left\{(\rho,m)\middle |
    \begin{aligned}
    &\rho\in C([0,T],L^1(\bbT^d)\cap\calP(\bbT^d)); \\
    &m\in L^1(\bbT^d\times[0,T];\bbR^d);\\
    &\rho(\cdot,0) = \rho_0, \, 
    \partial_t\rho-\nu\Delta\rho+\nabla\cdot m = 0
    \end{aligned}\right\}.
\label{eq: cstr set m}
\end{equation}
Here, the equations are satisfied in the sense of distributions.
Let $L$ be the Lagrangian
\begin{equation}
    L(x,v):=\sup_p\left\{ -\langle p,v \rangle - H(x,p) \right\},
\label{eq: lag}
\end{equation}
then $L$ is differentiable and strictly convex in $v$ by smoothness and strictly convexity of $H$ in $p$ (Assumption \ref{asp: throughout}).
We further define the dynamic cost $R$ as
\begin{equation}
    R({\rho},{m}):=\int_0^T\int_{\bbT^d} {\rho}(x,t)L\left(x,\frac{{m}(x,t)}{{\rho}(x,t)}\right)\dd x \dd t
\label{eq: def R}
\end{equation}
where when $\rho=0$, $L\left( x,\frac{m}{\rho} \right)$ is defined by convention 
\begin{equation}
    L\left( x,\frac{m}{\rho} \right) := \begin{cases}
        0, & m=0,\\ +\infty, & m\neq 0.
    \end{cases}
\end{equation}
It is easy to verify that $R$ is jointly convex in $(\rho,m)\in C$.
When $f,f_T$ derive from potentials $F,F_T$, denote the objective as 
\begin{equation}
    J(\rho,m):= R(\rho,m) + \int_0^T F(\rho(\cdot,t))\dd t + F_T(\rho(\cdot,T)).
\label{eq: obj m}
\end{equation}
Then $J$ is jointly convex in $(\rho,m)$ if $F$ and $F_T$ are also convex.
With the above notations, we have the following proposition for potential MFGs.
\begin{proposition}[Potential MFGs]
\label{thm: potential mfg}
    Assume that $f$ and $f_T$ derive from potentials $F$ and $F_T$, and $F, F_T$ are convex under $L^2$ metric. 
    If $(\rho^*,\phi^*)$ is of class $C^2(\bbT^d\times(0,T))$ with $\rho^*(x,0)=\rho_0(x)$, $\rho^*>0$,
    then $(\rho^*,\phi^*)$ is a solution to \eqref{eq: mfg} 
    if and only if 
    $m^*(x,t):=-\rho^*D_pH(x,\nabla\phi^*(x,t))$ and 
    \begin{equation}
    J(\rho^*,m^*) = \min_{(\rho,m)\in C} J(\rho,m)
    \end{equation}
\end{proposition}
\begin{proof}
$(\rho^*,\phi^*)$ being a solution to \eqref{eq: mfg} is sufficient due to the convexity of $F$ and $F_T$ and by the verification strategy used later in Lemma \ref{lem: prepare}. We omit the details here.

For the necessity, $(\rho^*,m^*)$ is the minimizer and therefore $(\rho^*,m^*)\in C$. In addition, $m^*=-\rho^* D_pH(\cdot,\nabla\phi^*)$. 
To show that $(\rho^*,\phi^*)$ solves the HJB equation in \eqref{eq: mfg},
we define
\begin{equation}
    \mathcal{L}(\rho,m;\phi^*):= 
    J(\rho,m) - \int_0^T\int_{\bbT^d}\phi^*(x,t)(\partial_t\rho - \nu\Delta\rho + \nabla\cdot m)(x,t)\dd x\dd t.
\end{equation}
Then for any $(\rho,m)\in C$, $\mathcal{L}(\rho,m;\phi^*) = J(\rho,m)$, and the optimality of $(\rho^*,m^*)$ gives
\begin{equation}
    \mathcal{L}(\rho^*+s(\rho-\rho^*),m^*+s(m-m^*);\phi^*) - \mathcal{L}(\rho^*,m^*;\phi^*)\geq 0.
\label{eq: rho m opt in calL}
\end{equation}
Because $\rho^*>0$, we have that for any $(\rho,m)\in C$ and $s\in[0,1]$, 
\begin{equation}
\begin{aligned}
    &R(\rho^*+s(\rho-\rho^*),m^*+s(m-m^*)) - R(\rho^*,m^*)\\
    =& s\int_0^T\int_{\bbT^d} (\rho-\rho^*)\left( L\left(x,\frac{m^*}{\rho^*}\right) - \frac{m^*}{\rho^*}D_vL\left(x,\frac{m^*}{\rho^*}\right) \right)(x,t)\dd x \dd t \\
    &+ s\int_0^T\int_{\bbT^d} (m-m^*)\cdot D_vL\left( x,\frac{m^*}{\rho^*}\right)(x,t) \dd x \dd t + o(s)
\end{aligned}
\label{eq: linearize R}
\end{equation}
Therefore, we have that for any $(\rho,m)\in C$ and $s\in[0,1]$,
\begin{equation}
\begin{aligned}
    &\mathcal{L}(\rho^*+s(\rho-\rho^*),m^*+s(m-m^*);\phi^*) - \mathcal{L}(\rho^*,m^*;\phi^*)\\
    =& s\int_0^T\int_{\bbT^d} (\rho-\rho^*)\left( \partial_t\phi^*+\nu\Delta\phi^* + L\left(x,\frac{m^*}{\rho^*}\right) - \frac{m^*}{\rho^*}D_vL\left(x,\frac{m^*}{\rho^*} \right) + f(x,\rho^*_t)\right)(x,t)\dd x\dd t  \\
    &+ s\int_0^T\int_{\bbT^d} (m-m^*) \left( D_vL\left( x,\frac{m^*}{\rho^*} \right) + \nabla\phi^* \right)(x,t)\dd x\dd t\\
    &+ s \int_{\bbT^d}(\rho-\rho^*) \Big(f_T(x,\rho^*_T) - \phi^*\Big)(x,T)\dd x + o(s)\\
    =& -s\int_0^T\int_{\bbT^d} (\rho-\rho^*)\Big( -\partial_t\phi^* - \nu\Delta\phi^* + H(x,\nabla\phi^*) - f(x,\rho^*_t)\Big)(x,t)\dd x\dd t  \\
    &+ s \int_{\bbT^d}(\rho-\rho^*) \Big(f_T(x,\rho^*_T) - \phi^*\Big)(x,T)\dd x + o(s),
\end{aligned}
\end{equation}
where the first equality is by definition of $\mathcal{L}$ and potential, together with \eqref{eq: linearize R} and integration by part, and the second equality is because $\frac{m^*}{\rho^*} = -D_pH(\cdot,\nabla\phi^*)$.
Then by optimality \eqref{eq: rho m opt in calL}, we conclude that $(\rho^*,\phi^*)$ solves the HJB equation.

\end{proof}
\begin{remark}
    For simplicity, we assume that $(\rho^*, \phi^*) \in C^2(\bbT^d \times (0,T))$ in Proposition \ref{thm: potential mfg}. For more general conditions under which the MFG system \eqref{eq: mfg} can be connected to a variational formulation, we refer the reader to Section 1.3.7.2 and Chapter 2 of \cite{Achdou2021MFGbook}, as well as the references therein.
\end{remark}

Since the potential formulation with monotone costs has a convex objective and linear constraint, many first-order optimization algorithms such as proximal gradient \cite{yu2024pgdeucmfg}, primal-dual \cite{liu2021nonlocal} and augmented Lagrangian \cite{benamou2015augmented,benamou2017variational} have been adapted to solve the problem besides fictitious play.
While potential MFGs enjoy many desired properties, there are also limitations.
For example, introducing $m=\rho v$ makes the Lipschitz constant of the objective function extremely large when $\rho$ is close to 0 and causes algorithms such as proximal gradient \cite{yu2024pgdeucmfg} to be unstable.
And there are MFGs not potential, for example $f(\cdot,\rho)=K*\rho$ where $K$ is not a symmetric kernel.

The convergence proofs in \cite{cardaliaguet2017ficplay,Lavigne2023gcg-ficplay} rely heavily on the decay of the functional $J(\rho, m)$, which is specific to potential games and does not extend to the general non-potential setting. One of the main contributions of this paper is to establish a convergence analysis of fictitious play for a broad class of MFGs without assuming the existence of a potential. 
In the potential case, \cite{cardaliaguet2017ficplay} proves that setting $\delta_k = \frac{1}{k+1}$ yields convergence, while \cite{Lavigne2023gcg-ficplay} connects fictitious play with the generalized conditional gradient method and analyzes how different choices of $\delta_k$ affect convergence rates. 
However, for non-potential games, the influence of $\delta_k$ on convergence has not been previously addressed. In Section \ref{sec: ficplay cvg}, we provide an explicit convergence rate estimate in terms of $\delta_k$ and offer practical guidance on choosing $\delta_k$ to accelerate convergence in the non-potential setting.

\section{Convergence and Acceleration of Fictitious Play}
\label{sec: ficplay cvg}

Given the fictitious play algorithm (Algorithm \ref{alg: ficplay}), several questions arise: Does it converge? Under what assumptions does convergence occur? And if it does, what is the rate of convergence? 
For potential MFGs, the objective function serves as a natural criterion for studying convergence.
Based on this, \cite{cardaliaguet2017ficplay} demonstrates the convergence of fictitious play for both first-order and second-order potential games, and \cite{Lavigne2023gcg-ficplay} offers a sublinear convergence of classical fictitious play through the lens of the generalized conditional gradient. \cite{Lavigne2023gcg-ficplay} also shows that a linear convergence rate can be achieved when the weight $\delta_k$ is chosen wisely.

In this section, we establish the linear convergence of fictitious play for non-potential MFGs. 
To that end, Section \ref{subsec: equilibrium} introduces the key concept of the best response. This allows us to reformulate the fictitious play algorithm and define the best response gain as a measure of the distance to the Nash equilibrium. 
In Section \ref{subsec: aspt, cvx}, we present the assumptions and definitions underlying the convergence analysis, including directional strong convexity as a relaxation of global strong convexity, and discuss their validity.
Section \ref{subsec: best response} then analyzes the stability of the best-response mapping, a central step in understanding the convergence behavior of fictitious play.
Leveraging these stability properties, we show that the gain of the best response functions is a Lyapunov function and prove the convergence of fictitious play for non-potential MFGs; details are provided in Section \ref{subsec: non pot mfg}.
Section \ref{subsec: acc btls} presents the backtracking line search for optimizing the weight $\delta_k$ to accelerate the algorithm. 
Finally, Section \ref{subsec: pot mfg} remarks on the convergence for potential MFGs, relating the gain of the best response to the dual gap in the generalized conditional gradient (Frank-Wolfe) algorithm.

\subsection{The best response and the equilibrium formulation}
\label{subsec: equilibrium}

When there is no ambiguity, we omit the input $x$ of $L,H,f$ and $f_T$.

We define the cost of $(\tilde{\rho},\tilde{m})$ given $\rho$ as 
\begin{equation}
J(\tilde{\rho},\tilde{m};\rho):=R(\tilde{\rho},\tilde{m}) 
+ \langle f(\rho),\tilde{\rho} \rangle_{Q_T} 
+ \langle f_T(\rho_T),\tilde{\rho}_T \rangle_{\bbT^d},
\label{eq: def J(,;rho)}
\end{equation}
where
\begin{equation*}
    \langle f(\rho),\tilde{\rho} \rangle_{Q_T}:=\int_0^T \int_{\bbT^d} f(x,\rho_t)\dd \tilde{\rho}_t(x)\dd t,\quad
    \langle f_T(\rho_T),\tilde{\rho}_T \rangle_{\bbT^d}:=\int_{\bbT^d} f_T(x,\rho_T)\dd\tilde{\rho}_T(x),
\end{equation*}
Here, the interaction cost and the terminal cost are linear in $\tilde{\rho}$.
We first connect the fictitious play update \eqref{eq: ficplay hjb}-\eqref{eq: ficplay fp} to an optimal control problem 

\begin{lemma} 
\label{lem: prepare}
For any $\rho \in \calP(\bbT^d)$, the optimal control problem 
\begin{equation*}
    \min_{(\tilde{\rho},\tilde{m})\in C} J(\tilde{\rho},\tilde{m};\rho)
\end{equation*}
admits a unique minimizer 
\begin{equation*}
    (\hat{\rho},\hat{m}) = (\hat{\rho},-\hat{\rho}D_pH(\cdot,\nabla\hat{\phi})),
\end{equation*}
where $(\hat{\rho},\hat{\phi})$ solves the decoupled system
\begin{subequations}
\label{eq: decoupled hjb, fp}
\begin{equation}
\left\{
\begin{aligned}
    -\partial_t\hat{\phi}(x,t) - \nu\Delta\hat{\phi}(x,t) + H(x,\nabla\hat{\phi}(x,t)) &= f(x,\rho(\cdot,t)), && (x,t) \in Q_T \\
    \hat{\phi}(x,T) &= f_T(x,\rho(\cdot,T)), && x \in \mathbb{T}^d
\end{aligned}
\right.
\label{eq:decoupled-system-a}
\end{equation}

\begin{equation}
\left\{
\begin{aligned}
    \partial_t\hat{\rho}(x,t) - \nu\Delta\hat{\rho}(x,t) - \nabla\cdot(\hat{\rho} D_p H(x,\nabla \hat{\phi}))(x,t) &= 0, && (x,t) \in Q_T \\
    \hat{\rho}(x,0) &= \rho_0(x), && x \in \mathbb{T}^d
\end{aligned}
\right.
\label{eq:decoupled-system-b}
\end{equation}
\end{subequations}

\end{lemma}

\begin{proof}
Lemma \ref{lem: ficplay welldefined} gives the existence of $(\hat{\rho},\hat{\phi})$.
We prove that $(\hat{\rho},\hat{\phi})$ is the unique minimizer by verification. Without loss of generality, consider $(\tilde{\rho},\tilde{m})\in C$ such that $\tilde{m}=\tilde{\rho}\tilde{v}$ for some $\tilde{v}$ and $(\tilde{\rho},\tilde{m})\neq(\hat{\rho},\hat{m})$.
Denote $\hat{v}:=-D_pH(\cdot,\nabla\hat{\phi})$ and $\hat{p}:=\nabla\hat{\phi}$.
Then by definition, 
\begin{equation}
    L(x,\hat{v}) = -\langle \hat{p},\hat{v} \rangle - H(x,\hat{p}),
\end{equation}
and by strictly convexity of $H$ in the Assumption \ref{asp: throughout}, for any $v\neq \hat{v}$, 
\begin{equation}
    L(x,v) > -\langle \hat{p},v \rangle - H(x,\hat{p}).
\end{equation}
Therefore, we have
\begin{equation}
\begin{aligned}
    R(\tilde{\rho},\tilde{m}) - R(\hat{\rho},\hat{m})
    =& \int_{0}^{T}\int_{\bbT^d} \tilde{\rho}L(x,\tilde{v}) - \hat{\rho}L(x,\hat{v})\dd x\dd t \\
    >&- \left\langle \nabla\hat{\phi}, \tilde{m}-\hat{m} \right\rangle_{Q_T}
    - \left\langle H(\nabla\hat{\phi}), \tilde{\rho}-\hat{\rho} \right\rangle_{Q_T} 
\end{aligned}
\label{eq: 14}
\end{equation}
and as a result,
\begin{equation}
\begin{aligned}
    J(\tilde{\rho},\tilde{m};\rho) - J(\hat{\rho},\hat{m};\rho)
    =&R(\tilde{\rho},\tilde{m}) - R(\hat{\rho},\hat{m})
    + \langle f(\rho),\tilde{\rho}-\hat{\rho} \rangle_{Q_T}
    + \langle f_T(\rho_T),\tilde{\rho}_T-\hat{\rho}_T \rangle_{\bbT^d}\\
    >&\left\langle -H(\nabla\hat{\phi})+f(\rho), \tilde{\rho}-\hat{\rho} \right\rangle_{Q_T} 
    - \left\langle \nabla\hat{\phi}, \tilde{m}-\hat{m} \right\rangle_{Q_T}
    + \langle f_T(\rho_T),\tilde{\rho}_T-\hat{\rho}_T \rangle_{\bbT^d}\\
    =&\left\langle -\partial_t\hat{\phi}-\nu\Delta\hat{\phi}, \tilde{\rho}-\hat{\rho} \right\rangle_{Q_T} 
    + \left\langle \hat{\phi}, \nabla\cdot(\tilde{m}-\hat{m}) \right\rangle_{Q_T}
    + \langle f_T(\rho_T),\tilde{\rho}_T-\hat{\rho}_T \rangle_{\bbT^d}\\
    =&\left\langle \hat{\phi}, (\partial_t-\nu\Delta)(\tilde{\rho}-\hat{\rho})+\nabla\cdot(\tilde{m}-\hat{m}) \right\rangle_{Q_T} = 0,
\end{aligned}
\label{eq: 13}
\end{equation}
where the last two equalities are due to \eqref{eq: decoupled hjb, fp}.
Therefore, we prove that $(\hat{\rho},\hat{m})$ is the unique minimizer.

\end{proof}

Since for given $\rho$, the minimizer to $J(\tilde{\rho},\tilde{m};\rho)$ always exists and is unique, we define the minimizer $(\hat{\rho},\hat{m})$ as the best response of $\rho$.

\begin{definition}[Best response and the gain of the best response]
For any $(\rho,m)\in C$, we define $(\hat{\rho},\hat{m})$ as the best response of $\rho$ if it solves the minimization problem
\begin{equation}
(\hat{\rho},\hat{m})=\argmin_{(\tilde{\rho},\tilde{m})\in C} J(\tilde{\rho},\tilde{m};\rho),
\label{eq: def best response}
\end{equation}
and define the gain of the best response as the gap between $J$ evaluated at $(\rho,m)$ and its best response:
\begin{equation}
    g(\rho,m):= J(\rho,m;\rho) - \min_{(\tilde{\rho},\tilde{m})\in C} J(\tilde{\rho},\tilde{m};\rho)
    =J(\rho,m;\rho) - J(\hat{\rho},\hat{m};\rho).
\label{eq: def exploit}
\end{equation}
\end{definition}

\begin{remark}
Our definition of ``best response'' is related to but different from the ``best reply strategy'' in game theory. The relation between MFG and best reply strategy is studied in \cite{degond2014large,degond2017mfgbrs}.
Additionally, our concept of ``the gain of the best response'' aligns with the notion of ``exploitability'' in related literature. 
For instance, \cite{hadikhanloo2017anonymous} introduces this concept for first-order monotone MFGs and uses it to prove convergence. 
Similarly, \cite{geist2022curlmfg} employs exploitability to demonstrate first-order convergence of fictitious play in continuous-time discrete MFGs.
\end{remark}

The best response is the central concept to analyze the convergence of fictitious play for non-potential MFGs.
Lemma \ref{lem: prepare} shows that the fictitious play (Algorithm \ref{alg: ficplay}) is equivalent to iteratively solving for the best response and computing the weighted average of the replies. To be precise, it is equivalent to the following scheme,
\begin{align}
    &(\hat{\rho}^{(k)},\hat{m}^{(k)})=\argmin_{(\tilde{\rho},\tilde{m})\in C} J(\tilde{\rho},\tilde{m};\rho^{(k-1)}),
    \label{eq: ficplay lin}\\
    &(\rho^{(k)},m^{(k)})=(1-\delta_k)(\rho^{(k-1)},m^{(k-1)}) + \delta_k (\hat{\rho}^{(k)},\hat{m}^{(k)}),
    \label{eq: ficplay avg rhom}
\end{align} 
and the $\hat{\phi}^{(k)}$ in Algorithm \ref{alg: ficplay} is the Lagrangian multiplier of \eqref{eq: ficplay lin}.
In addition, a mean-field Nash Equilibrium is the fixed point of the best response mapping, or equivalently, the root of the gain of the best response, as detailed in the following Lemma.
\begin{lemma}
\label{lem: NE - dual gap=0}
    If $(\rho^*,m^*)\in C$ and $\rho^*>0$, then the following statements are equivalent,
    \begin{enumerate}
        \item $(\rho^*,m^*)=(\rho^*,-\rho^*D_pH(\cdot,\nabla\phi^*))$ where $(\rho^*,\phi^*)$ solves \eqref{eq: mfg}
        \item $(\rho^*,m^*)=\argmin_{(\tilde{\rho},\tilde{m})\in C} J(\tilde{\rho},\tilde{m};\rho^*)$,
        \item $g(\rho^*,m^*)=0$.
    \end{enumerate}
\end{lemma}

\begin{proof}
    $1 \Leftrightarrow 2$ is from Lemma \ref{lem: prepare}, and $2 \Leftrightarrow 3$ is by definition.
\end{proof}

In the rest of this paper, we refer to \eqref{eq: mfg} the PDE formulation, 
\begin{equation}
    \inf_{(\rho,m)\in C} J(\rho,m)
\label{eq: pot mfg}
\end{equation}
the variational formulation and 
\begin{equation}
    (\rho,m)=\argmin_{(\tilde{\rho},\tilde{m})\in C} J(\tilde{\rho},\tilde{m};\rho) \qquad \text{ or } \qquad g(\rho,m)=0
\label{eq: eq mfg}
\end{equation}
the equilibrium formulation of an MFG.

\subsection{Directional strong convexity and other assumptions}
\label{subsec: aspt, cvx}

The following assumptions and definitions appear frequently in our analysis.
\begin{assumption}
    $f$ satisfies one of the following assumptions
    \begin{enumerate}
        \item $f$ is nonlocal and $\|f(\cdot,\rho)-f(\cdot,\rho')\|_{\bbT^d}\leq L_f \|\rho-\rho'\|_{\bbT^d}$ for any $\rho,\rho'\in\calP(\bbT^d)$.
        \item $f$ is local and $|f(x,\rho)-f(x,\rho')|\leq L_f|\rho-\rho'|$ for any $x\in\bbT^d$ and $\rho,\rho'\in[0,+\infty)$.
    \end{enumerate}
\label{aspt: lip cts}
\end{assumption}
\begin{assumption}
    $f_T(x,\rho_T)=f_T(x)$ is independent of $\rho_T$.
\label{aspt: independent cost}
\end{assumption}

\begin{remark}[Independent terminal cost]
    The assumption that $f_T(x,\rho_T)$ being independent of $\rho_T$ is quite common in MFG literature, see \cite{cacase2021policyitercvg,camilli2022policyitercvgrate,Achdou2013fdcvg,achdou2012multigrid,Lavigne2023gcg-ficplay,inoue2023ficplaylsmfg} for example.
    When this assumption holds, the terminal condition of the HJB equation is independent of $\rho_T$. But as long as $f(x,\rho)$ depends on $\rho$, the MFG system \eqref{eq: mfg} is still coupled and the problem is not oversimplified. 
\end{remark}

\begin{definition}[Strong monotonicity]
\label{def: str monotone}
    $f:\bbT^d\times\mathcal{P}(\bbT^d)$ is $\lambda_f$-strongly monotone if $\lambda_f>0$ and for any $\rho,\rho'\in\mathcal{P}(\bbT^d)\cap L^2(\bbT^d)$,
    \begin{equation*}
        \int_{\bbT^d}(f(x,\rho)-f(x,\rho')) (\rho-\rho')(x) \dd x \geq \lambda_f \int_{\bbT^d} |\rho(x)-\rho'(x)|^2\dd x.
    \end{equation*}
\end{definition}

\begin{remark}[Strong monotonicity of $f$]
    $f$ being strongly monotone is a common assumption in game theory and is related to the price of anarchy \cite{Gao2024finitemfg}.
\end{remark}

The next definition is the key to convergence rate estimation.
\begin{definition}[Directional strong convexity]
\label{def: cvx of r}
    For given $(\rho^i,m^i)\in C,i=0,1$, denote $({\rho}^\theta,{m}^\theta):=
    (1-\theta)({\rho}^0,{m}^0) + \theta({\rho}^1,{m}^1)$, and 
    \begin{equation}
        y(\theta):=R\left( \rho^\theta,m^\theta \right),\theta\in[0,1].
    \label{eq: def y(epsilon)}
    \end{equation}
    We define that $R$ is directional $\lambda_R$-strongly convex between $(\rho^0,m^0)$ and $(\rho^1,m^1)$ if $\lambda_R>0$ and $y$ is $(\lambda_R\|\rho^0-\rho^1\|_{Q_T}^2)$-strongly convex in $[0,1]$.
\end{definition}

By the strong convexity of $y$, if $\rho^0>0,\rho^1>0$ and $R$ is directional $\lambda_R$-strongly convex between $(\rho^0,m^0)$ and $(\rho^1,m^1)$, then
\begin{equation*}
    y'(0) + \frac{\lambda_R}{2}\|\rho^0-\rho^1\|_{Q_T}^2
    \leq R(\rho^1,m^1) - R(\rho^0,m^0),
\end{equation*}
where 
\begin{equation*}
    y'(\theta) = \left\langle L\left(x,\frac{m^\theta}{\rho^\theta}\right)
    -\frac{m^\theta}{\rho^\theta}D_vL\left(x,\frac{m^\theta}{\rho^\theta}\right), \rho^1-\rho^0 \right\rangle_{Q_T}
    + \left\langle D_vL\left(x,\frac{m^\theta}{\rho^\theta}\right),m^1-m^0 \right\rangle_{Q_T}.
\end{equation*}

We remark that the standard definition of strong convexity requires the existence of a uniform $\lambda_R$ such that the function is $\lambda_R$-convex along any $(\rho^0,m^0)$ and $(\rho^1,m^1)$. 
However, the dynamic cost $R$ does not satisfy this strict definition even when $L$ is strongly convex. 
In contrast, our proof only relies on the local existence of $\lambda_R$, which holds for any distinct pair $(\rho^0, m^0), (\rho^1, m^1) \in C$ under mild conditions.
To illustrate this, assume that $L$ is $C^2$ and uniformly $\lambda_L$-convex in $v$ ($\lambda_L>0$).
For any $(\rho^0,m^0),(\rho^1,m^1)\in C$, and $\rho^0>0,\rho^1>0$, denote 
\begin{equation*}
    D_{vv}^{\theta}(x,t):= D_{vv}L\left(x,\frac{m^\theta(x,t)}{\rho^\theta(x,t)}\right).
\end{equation*}
Then
\begin{equation}
\begin{aligned}
    y''(\theta)=&\int_{Q_T}\frac{1}{(\rho^{\theta})^3} 
    [ (\rho^1-\rho^0)m^{\theta} - \rho^{\theta}(m^1-m^0) ]^{\top} D_{vv}^{\theta}[ (\rho^1-\rho^0)m^{\theta} - \rho^{\theta}(m^1-m^0) ]   \\
    \geq &\lambda_L\int_{Q_T}\frac{(\rho^0\rho^1)^2}{(\rho^{\theta})^3}\left| \frac{m^0}{\rho^0}-\frac{m^1}{\rho^1} \right|^2
    \geq \lambda_L \int_{Q_T}\min\left\{ \frac{(\rho^0)^2}{\rho^1},\frac{(\rho^1)^2}{\rho^0} \right\}\left| \frac{m^0}{\rho^0}-\frac{m^1}{\rho^1} \right|^2.
\end{aligned}
\end{equation}
Since $\rho^0,\rho^1>0$, the lower bound is zero only if 
\begin{equation}
\frac{m^0}{\rho^0} = \frac{m^1}{\rho^1} \quad \text{a.e.}.
\label{eq: m/rho =}
\end{equation}
Notice that $\rho^\theta$ initialize from $\rho_0$ and is driven by $\frac{m^\theta}{\rho^\theta}$, \eqref{eq: m/rho =} implies $(\rho^0, m^0) = (\rho^1, m^1)$ a.e.
Therefore, given the assumption that $\rho^0 > 0$ and $\rho^1 > 0$, we may assume, without loss of generality, that
\begin{equation*}
\int_{Q_T}\min\left\{ \frac{(\rho^0)^2}{\rho^1},\frac{(\rho^1)^2}{\rho^0} \right\}\left| \frac{m^0}{\rho^0}-\frac{m^1}{\rho^1} \right|^2 = \underline{c}>0.
\end{equation*}
As a result, we obtain $y''(\theta) \geq \lambda_L \underline{c} > 0$ for all $\theta \in [0,1]$, implying that $R$ is directionally strongly convex between $(\rho^0, m^0)$ and $(\rho^1, m^1)$. 
This also highlights that the value of $\lambda_R$ may become small, posing numerical challenges when $\min\left\{ \frac{(\rho^0)^2}{\rho^1}, \frac{(\rho^1)^2}{\rho^0} \right\}$ is uniformly small.

On the other hand, if $f$ is $2\lambda_f$-strongly monotone, we can modify the definition of $R$ accordingly. Let $\lambda_R = \lambda_f$ and define the modified cost function as
\begin{equation*}
    \tilde{R}(\tilde{\rho},\tilde{m}):= R(\tilde{\rho},\tilde{m})+\frac{\lambda_R}{2}\|\tilde{\rho}\|_{Q_T}^2,\quad
    \tilde{f}(x,\rho):= f(x,\rho) - \lambda_R\rho(x), \quad
    \tilde{J}(\tilde{\rho},\tilde{m};\rho):=\tilde{R}(\tilde{\rho},\tilde{m}) + \langle \tilde{f}(\rho),\tilde{\rho} \rangle_{Q_T} + \langle f_T,\tilde{\rho}_T\rangle_{\bbT^d}.
\end{equation*}
With this formulation, the modified cost $\tilde{R}$ is $\lambda_R$-convex between any $(\rho^0, m^0), (\rho^1, m^1) \in C$, and $\tilde{f}$ is $\lambda_f$-monotone. The corresponding best response is given by
\begin{equation}
    (\hat{\rho},\hat{m})\in\argmin_{(\tilde{\rho},\tilde{m})\in C} \tilde{J}(\tilde{\rho},\tilde{m};\rho),
    \label{eq: best response mdf}
\end{equation} 
which leads to the modified update system:
\begin{equation}
\begin{aligned}
    -\partial_t\hat{\phi}(x,t)-\nu\Delta\hat{\phi}(x,t)+H(x,\nabla\hat{\phi}(x,t))=f(x,\rho(\cdot,t)) + \lambda_f(\hat{\rho}-\rho),
    \quad &(x,t)\in Q_T,  \\
    \partial_t\hat{\rho}(x,t)-\nu\Delta\hat{\rho}(x,t)-\nabla\cdot(\hat{\rho}D_pH(x,\nabla\hat{\phi}))(x,t)=0,
    \quad &(x,t)\in Q_T,  \\
    \hat{\phi}(x,T)=f_T(x,\rho(T)),\quad \hat{\rho}(x,0)=\rho_0(x), 
    \quad &x\in\bbT^d.    
\end{aligned}
\end{equation}
When $\rho$ is close to the Nash equilibrium, the difference $\hat{\rho} - \rho$ becomes small. In this regime, the fictitious play update \eqref{eq: ficplay hjb}–\eqref{eq: ficplay fp} approximately solves the modified best response problem \eqref{eq: best response mdf}, where directional strong convexity of $R$ is guaranteed. Notably, the value of $\lambda_R$ is constant and independent of $\rho$ and $m$.

\subsection{The stability of the best response}
\label{subsec: best response}

This section examines the stability of the best response mapping. 
Throughout, we assume there exists a Nash Equilibrium $(\rho^*,m^*,\phi^*)$, such that  $(\rho^*,m^*)=\argmin_{(\tilde{\rho},\tilde{m})\in C}J(\tilde{\rho},\tilde{m};\rho^*)$ and $m^*=-\rho^*D_pH(\nabla\phi^*)$.

We first give a qualitative description of the best response mapping.
\begin{proposition}
\label{lem: stability}
    For given $\rho^i,i=0,1$, let $(\hat{\rho}^i,\hat{m}^i)=\argmin_{(\tilde{\rho},\tilde{m})\in C} J(\tilde{\rho},\tilde{m};\rho^i)$, if Assumption \ref{aspt: independent cost} holds, then 
    \begin{equation}
        \langle \hat{\rho}^0-\hat{\rho}^1, f(\rho^0)-f(\rho^1) \rangle_{Q_T} + D = 0,
    \label{eq: lem stb 1}
    \end{equation}
    where $D\geq 0$ and $D=0$ if and only if $\hat{\rho}^0=\hat{\rho}^1$.
\end{proposition}
\begin{proof}
    From Lemma \ref{lem: prepare} and optimality of $\hat{\rho}^i$, there exist $\hat{\phi}^i$ such that
    \begin{equation}
    \begin{aligned}
        &\langle \hat{\rho}^0-\hat{\rho}^1, f(\rho^0)-f(\rho^1) \rangle_{Q_T}\\
        =&\underbrace{\langle \hat{\rho}^0-\hat{\rho}^1, -\partial_t(\hat{\phi}^0-\hat{\phi}^1) -\nu\Delta(\hat{\phi}^0-\hat{\phi}^1) \rangle_{Q_T}}_{A} \\
        &- \langle \hat{\rho}^0, H(\nabla\hat{\phi}^1)-H(\nabla\hat{\phi}^0) \rangle_{Q_T}
        - \langle \hat{\rho}^1, H(\nabla\hat{\phi}^0)-H(\nabla\hat{\phi}^1) \rangle_{Q_T}
    \end{aligned}
    \end{equation}
    From the boundary conditions and FP equations, we have
    \begin{equation}
    \begin{aligned}
        A=&\langle \hat{\phi}^0-\hat{\phi}^1, \nabla\cdot(\hat{\rho}^0D_pH(\nabla\hat{\phi}^0)) - \nabla\cdot(\hat{\rho}^1D_pH(\nabla\hat{\phi}^1))\rangle_{Q_T}\\
        =&\langle \hat{\rho}^0, D_pH(\nabla\hat{\phi}^0)\nabla(\hat{\phi}^1-\hat{\phi}^0) \rangle_{Q_T} 
        + \langle \hat{\rho}^1, D_pH(\nabla\hat{\phi}^1)\nabla(\hat{\phi}^0-\hat{\phi}^1) \rangle_{Q_T} 
    \end{aligned}
    \end{equation}
    Let 
    \begin{equation}
    \begin{aligned}
        D:=&\langle \hat{\rho}^0, H(\nabla\hat{\phi}^1)-H(\nabla\hat{\phi}^0)-D_pH(\nabla\hat{\phi}^0)\nabla(\hat{\phi}^1-\hat{\phi}^0)\rangle_{Q_T}\\
        &+\langle \hat{\rho}^1, H(\nabla\hat{\phi}^0)-H(\nabla\hat{\phi}^1)-D_pH(\nabla\hat{\phi}^1)\nabla(\hat{\phi}^0-\hat{\phi}^1)\rangle_{Q_T}.
    \end{aligned}
    \end{equation}
    Combining the above gives \eqref{eq: lem stb 1} and by strict convexity of $H$, $D\geq0$. 
    If $\hat{\rho}^0=\hat{\rho}^1$ then by definition $D=0$.
    If $D=0$, then by positivity of $\hat{\rho}^0$ and $\hat{\rho}^1$, we have $\nabla\hat{\phi}^0=\nabla\hat{\phi}^1$.
    Since $\hat{\rho}^0$ and $\hat{\rho}^1$ start from the same distribution and are driven by the same vector field $D_pH(\nabla\hat{\phi}^0)=D_pH(\nabla\hat{\phi}^1)$, we have that $\hat{\rho}^0=\hat{\rho}^1$.
\end{proof}

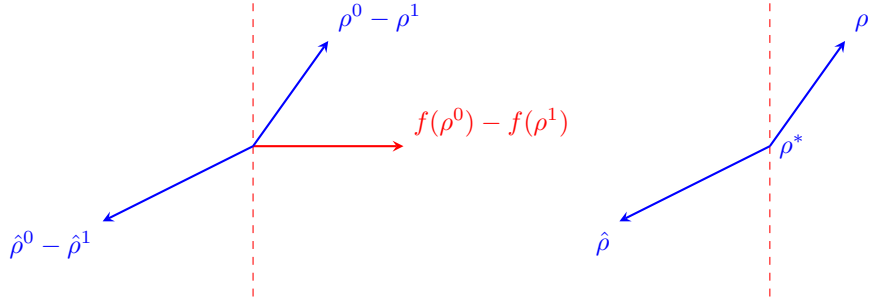
\begin{figure}[htbp]
    \centering
    \begin{tikzpicture}[scale=2, >=stealth]
    \draw[->, thick, red] (0,0) -- (1,0) node[anchor=south west] {$f(\rho^0)-f(\rho^1)$};
    \draw[dashed,red] (0,-1) -- (0,1);
    \draw[->, thick, blue] (0,0) -- (0.5, 0.7) node[anchor=south west] {$\rho^0-\rho^1$};
    \draw[->, thick, blue] (0,0) -- (-1, -0.5) node[anchor=north east] {$\hat{\rho}^0-\hat{\rho}^1$};
    \end{tikzpicture}
    \begin{tikzpicture}[scale=2, >=stealth]
    \draw[dashed,red] (0,-1) -- (0,1);
    \draw[->, thick, blue] (0,0) -- (0.5, 0.7) node[anchor=south west] {$\rho$};
    \draw[->, thick, blue] (0,0) -- (-1, -0.5) node[anchor=north east] {$\hat{\rho}$};
    \node[blue, anchor=west] at (0,0) {$\rho^*$};
    \end{tikzpicture}
    \caption{The illustration of the relations between $\rho^0-\rho^1$, $\hat{\rho}^0-\hat{\rho}^1$ and $f(\rho^0)-f(\rho^1)$ (left) and the relations between $\rho,\hat{\rho}$ and $\rho^*$ (right).}
    \label{fig: illu}
\end{figure}

Proposition \ref{lem: stability} shows that when $f$ is monotone, the vectors $\rho^0 - \rho^1$ and $\hat{\rho}^0 - \hat{\rho}^1$ always lie on opposite sides of the “plane” defined by $f(\rho^0) - f(\rho^1)$, as illustrated in Figure \ref{fig: illu} (left). Furthermore, since $\hat{\rho}^0 > 0$ and $\hat{\rho}^1 > 0$, the vectors $\hat{\rho}^0 - \hat{\rho}^1$ and $f(\rho^0) - f(\rho^1)$ are “orthogonal” if and only if $\hat{\rho}^0 - \hat{\rho}^1$ degenerates.
The experiment in Section \ref{subsec: num cvg} verifies that the quantity $D$ remains nonnegative and reaches zero only when $\hat{\rho} = \rho^*$.

Specializing Proposition \ref{lem: stability} to the case where $\rho^0 = \rho$ is any density evolution and $\rho^1 = \rho^*$ is the Nash equilibrium, Figure \ref{fig: illu} (right) illustrates that $\rho$ and $\hat{\rho}$ always lie on opposite sides of $\rho^*$. 
This geometric interpretation helps explain why repeated application of the best response mapping may fail to converge to the Nash equilibrium $\rho^*$. The numerical experiment in Section \ref{subsec: num fixpt} confirms this behavior, showing that fixed-point iteration can oscillate between two states and fail to converge.

While Proposition \ref{lem: stability} gives a qualitative description, to understand the best response mapping and establish the convergence rate for fictitious play, we need a more careful stability analysis.
To be precise, we next show that the best response mapping is stable in three aspects:
\begin{itemize}
    \item Lemma \ref{lem: best response dist bounded by exploit}: the gain of the best response at $\rho$ is small, 
    only if $\rho$ is close to its best response $\hat{\rho}$.
    \item Lemma \ref{lem: NE < BR}: if $\rho$ is close to its best response, then $\rho$ is close to the fixed point $\rho^*$.
    \item Lemma \ref{lem: optctl map is cts}: if the costs $f(\rho^0)$ and $f(\rho^1)$ are close, then the best responses $\hat{\rho}^0$ and $\hat{\rho}^1$ are close. Specifically, if $f$ is Lipschitz continuous (Assumption \ref{aspt: lip cts}) and $\rho$ is close to the fixed point $\rho^*$, then $\rho$ is close to its best response $\hat{\rho}$.
\end{itemize}
Lemma \ref{lem: NE < BR} and Lemma \ref{lem: optctl map is cts} together give the stability of the fixed point of the best response mapping.
And Lemma \ref{lem: optctl map is cts} also shows the stability of the best response mapping. 
Combining the three lemmas, we have that
\begin{equation*}
    c_{(\rho,m)} g(\rho,m) \geq \|\rho-\hat{\rho}\|_{Q_T} \sim \|\rho-\rho^*\|_{Q_T},
\end{equation*}
for some $c_{(\rho,m)}\in\bbR$.

We first prove that the distance between $\rho$ and its best response $\hat{\rho}$ is bounded by the gain of the best response at $\rho$. 
The key property for this lemma to hold is the constraint being linear.
A similar technique is used in \cite{Lavigne2023gcg-ficplay}.
\begin{lemma}
\label{lem: best response dist bounded by exploit}
    For given $(\rho,m)\in C$, let its best response be $(\hat{\rho},\hat{m})$. 
    For any $(\tilde{\rho},\tilde{m})\in C$, 
    if $R$ is directional $\lambda_R$-strongly convex between $(\tilde{\rho},\tilde{m})$ and $(\hat{\rho},\hat{m})$,
    then 
    \begin{equation}
        J(\tilde{\rho},\tilde{m};\rho)-J(\hat{\rho},\hat{m};\rho)\geq\frac{\lambda_R}{2}\|\hat{\rho}-\tilde{\rho}\|_{Q_T}^2.
    \label{eq: lem best response 1}
    \end{equation}
    As a consequence, if $R$ is directional $\lambda_R$-strongly convex between $(\rho,m)$ and its best response $(\hat{\rho},\hat{m})$, then the distance between $\rho$ and its best response is equivalent to the gain of the best response at $\rho$, i.e.
    \begin{equation}
        \frac{\lambda_R}{2}\|\hat{\rho}-\rho\|^2_{Q_T}\leq g(\rho,m).
    \label{eq: lem best response 2}
    \end{equation}
    
\end{lemma}

\begin{proof}
    By definition, we have
    \begin{equation}
        J(\tilde{\rho},\tilde{m};\rho)-J(\hat{\rho},\hat{m};\rho)
        =R(\tilde{\rho},\tilde{m}) - R(\hat{\rho},\hat{m})
        +\langle f(\rho),\tilde{\rho}-\hat{\rho}\rangle_{Q_T}
        +\langle f_T(\rho_T),\tilde{\rho}_T-\hat{\rho}_T\rangle_{\bbT^d}
    \end{equation}
    By Lemma \ref{lem: ficplay welldefined}, $\hat{\rho}>0$. Therefore, the directional strong convexity of $R$ gives that
    \begin{equation}
        R(\tilde{\rho},\tilde{m}) - R(\hat{\rho},\hat{m})
        \geq \left\langle L(x,\frac{\hat{m}}{\hat{\rho}})
                    -\frac{\hat{m}}{\hat{\rho}}D_vL(x,\frac{\hat{m}}{\hat{\rho}}),
            \tilde{\rho}-\hat{\rho} \right\rangle_{Q_T}
            +\left\langle D_vL(x,\frac{\hat{m}}{\hat{\rho}}),
            \tilde{m}-\hat{m} \right\rangle_{Q_T}
            +\frac{\lambda_R}{2}\|\tilde{\rho}-\hat{\rho}\|_{Q_T}^2.
    \label{eq: lem best response cvx}
    \end{equation}
    By optimality of $(\hat{\rho},\hat{m})$ and according to Lemma \ref{lem: prepare}, there exists $\hat{\phi}$ such that
    \begin{equation}
        L(x,\frac{\hat{m}}{\hat{\rho}})
        -\frac{\hat{m}}{\hat{\rho}}D_vL(x,\frac{\hat{m}}{\hat{\rho}})
        +f(x,\rho)=-\partial_t\hat{\phi}-\nu\Delta\hat{\phi},\quad
        D_vL(x,\frac{\hat{m}}{\hat{\rho}})=-\nabla\hat{\phi}.
    \end{equation}
    Combining the above gives us
    \begin{equation}
        J(\tilde{\rho},\tilde{m};\rho)-J(\hat{\rho},\hat{m};\rho)
        \geq A
            +\frac{\lambda_R}{2}\|\tilde{\rho}-\hat{\rho}\|_{Q_T}^2
    \end{equation}
    where
    \begin{equation}
        A=
        \langle -\partial_t\hat{\phi}-\nu\Delta\hat{\phi},
                    \tilde{\rho}-\hat{\rho}\rangle_{Q_T}
            +\langle -\nabla\hat{\phi},
                    \tilde{m}-\hat{m} \rangle_{Q_T}
            +\langle f_T(\rho_T),\tilde{\rho}_T-\hat{\rho}_T\rangle_{\bbT^d}
    \end{equation}
    Because $(\tilde{\rho},\tilde{m}),(\hat{\rho},\hat{m})\in C$, and $\hat{\phi}(x,T)=f_T(x,\rho_T)$, with integration by part, we have
    \begin{equation}
    \begin{aligned}
        A=&\langle \hat{\phi}, \partial_t(\tilde{\rho}-\hat{\rho}) - \nu\Delta(\tilde{\rho}-\hat{\rho}) + \nabla\cdot(\tilde{m}-\hat{m}) \rangle_{Q_T}
        =0.
    \end{aligned}
    \end{equation}
    Therefore \eqref{eq: lem best response 1} is true.
    Taking $(\tilde{\rho},\tilde{m})=(\rho,m)$ and with the definition of gain of the best response in \eqref{eq: def exploit}, we conclude \eqref{eq: lem best response 2}.
    
\end{proof}

The next lemma shows that if the difference between $\rho$ and its best response $\hat{\rho}$ is small, then $\rho$ is close to the Nash Equilibrium $\rho^*$.
\begin{lemma}
\label{lem: NE < BR}
    For given $(\rho,m)\in C$, let its best response be $(\hat{\rho},\hat{m})$ and let $(\rho^*,m^*,\phi^*)$ be a Nash Equilibrium.
    Under Assumptions \ref{aspt: lip cts}-\ref{aspt: independent cost}, if in addition $f$ is $\lambda_f$-strongly monotone, then 
    \begin{equation}
        \|\rho-\rho^*\|_{Q_T}\leq \frac{L_f}{\lambda_f}\|\rho-\hat{\rho}\|_{Q_T}.
    \label{eq: lem NE < BR}
    \end{equation}
\end{lemma}
\begin{proof}
    Notice that
    \begin{equation}
        \langle \rho-\rho^*,f(\rho)-f(\rho^*) \rangle_{Q_T}
        = \langle \rho-\hat{\rho}, f(\rho)-f(\rho^*) \rangle_{Q_T}
        + \langle \hat{\rho}-\rho^*, f(\rho)-f(\rho^*) \rangle_{Q_T},
    \end{equation}
    Since $(\rho^*,m^*,\phi^*)$ is a Nash Equilibrium, by Lemma \ref{lem: NE - dual gap=0}, we know that $(\rho^*,m^*)=\argmin_{(\tilde{\rho},\tilde{m})\in C}J(\tilde{\rho},\tilde{m};\rho^*)$.
    Applying Proposition \ref{lem: stability} then gives
    \begin{equation}
        \langle \rho-\rho^*,f(\rho)-f(\rho^*) \rangle_{Q_T}
        \leq  \langle \rho-\hat{\rho}, f(\rho)-f(\rho^*) \rangle_{Q_T}.
    \end{equation}    
    Therefore by strong monotonicity and Lipschitz continuity of $f$, we have
    \begin{equation}
        \lambda_f\|\rho-\rho^*\|_{Q_T}\leq L_f\|\rho-\hat{\rho}\|_{Q_T},
    \end{equation} 
    which concludes \eqref{eq: lem NE < BR}.
\end{proof}

In the end, we show the stability of the best response mapping and as a consequence, if $\rho$ is close to $\rho^*$, then $\rho$ is close to its best response.
\begin{lemma}
\label{lem: optctl map is cts}
    Let $(\hat{\rho}^i,\hat{m}^i)$ be the best response to given $(\rho^i,m^i),i=0,1$.
    Under the Assumption \ref{aspt: independent cost}, if $R$ is $\lambda_R$-strongly convex between $(\hat{\rho}^0,\hat{m}^0)$ and $(\hat{\rho}^1,\hat{m}^1)$,
    then 
    \begin{equation}
        \langle \hat{\rho}^0-\hat{\rho}^1, f(\rho^0)-f(\rho^1) \rangle_{Q_T} + D = 0 ,
    \label{eq: lem stability}
    \end{equation}
    where 
    \begin{equation}
        D \geq \lambda_R \|\hat{\rho}^0-\hat{\rho}^1\|_{Q_T}^2.
    \label{eq: lem est D}
    \end{equation}
    As a consequence
    \begin{equation}
        \|\hat{\rho}^0-\hat{\rho}^1\|_{Q_T}  
    \leq \frac{1}{\lambda_R}\|f({\rho}^0)-f({\rho}^1)\|_{Q_T}.
    \label{eq: lem optctl map is cts}
    \end{equation}
\end{lemma}
\begin{proof}
    Denote $(\hat{\rho}^\theta,\hat{m}^\theta):=
    (1-\theta)(\hat{\rho}^0,\hat{m}^0) + \theta(\hat{\rho}^1,\hat{m}^1)$.
    Then the optimality of $(\hat{\rho}^0,\hat{m}^0)$ gives
    \begin{equation}
        R(\hat{\rho}^0,\hat{m}^0)
        +\langle f(\rho^0),\hat{\rho}^0 \rangle_{Q_T} 
        + \langle f_T,\hat{\rho}^0_T \rangle_{\bbT^d} 
        \leq R(\hat{\rho}^\theta,\hat{m}^\theta)
        +\langle f(\rho^0),\hat{\rho}^\theta \rangle_{Q_T} 
        + \langle f_T,\hat{\rho}^\theta_T \rangle_{\bbT^d},
    \end{equation}
    which is equivalent to
    \begin{equation}
        R(\hat{\rho}^0,\hat{m}^0) 
        \leq R(\hat{\rho}^\theta,\hat{m}^\theta)
            +\theta \langle f(\rho^0),\hat{\rho}^1-\hat{\rho}^0 \rangle_{Q_T} 
            +\theta \langle f_T,\hat{\rho}^1_T-\hat{\rho}^0_T \rangle_{\bbT^d}.
    \label{eq: lem cts - opt}
    \end{equation}
    By the directional strong convexity of $R$, we have
    \begin{equation}
        R(\hat{\rho}^\theta,\hat{m}^\theta) 
            + \frac{\theta(1-\theta)\lambda_R}{2}\|\hat{\rho}^0-\hat{\rho}^1\|^2_{Q_T}
        \leq (1-\theta)R(\hat{\rho}^0,\hat{m}^0) 
            + \theta R(\hat{\rho}^1,\hat{m}^1). 
    \label{eq: lem cts - cvx}
    \end{equation}
    Adding both sides of \eqref{eq: lem cts - opt} and \eqref{eq: lem cts - cvx} gives
    \begin{equation}
        \frac{(1-\theta)\lambda_R}{2}\|\hat{\rho}^0-\hat{\rho}^1\|^2_{Q_T}
        \leq -R(\hat{\rho}^0,\hat{m}^0) + R(\hat{\rho}^1,\hat{m}^1) 
            +\langle f(\rho^0),\hat{\rho}^1-\hat{\rho}^0 \rangle_{Q_T} 
            +\langle f_T,\hat{\rho}^1_T-\hat{\rho}^0_T \rangle_{\bbT^d}.
    \end{equation}
    Then we let $\theta\to 0$ and obtain
    \begin{equation}
        \frac{\lambda_R}{2}\|\hat{\rho}^0-\hat{\rho}^1\|^2_{Q_T}
        \leq -R(\hat{\rho}^0,\hat{m}^0) + R(\hat{\rho}^1,\hat{m}^1) 
            +\langle f(\rho^0),\hat{\rho}^1-\hat{\rho}^0 \rangle_{Q_T} 
            +\langle f_T,\hat{\rho}^1_T-\hat{\rho}^0_T \rangle_{\bbT^d}.
    \label{eq: lem cts - 0}
    \end{equation}
    Similarly, the optimality of $(\hat{\rho}^1,\hat{m}^1)$ gives
    \begin{equation}
        \frac{\lambda_R}{2}\|\hat{\rho}^0-\hat{\rho}^1\|^2_{Q_T}
        \leq -r(\hat{\rho}^1,\hat{m}^1) + r(\hat{\rho}^0,\hat{m}^0) 
            +\langle f(\rho^1),\hat{\rho}^0-\hat{\rho}^1 \rangle_{Q_T} 
            +\langle f_T,\hat{\rho}^0_T-\hat{\rho}^1_T \rangle_{\bbT^d}.
    \label{eq: lem cts - 1}
    \end{equation}
    Combining \eqref{eq: lem cts - 0} and \eqref{eq: lem cts - 1}, we have
    \begin{equation}
        \lambda_R\|\hat{\rho}^0-\hat{\rho}^1\|^2_{Q_T}
        \leq -\langle f(\rho^0)-f(\rho^1),\hat{\rho}^0-\hat{\rho}^1 \rangle_{Q_T}
        \leq \|f(\rho^1)-f(\rho^0)\|_{Q_T}\|\hat{\rho}^1-\hat{\rho}^0\|_{Q_T}, 
    \end{equation}
    which concludes the proof.   
\end{proof}

\begin{corollary}
    For given $(\rho,m)\in C$, let its best response be $(\hat{\rho},\hat{m})$ and let $(\rho^*,m^*,\phi^*)$ be a Nash Equilibrium.
    Under the Assumption \ref{aspt: lip cts}-\ref{aspt: independent cost}, if $R$ is $\lambda_R$-strongly convex between $({\rho},{m})$ and $({\rho}^*,{m}^*)$,
    then
    \begin{equation}
        \|\rho-\hat{\rho}\|_{Q_T}\leq \frac{L_f+\lambda_R}{\lambda_R}\|\rho-\rho^*\|_{Q_T}.
    \end{equation}
\end{corollary}

\begin{proof}
    By Assumption \ref{aspt: lip cts}, $\|f(\rho^1)-f(\rho^0)\|_{Q_T}\leq L_f\|\rho^1-\rho^0\|_{Q_T}$.
    Taking $\rho^0=\rho$ and $\rho^1=\rho^*$ in Lemma \ref{lem: optctl map is cts} and by triangle inequality concludes the proof.
\end{proof}

\subsection{Non-potential MFGs}
\label{subsec: non pot mfg}

In this section, we adopt the gain of the best response $g(\rho,m)$ as the Lyapunov function to establish the convergence of scheme \eqref{eq: ficplay lin}-\eqref{eq: ficplay avg rhom} and thus fictitious play (Algorithm \ref{alg: ficplay}) for non-potential MFGs. 

Under certain conditions, Lemma \ref{lem: best response dist bounded by exploit} implies that,
\begin{equation}
    \frac{\lambda_R}{2}\|\rho-\hat{\rho}\|_{Q_T}^2 
\leq g(\rho,m),
\end{equation}
and combining Lemma \ref{lem: optctl map is cts} and Lemma \ref{lem: NE < BR} gives 
\begin{equation}
    \frac{\lambda_f}{L_f} \|\rho-\rho^*\|_{Q_T}
\leq \|\rho-\hat{\rho}\|_{Q_T} \leq
\frac{L_f+\lambda_R}{\lambda_R}\|\rho-\rho^*\|_{Q_T}.
\end{equation}
Therefore, any of $\|\rho-\hat{\rho}\|_{Q_T}$, $\|\rho-\rho^*\|_{Q_T}$ and $g(\rho,m)$ converges to 0 implies $\rho$ converges to the Nash Equilibrium $\rho^*$.

While it is not straightforward to prove that $\|\rho-\hat{\rho}\|_{Q_T}\to0$ or $\|\rho-\rho^*\|_{Q_T}\to0$ for fictitious play, Lemma \ref{lem: dual gap recur} in this section shows that the gain of the best response $g$ decays along the iteration.
To be precise, we show that
\begin{equation}
    g(\rho^{(k)},m^{(k)})-g(\rho^{(k-1)},m^{(k-1)})\leq-\delta_kg(\rho^{(k-1)},m^{(k-1)}) + z_k.
\end{equation}
According to Lemma \ref{lem: optctl map is cts}, the residue $z_k=c_1\delta_k^2\|\hat{\rho}^{(k)}-\rho^{(k-1)}\|^2_{Q_T}$ for some constant $c_1$. Therefore we can control the sign of the right-hand side by the weight $\delta_k$ and obtain the decay and convergence of $g$.
In addition, Lemma \ref{lem: best response dist bounded by exploit} shows that $\|\hat{\rho}^{(k)}-\rho^{(k-1)}\|^2_{Q_T}\leq c_2g(\rho^{(k-1)},m^{(k-1)})$ for some constant $c_2$.
Therefore $z_k\leq c_1c_2\delta_k^2 g(\rho^{(k-1)},m^{(k-1)})$ and
\begin{equation}
    g(\rho^{(k)},m^{(k)})\leq(1-\delta_k+c_1c_2\delta_k^2)g(\rho^{(k-1)},m^{(k-1)}),
\end{equation}
which gives linear convergence.

We first prove the decay property of $g$. 
A similar technique was used in \cite{hadikhanloo2017anonymous} for anonymous games.
\begin{lemma}[The decay property of $g$]
\label{lem: dual gap recur}
    Let the sequence $\{({\rho}^{(k)},{m}^{(k)})\}, \{(\hat{\rho}^{(k)},\hat{m}^{(k)})\}$ be generated by the numerical scheme \eqref{eq: ficplay lin}-\eqref{eq: ficplay avg rhom}.
    If Assumption \ref{aspt: independent cost} holds and $f$ is monotone, 
    then 
    \begin{equation}
        g(\rho^{(k)},m^{(k)})
    \leq (1-\delta_k)g(\rho^{(k-1)},m^{(k-1)}) + D,\text{ for } k=1, 2, \cdots.
    \end{equation}
    where
    $$ D= -\langle f(\rho^{(k-1)})-f(\rho^{(k)}), \hat{\rho}^{(k)} - \hat{\rho}^{(k+1)}\rangle_{Q_T} \geq 0. $$
    If in addition, Assumptions \ref{aspt: lip cts} holds,
    $R$ is directional $\lambda_R$-strongly convex between $({\rho}^{(k-1)},{m}^{(k-1)})$ and $(\hat{\rho}^{(k)},\hat{m}^{(k)})$, and between $(\hat{\rho}^{(k)},\hat{m}^{(k)})$ and $(\hat{\rho}^{(k+1)},\hat{m}^{(k+1)})$,  
    then
    \begin{equation}
        D \leq \delta_k^2\frac{L_f^2}{\lambda_R}\|\hat{\rho}^k-\rho^{(k-1)}\|_{Q_T}^2 \leq \delta_k^2\frac{2L_f^2}{\lambda_R^2}g(\rho^{(k-1)},m^{(k-1)}).
    \end{equation}
\end{lemma}

\begin{proof}
The definition of the best response gives
\begin{equation}
    J(\hat{\rho}^{(k)},\hat{m}^{(k)};
           \rho^{(k-1)})
        \leq J(\hat{\rho}^{(k+1)},\hat{m}^{(k+1)};
           \rho^{(k-1)}).
\end{equation}
Therefore, we have
\begin{equation}
\begin{aligned}
    &g(\rho^{(k)},m^{(k)}) 
        - g(\rho^{(k-1)},m^{(k-1)})\\
    \leq & J(\rho^{(k)},m^{(k)}; 
        \rho^{(k)})
        -J(\rho^{(k-1)},m^{(k-1)};
           \rho^{(k-1)})\\
    &+J(\hat{\rho}^{(k+1)},\hat{m}^{(k+1)};
           \rho^{(k-1)})
        -J(\hat{\rho}^{(k+1)},\hat{m}^{(k+1)};
           \rho^{(k)})\\
    =&\underbrace{J(\rho^{(k)},m^{(k)}; 
             \rho^{(k)})
        -J(\rho^{(k-1)},m^{(k-1)};
           \rho^{(k-1)})
    + \langle f(\rho^{(k-1)})-f(\rho^{(k)}), \hat{\rho}^{(k)}\rangle_{Q_T}}_A+D,\\
\end{aligned}
\end{equation}
where $D=-\langle f(\rho^{(k-1)})-f(\rho^{(k)}), \hat{\rho}^{(k)}-\hat{\rho}^{(k+1)}\rangle_{Q_T}\geq 0$, by Proposition \ref{lem: stability}. And by Lemma \ref{lem: optctl map is cts}, if Assumption \ref{aspt: lip cts} holds and $R$ is directional convex, then 
\begin{equation}
    D\leq \frac{1}{\lambda_R}\|f(\rho^{(k-1)})-f(\rho^{(k)})\|_{Q_T}^2
    \leq \delta_k^2\frac{L_f^2}{\lambda_R}\|\hat{\rho}^{(k)}-\rho^{(k-1)}\|_{Q_T}^2 \leq \delta_k^2\frac{2L_f^2}{\lambda_R^2}g(\rho^{(k-1)},m^{(k-1)}).
\end{equation}
It suffices to show that $A\leq-\delta_k g(\rho^{(k-1)},m^{(k-1)})$.
By definition of $J$, we insert $\pm\langle f(\rho^{(k-1)}),\rho^{(k)} \rangle_{Q_T}$ to $A$ and obtain
\begin{equation*}
\begin{aligned}
    A = & \underbrace{R(\rho^{(k)},m^{(k)}) - R(\rho^{(k-1)},m^{(k-1)})
        +\langle f(\rho^{(k-1)}),{\rho}^{(k)}-\rho^{(k-1)} \rangle_{Q_T}
        +\langle f_T, {\rho}^{(k)}_T-\rho^{(k-1)}_T \rangle_{\bbT^d}}_{A_1}\\
    & +\underbrace{\langle f(\rho^{(k)})-f(\rho^{(k-1)}),\rho^{(k)} \rangle_{Q_T}
    +\langle f(\rho^{(k-1)})-f(\rho^{(k)}), \hat{\rho}^{(k)}\rangle_{Q_T}
    }_{A_2}\\
\end{aligned}
\end{equation*}
Notice that $(\rho^{(k)},m^{(k)})$ is a convex combination of $(\rho^{(k-1)},m^{(k-1)})$ and $(\hat{\rho}^{(k)},\hat{m}^{(k)})$. The convexity of $R$ and definition of $g$ thus gives
\begin{equation}
\begin{aligned}
    A_1\leq& \delta_k (R(\hat{\rho}^{(k)},\hat{m}^{(k)}) - R(\rho^{(k-1)},m^{(k-1)}))
        +\delta_k \langle f(\rho^{(k-1)}),\hat{\rho}^{(k)}-\rho^{(k-1)} \rangle_{Q_T}
        +\delta_k \langle f_T, \hat{\rho}^{(k)}_T-\rho^{(k-1)}_T \rangle_{\bbT^d}\\
    = & -\delta_k g(\rho^{(k-1)},m^{(k-1)}).    
\end{aligned}
\end{equation}
By writing $\hat{\rho}^{(k)}$ as a linear combination of $\rho^{(k-1)}$ and $\rho^{(k)}$ and applying the monotonicity of $f$, we have
\begin{equation}
    A_2=-\frac{(1-\delta_k)}{\delta_k}\langle f(\rho^{(k)}) - f(\rho^{(k-1)}),\rho^{(k)}-\rho^{(k-1)} \rangle_{Q_T}\leq 0.  
\label{eq: monotone term}
\end{equation}
Combining $A_1$ and $A_2$ proves the lemma.
\end{proof}

Before stating the main theorem, we remark on the cases when $f_T$ depends on $\rho_T$ or when $f$ is not monotone.
\begin{remark}[Lack of Assumption \ref{aspt: independent cost}]
\label{rem: nonind cost}
    If Assumption \ref{aspt: independent cost} does not hold, i.e. the terminal cost also depends on $\rho$, but $f$ and $f_T$ are monotone, then under sufficient regularity and convexity assumptions, we have
    \begin{equation}
        g(\rho^{(k)},m^{(k)})-g(\rho^{(k-1)},m^{(k-1)})
        \leq -\delta_k g(\rho^{(k-1)},m^{(k-1)}) + D,
    \end{equation}
    where 
    $$D=-\langle f(\rho^{(k-1)})-f(\rho^{(k)}), \hat{\rho}^{(k)} - \hat{\rho}^{(k+1)}\rangle_{Q_T}
    -\langle f_T(\rho^{(k-1)})-f_T(\rho^{(k)}), \hat{\rho}_T^{(k)} - \hat{\rho}_T^{(k+1)}\rangle_{\bbT^d}.$$
    As long as $D < \delta_k g(\rho^{(k-1)},m^{(k-1)})$, we still have the decay of $g$.
\end{remark}

\begin{remark}[Non-monotone cost]
\label{rem: nonmono cost}
    If $f$ is weakly monotone, i.e. there exists some $\mu_f\geq 0$ such that for any $\rho,\rho'\in\calP(\bbT^d)\cap L^2(\bbT^d)$,
    \begin{equation*}
        \int_{\bbT^d}(f(x,\rho)-f(x,\rho'))\dd (\rho-\rho')(x)\geq -\mu_f \|\rho-\rho'\|_{\bbT^d}^2,
    \end{equation*}  
    then the term $A_2$ in \eqref{eq: monotone term} can be bounded by
    \begin{equation}
        A_2\leq \mu_f\frac{1-\delta_k}{\delta_k}\|\rho^{(k)}-\rho^{(k-1)}\|_{Q_T}^2
        = \mu_f\delta_k(1-\delta_k)\|\hat{\rho}^{(k)}-\rho^{(k-1)}\|_{Q_T}^2
        \leq\mu_f\delta_k\|\hat{\rho}^{(k)}-\rho^{(k-1)}\|_{Q_T}^2.
    \end{equation}
    Under the same regularity and convexity assumptions of Lemma \ref{lem: dual gap recur},
    \begin{equation}
        A_2\leq \frac{2\mu_f}{\lambda_R}\delta_k g(\rho^{(k-1)},m^{(k-1)}),
    \end{equation}
    and the recurrent relation of the gain function becomes 
    \begin{equation}
        g(\rho^{(k)},m^{(k)})-g(\rho^{(k-1)},m^{(k-1)})
        \leq \left(-\delta_k\left(1-\frac{2\mu_f}{\lambda_R}\right)
        +\delta_k^2\frac{2L_f^2}{\lambda_R^2}\right) g(\rho^{(k-1)},m^{(k-1)}).
    \end{equation}
    As long as $\mu_f\leq\frac{\lambda_R}{2}$, one can choose the weight $\delta_k$ such that $-\delta_k\left(1-\frac{2\mu_f}{\lambda_R}\right)
        +\delta_k^2\frac{2L_f^2}{\lambda_R^2}<0$ and thus $g$ decays.
\end{remark}

With Lemma \ref{lem: dual gap recur}, we state the convergence of fictitious play for general MFGs possibly without a potential structure as follows.
\begin{theorem}
\label{thm: non-pot cvg rate}
    Let the sequence $\{({\rho}^{(k)},{m}^{(k)})\}, \{(\hat{\rho}^{(k)},\hat{m}^{(k)})\}$ be generated by the numerical scheme \eqref{eq: ficplay lin}-\eqref{eq: ficplay avg rhom}.
    Under Assumption \ref{aspt: lip cts}-\ref{aspt: independent cost}, if in addition $f$ is monotone and $R$ is directional $\lambda_R$-strongly convex between $(\hat{\rho}^{(k)},\hat{m}^{(k)})$ and $(\hat{\rho}^{(k+1)},\hat{m}^{(k+1)})$, and between $({\rho}^{(k)},{m}^{(k)})$ and $(\hat{\rho}^{(k+1)},\hat{m}^{(k+1)})$, for any $k$,
    then we have
    \begin{equation}
        g(\rho^{(K)},m^{(K)})
        \leq \exp\left(\sum_{k=1}^K (c\delta_{k}^{2} -\delta_k)\right)
        g(\rho^{(0)},m^{(0)}),
    \label{eq: thm non-pot general}
    \end{equation}
    where $c=\frac{2L_f^2}{\lambda_R^2}$.
    Specifically, taking weight $\delta_k=\frac{\alpha}{k+\alpha}$ with a fixed $\alpha>0$ gives sublinear convergence
    \begin{equation}
        g(\rho^{(K)},m^{(K)})
        \leq g(\rho^{(0)},m^{(0)})
        \frac{(\alpha+1)^\alpha\exp(c\alpha)}{(K+\alpha+1)^\alpha}, 
    \label{eq: thm non-pot sublinear}
    \end{equation}
    and taking weight $\delta_k=\delta\leq\min\{1,\frac{1}{c}\}$ gives linear convergence
    \begin{equation}
        g({\rho}^{(K)},{m}^{(K)})
        \leq \exp\left( -\delta(1-c\delta)K \right) g({\rho}^{(0)},{m}^{(0)}).
    \label{eq: thm non-pot linear}
    \end{equation}
\end{theorem}

\begin{proof}
    Applying Lemma \ref{lem: optctl map is cts} and Lemma \ref{lem: best response dist bounded by exploit} to Lemma \ref{lem: dual gap recur}, and combining the Lipschitz continuity of $f$ give
    \begin{equation}
        g(\rho^{(k)},m^{(k)})
        \leq     
        \left(1-\delta_k + \frac{2L_f^2}{\lambda_R^2}\delta_k^2
        \right)g(\rho^{(k-1)},m^{(k-1)})
    \label{eq: thm non-pot recur}
    \end{equation}
    Using the fact $1+x \le e^{x}$ for all $x\in \mathbb R$, we have
    \begin{equation}
        g(\rho^{(k)},m^{(k)})
        \leq     
        \exp(c\delta_k^2-\delta_k)g(\rho^{(k-1)},m^{(k-1)}),
    \label{eq: 75}
    \end{equation}
    Telescoping from $k=1$ to $K$ gives \eqref{eq: thm non-pot general}.

    Taking $\delta_k=\frac{\alpha}{k+\alpha},\alpha>0$ gives 
    $$ \sum_{k=1}^K \delta_k \geq \int_1^{K+1}\frac{\alpha}{x+\alpha}\dd x = \alpha \ln\frac{K+\alpha+1}{\alpha+1}, $$
    $$ \sum_{k=1}^K \delta_k^2 \leq \alpha^2  \int_0^\infty \frac{1}{(x+\alpha)^2}\dd x =\alpha, $$
    and therefore \eqref{eq: thm non-pot sublinear}.

    Taking $\delta_k=\delta$ in \eqref{eq: 75} and telescoping gives \eqref{eq: thm non-pot linear}.
\end{proof}

\begin{corollary}
\label{cor: non-pot cvg rate}
    Let $(\rho^*,m^*)$ be the Nash Equilibrium
    and the sequence $\{({\rho}^{(k)},{m}^{(k)})\}, \{(\hat{\rho}^{(k)},\hat{m}^{(k)})\}$ be generated by the numerical scheme \eqref{eq: ficplay lin}-\eqref{eq: ficplay avg rhom}.
    Under Assumptions \ref{aspt: lip cts}-\ref{aspt: independent cost}, and assumptions of Theorem \ref{thm: non-pot cvg rate}, if in addition $f$ is $\lambda_f$ monotone, then we have
    \begin{equation}
        \|\rho^{(K)}-\rho^*\|_{Q_T}
        \leq \frac{L_f}{\lambda_f}\sqrt{\frac{2}{\lambda_R}g(\rho^0,m^0)}
        \exp\left(\half\sum_{k=1}^K (c\delta_{k}^{2} -\delta_k)\right),
    \label{eq: cor non-pot general}
    \end{equation}
    where $c=\frac{2L_f^2}{\lambda_R^2}$.
\end{corollary}
The proof is straightforward from Theorem \ref{thm: non-pot cvg rate}, Lemma \ref{lem: best response dist bounded by exploit} and Lemma \ref{lem: NE < BR}.

\subsection{Acceleration: back-tracking line search}
\label{subsec: acc btls}

Theorem \ref{thm: non-pot cvg rate} demonstrates that the weight $\delta_k$ fundamentally determines the convergence rate. 
Recall that the key to the convergence of fictitious play is Lemma \ref{lem: dual gap recur}, which states
$$g(\rho^{(k)},m^{(k)}) - g(\rho^{(k-1)},m^{(k-1)}) \leq -\delta_k g(\rho^{(k-1)},m^{(k-1)}) + D,$$
and
$$ D = -\langle f(\rho^{(k-1)}) - f(\rho^{(k)}), \hat{\rho}^{(k)} - \hat{\rho}^{(k+1)} \rangle_{Q_T} \geq 0.$$
To achieve convergence, we require the improvement in the gain function to be positive
$$\delta_k g(\rho^{(k-1)},m^{(k-1)}) - D > 0,$$
and to ensure fast convergence, we need to maximize this improvement.
When the best-response mapping is continuous, $D = \mathcal{O}(\delta_k^2)$, and a sufficiently small $\delta_k$ guarantees positive improvement. Thus, decaying weights such as $\delta_k = \frac{1}{k+1}$ in \cite{cardaliaguet2017ficplay} or $\delta_k = \frac{2}{k+2}$ will eventually lead to convergence. However, since the improvement in the gain function is $\mathcal{O}(\delta_k)$, with these choices of $\delta_k$, the convergence slows down as $k$ increases.

Therefore, it is essential to choose the largest feasible $\delta_k$ that still ensures a positive improvement.
In this section, we adopt a backtracking line search strategy from optimization, allowing for a dynamic selection of an optimal weight $\delta_k$ to enhance the convergence rate. The numerical experiment in Section \ref{subsec: num btls} illustrates the effectiveness of this acceleration.

Algorithm \ref{alg: btls} outlines the process, where $a \leftarrow b$ denotes assigning the value of $b$ to $a$. Each iteration within the while-loop requires solving the HJB and FP equations once.
The backtracking line search aims to find the smallest non-negative integer $n$ such that $\delta_k = \delta \beta^n$ and $D \leq \zeta \delta_k g(\rho^{(k-1)}, m^{(k-1)})$, where $\delta \in (0,1]$ and $\beta, \zeta \in (0,1)$ are preselected parameters.
The parameter $\zeta$ controls the amount of gain improvement, specifically $(1 - \zeta) \delta_k g(\rho^{(k-1)}, m^{(k-1)})$. A smaller $\zeta$ yields greater improvement but may require more iterations to satisfy or may be difficult to achieve. The parameter $\beta$ determines the shrinkage rate for $\delta_k$ when the target improvement is unmet. A smaller $\beta$ results in fewer iterations to find a valid $\delta_k$ but may lead to overly modest improvements.

\begin{algorithm}[ht]
\caption{Back-tracking Line Search Oracle}
\begin{algorithmic}
    \STATE{Given} $(\rho^{(k-1)},m^{(k-1)})$, its best response $(\hat{\rho}^{(k)},m^{(k)})$ and consequently the gain $g(\rho^{(k-1)},m^{(k-1)})$.
    \STATE{Parameters} $0< \beta,\zeta < 1.$
    \STATE{Initialization} $\delta_k\leftarrow\delta, 0< \delta \leq 1.$
    \WHILE{$D > \zeta\delta_k g(\rho^{(k-1)},m^{(k-1)})$}
    
        \STATE{$\delta_k\leftarrow \beta\delta_k$}

        \STATE{$(\rho^{(k)},m^{(k)})\leftarrow (1-\delta)(\rho^{(k)},m^{(k)}) + \delta (\hat{\rho}^{(k)},m^{(k)})$}

        \STATE{Compute $(\hat{\rho}^{(k+1)},m^{(k+1)})$ (the best response of $(\rho^{(k)},m^{(k)})$) and $D$.}
        
    \ENDWHILE
    \STATE{Output} $\delta_k.$
\end{algorithmic}
\label{alg: btls}
\end{algorithm}

We close this section with a few remarks on the weights.

\begin{remark}[Choices of weights]
    According to Theorem \ref{thm: non-pot cvg rate}, the weight (step-size) $\delta_k=\frac{\alpha}{k+\alpha},(\alpha>0)$ gives a convergence rate $\frac{1}{k^{\alpha}}$.     
    But this does not suggest that a larger $\alpha$ is in practice better because the constant also grows exponentially with $\alpha$.
    This is also noticed and discussed in the conditional gradient method \cite{nesterov2018gcg-primaldual} and in the potential games setting \cite{Lavigne2023gcg-ficplay}.
\end{remark}

\begin{remark}[Warm restart is not the optimal strategy]
\label{rem: restart}
    \cite{hadikhanloo2019finite} observes that fictitious play with $\delta_k = \frac{1}{k}$ converges very slowly as the number of iterations grows, primarily because excessive weight is given to inaccurate initial guesses and historical data. 
    Without a theoretical explanation, \cite{hadikhanloo2019finite} suggests and empirically demonstrates that when $({\rho}^{(k)}, {m}^{(k)})$ is near equilibrium, restarting fictitious play with $({\rho}^{(k)}, {m}^{(k)})$ keeps the value of $k$ small, thereby accelerating convergence. 
    Our analysis reveals that this slow convergence for large $k$ stems from the sublinear convergence rate. 
    Restarting the algorithm with a small $k$ results in a large $\delta_k$, which could prevent the decay property of the gain function $g$ from holding and lead to a large increase of $g$. Thus, the restart strategy primarily aids convergence asymptotically as $k$ increases.
\end{remark}

\begin{remark}[Back-tracking line search for potential MFGs]
    \cite{Lavigne2023gcg-ficplay} also proposes a back-tracking line search strategy to determine the weight $\delta_k$. However, their strategy requires the existence of the potential and is only applicable to potential MFGs.
\end{remark}

\subsection{Potential MFGs}
\label{subsec: pot mfg}

Finally, we turn our attention to potential MFGs. We begin by examining the relationship between our Lyapunov function $g$ and the one used in \cite{Lavigne2023gcg-ficplay} in the context of potential MFGs. Next, we present a simplified convergence proof tailored to the potential MFGs and highlight the key differences from the general, non-potential case. We conclude by comparing fictitious play with other optimization-based algorithms designed for potential MFGs.

Throughout Section \ref{subsec: pot mfg}, we assume that Assumptions \ref{aspt: lip cts} and \ref{aspt: independent cost} hold. In addition, we introduce the following assumption:
\begin{assumption}
$f$ is derived from a functional $F : \calP(\bbT^d) \to \bbR$. 
\label{asp: f is potential}    
\end{assumption}
Assumption \ref{aspt: independent cost} further implies that $f_T$ is also derived from a potential of linear form
\begin{equation*}
    F_T : \calP(\bbT^d) \to \bbR, \quad \rho \mapsto \int_{\bbT^d} f_T(x) \, \dd\rho(x).
\end{equation*}
Under these conditions, and by Proposition \ref{thm: potential mfg}, the fixed-point formulation of the MFG system \eqref{eq: mfg} can be equivalently expressed as the convex optimization problem given in the variational formulation \eqref{eq: pot mfg}.

\cite{Lavigne2023gcg-ficplay} points out that the fictitious play \eqref{eq: ficplay lin}-\eqref{eq: ficplay avg rhom} is equivalent to applying the generalized conditional gradient (GCG) to the variational formulation \eqref{eq: pot mfg}. 
GCG \cite{Bredies2009gcg-itershrink, Kunisch2024gcg-fastcvg} is a projection-free optimization algorithm and can be viewed as a variant of the Frank-Wolfe method \cite{jaggi13gcg-revisiting-cg}.
Consider the optimization problem
\begin{equation}
    \min_{z\in Z} \quad j(z) := r(z) + s(z)
\end{equation}
where $Z \subset \bbR^d$ is a convex set, $r$ is possibly non-differentiable, and $s$ is differentiable. At iteration $k$, the GCG algorithm linearizes $s$ at $z^{(k-1)}$ and solves
\begin{equation}
    \hat{z}^{(k)} = \argmin_{\tilde{z}\in Z} ~ r(\tilde{z}) + \langle\nabla s(z^{(k-1)}),\tilde{z}\rangle,
\label{eq: gcg lin}
\end{equation}
followed by the update
\begin{equation}
    z^{(k)} = (1-\delta_k)z^{(k-1)} + \delta_k \hat{z}^{(k)}.
\label{eq: gcg avg}
\end{equation}
Assuming $z^* \in Z$ is a minimizer, the primal gap is defined as
\begin{equation}
    h(z) := j(z) - j(z^*)
\end{equation}
and the dual gap is defined as
\begin{equation}
    g(z) := \max_{\tilde{z}\in C} \{r(z) - r(\tilde{z}) + \langle\nabla s(z),z-\tilde{z}\rangle\} 
    = r(z) + \langle \nabla s(z),z \rangle - \min_{\tilde{z}\in C} \{r(\tilde{z}) + \langle \nabla s(z), \tilde{z} \rangle \}
\end{equation}

In the MFG setting, at each iteration $(\rho^{(k-1)}, m^{(k-1)})$, solving for the best response \eqref{eq: ficplay lin} corresponds to \eqref{eq: gcg lin} by
linearizing the interaction and terminal costs. 
And \eqref{eq: ficplay avg rhom} corresponds to \eqref{eq: gcg avg}.
Assume $(\rho^*,m^*)$ is a Nash Equilibrium with $\rho^*>0$ and satisfies
\begin{equation*}
    (\rho^*,m^*) = \argmin_{(\rho,m)\in C} J(\rho,m).
\end{equation*}
Since the objective $J$ is convex, the primal gap
\begin{equation}
    h(\rho,m):= J(\rho,m) - J(\rho^*,m^*)
\label{eq: def h}
\end{equation}
quantifies the suboptimality of a candidate pair $(\rho, m)$ relative to the Nash equilibrium. The dual gap in this context corresponds exactly to the gain of the best response $g$ defined in \eqref{eq: def exploit}.

The proof techniques in \cite{Lavigne2023gcg-ficplay} draw heavily from the analysis of the GCG method in \cite{Kunisch2024gcg-fastcvg}, relying primarily on the decay of the primal gap $h$. Although the dual gap $g$ is also well-studied in the literature on Frank-Wolfe and GCG methods, its analysis is typically grounded in its relationship to the primal gap.
In non-potential MFGs, where such a relationship may not exist, directly analyzing the dual gap becomes essential. To this end, Section \ref{subsec: best response} of this paper establishes a crucial connection between the fixed point of the best response mapping and the dual gap (i.e., the gain of the best response), providing the foundation for extending convergence analysis from potential to general MFGs.

Although our convergence analysis is formulated directly in terms of the dual gap $g$, we first show that in the potential setting, $g$ is equivalent to the primal gap $h$. This connection serves as a consistency check and helps unify the analysis across both frameworks.
\begin{proposition}
\label{lem: primal = dual}
    Under Assumptions \ref{aspt: lip cts}, \ref{aspt: independent cost} and \ref{asp: f is potential}, if $F$ is convex on $\calP(\bbT^d)$ under $L^2$ metric and $R$ is $\lambda_R$-strongly convex between $(\rho,m)$ and its best response, 
    then
    \begin{equation}
        h(\rho,m)\leq g(\rho,m) \leq \max\left\{4,\frac{4L_f}{\lambda_R}\right\} h(\rho,m).
    \end{equation}
\end{proposition}

\begin{proof}
    We first show that $h(\rho,m)\leq g(\rho,m)$.
    \begin{equation}
    \begin{aligned}
        h(\rho,m)
        &= R(\rho,m)-R(\rho^*,m^*) + \int_0^T F(\rho_t)-F(\rho_t^*)\dd t + \langle f_T,\rho_T-\rho^*_T \rangle_{\bbT^d}\\
        &\leq R(\rho,m)-R(\rho^*,m^*) + \langle f(\rho),\rho-\rho^*\rangle_{Q_T} + \langle f_T,\rho_T-\rho^*_T \rangle_{\bbT^d}\\
        &\leq g(\rho,m).
    \end{aligned}
    \end{equation}
    Here, the first and third lines are by definition and the second line is by convexity of $F$.

    Next, we denote $(\hat{\rho},\hat{m})$ as the best response of $(\rho,m)\in C$ and
    $$(\rho^+,m^+)=(1-\delta)({\rho},{m}) + \delta(\hat{\rho},\hat{m}),$$ where $\delta\in(0,1]$.
    Then the convexity of $R$ and the average step give
    \begin{equation}
        R(\rho^+,m^+)-R(\rho,m)
        \leq \delta(R(\hat{\rho},\hat{m})-R(\rho,m)).
    \label{eq: r cvx}
    \end{equation}
    The Lipschitz continuity of $f$ (Assumption \ref{aspt: lip cts}) and the average step give 
    \begin{equation}
    \begin{aligned}
        \int_0^T F(\rho^+)-F(\rho) \dd t
        \leq & \langle f(\rho),\rho^+-\rho\rangle_{Q_T} + \frac{L_f}{2}\| \rho^+-\rho \|^2_{Q_T}\\
        = & \delta \langle f(\rho),\hat{\rho}-\rho\rangle_{Q_T} 
            + \delta^2\frac{L_f}{2} \| \hat{\rho}-\rho \|^2_{Q_T}\\
    \end{aligned}
    \label{eq: s L-smooth}
    \end{equation}
    By combining \eqref{eq: def h}, \eqref{eq: r cvx} and \eqref{eq: s L-smooth}, we have
    \begin{equation}
    \begin{aligned}
        &h(\rho^+,m^+) - h(\rho,m)\\
        \leq & \delta\Big( R(\hat{\rho},\hat{m}) 
                + \langle f(\rho),\hat{\rho} \rangle_{Q_T}
                + \langle f_T, \hat{\rho}_T \rangle_{\bbT^d} \Big)
         - \delta \Big( R(\rho,m) 
                + \langle f(\rho),\rho \rangle_{Q_T}
                + \langle f_T, \rho_T \rangle_{\bbT^d} \Big)
         + \delta^2\frac{L_f}{2} \| \hat{\rho}-\rho \|^2_{Q_T}\\
        =& -\delta g(\rho,m) + \delta^2\frac{L_f}{2}\|\hat{\rho}-\rho\|^2.
    \end{aligned}
    \label{eq: key for primal=dual}
    \end{equation}
    Since $h\geq 0$ and by Lemma \ref{lem: best response dist bounded by exploit}, we have
    \begin{equation}
        \left( \delta - \delta^2\frac{L_f}{\lambda_R} \right) g(\rho,m)
        \leq h(\rho,m).
    \end{equation}
    Taking $\delta=\min\left\{ \frac{1}{2},\frac{\lambda_R}{2L_f} \right\}$ concludes the lemma.
\end{proof}

In the rest of this section, we give a simplified convergence analysis for potential MFGs.
To prove the convergence for potential MFGs, we establish the decay property of the primal gap, which does not require the directional strong convexity of $R$.
\begin{lemma}[The decay property of $h$]
\label{lem: primal gap recur}
    Under assumptions \ref{aspt: lip cts}, \ref{aspt: independent cost} and \ref{asp: f is potential}, if $F$ is convex on $\calP(\bbT^d)$ under $L^2$ metric, then the sequence generated by fictitious play \eqref{eq: ficplay lin}-\eqref{eq: ficplay avg rhom} satisfies
    \begin{equation}
        h(\rho^{(k)},m^{(k)})
        \leq (1-\delta_k)h(\rho^{(k-1)},m^{(k-1)})
            + \delta_k^2\frac{L_f}{2} \| \hat{\rho}^{(k)}-\rho^{(k-1)} \|^2_{Q_T}.
    \end{equation}
\end{lemma}

\begin{proof}
    Repeating the analysis in \eqref{eq: r cvx}, \eqref{eq: s L-smooth} and \eqref{eq: key for primal=dual} gives us
    \begin{equation}
    \begin{aligned}
        h(\rho^{(k)},m^{(k)}) - h(\rho^{(k-1)},m^{(k-1)})
        \leq & \delta_k\left( R(\hat{\rho}^{(k)},\hat{m}^{(k)}) 
                + \langle f(\rho^{(k-1)}),\hat{\rho}^{(k)} \rangle_{Q_T}
                + \langle f_T, \hat{\rho}_T^{(k)} \rangle_{\bbT^d} \right)\\
        & - \delta_k \left( R(\rho^{(k-1)},m^{(k-1)}) 
                + \langle f(\rho^{(k-1)}),\rho^{(k-1)} \rangle_{Q_T}
                + \langle f_T, \rho_T^{(k-1)} \rangle_{\bbT^d} \right)\\
        & + \delta_k^2\frac{L_f}{2} \| \hat{\rho}^{(k)}-\rho^{(k-1)} \|^2_{Q_T}\\
    \end{aligned}
    \label{eq: 89}
    \end{equation}
    The optimality of $(\hat{\rho}^{(k)},\hat{m}^{(k)})$ gives
    \begin{equation}
        R(\hat{\rho}^{(k)},\hat{m}^{(k)}) 
                + \langle f(\rho^{(k-1)}),\hat{\rho}^{(k)} \rangle_{Q_T}
                + \langle f_T, \hat{\rho}_T^{(k)} \rangle_{\bbT^d}
        \leq R(\rho^*,m^*) 
                + \langle f(\rho^{(k-1)}),\rho^* \rangle_{Q_T}
                + \langle f_T, \rho_T^* \rangle_{\bbT^d}
    \end{equation}
    and the convexity of $F$ gives
    \begin{equation}
        \langle f(\rho^{(k-1)}),\rho^* -\rho^{(k-1)} \rangle_{Q_T}
        \leq \int_0^T F(\rho^*_t)-F(\rho^{(k-1)}_t)\dd t.
    \end{equation}
    Plugging into \eqref{eq: 89}, we have
    \begin{equation}
        h(\rho^{(k)},m^{(k)}) - h(\rho^{(k-1)},m^{(k-1)})
        \leq  -\delta_k h(\rho^{(k-1)},m^{(k-1)})
        + \delta_k^2\frac{L_f}{2} \| \hat{\rho}^{(k)}-\rho^{(k-1)} \|^2_{Q_T},
    \end{equation}
    which concludes the lemma.
\end{proof}

By coercivity of the objective function $J$, we can restrict our focus to a subset of $C$ with diameter $D<\infty$.
According to \cite{braun2023gcg-cgbook}, $\delta_k=\frac{2}{k+2}$ is known as the open-loop step-size and is a standard choice in Frank-Wolfe when the information of the function is not available.
Taking $\delta_k=\frac{2}{k+2}$ in Lemma \ref{lem: primal gap recur}, we have the convergence for potential games.
We emphasize that the sublinear convergence of the primal gap does not require the directional strong convexity of $R(\rho,m)$ in $\rho$. 
\begin{theorem}
\label{thm: pot cvg rate}
    Under assumptions \ref{aspt: lip cts}, \ref{aspt: independent cost} and \ref{asp: f is potential}, if $F$ is convex on $\calP(\bbT^d)$ under $L^2$ metric and suppose that the sequence generated by fictitious play \eqref{eq: ficplay lin}–\eqref{eq: ficplay avg} with step size $\delta_k = \frac{2}{k+2}$ remains within a subset of $C$ with diameter $D$. Then the sequence satisfies the bound
    \begin{equation}
        h({\rho}^{(k)},{m}^{(k)})
        \leq \frac{2}{(k+2)(k+1)}h({\rho}^{(0)},{m}^{(0)}) + \frac{2L_fD^2}{k+2}.
    \label{eq: cvg rate}
    \end{equation}
\end{theorem}
\begin{proof}

    Since the upper-bound of $\| \hat{\rho}^{(k)}-{\rho}^{(k-1)} \|$ is D, from Lemma \ref{lem: primal gap recur} we have
    \begin{equation}
        h({\rho}^{(k)},{m}^{(k)})
        \leq \frac{k}{k+2}h({\rho}^{(k-1)},{m}^{(k-1)})
            + \frac{2L_fD^2}{(k+2)^2}.
    \label{eq: recur}
    \end{equation}  
    Telescoping \eqref{eq: recur} concludes the theorem.
\end{proof}

When the Mean Field Game (MFG) is potential, other optimization algorithms have been employed to solve the variational formulation. 
The convergence of primal-dual \cite{liu2021nonlocal} depends on the regularity of the constraint. For MFG, the constraint involves differential operators that are unbounded.
When solving the discretized system, the condition number of the derivative operator grows to infinity as the mesh step decreases to 0 and the primal-dual will take more iterations to converge.
Compared to primal-dual \cite{liu2021nonlocal}, our convergence is independent of the regularity of the constraint and therefore the mesh step. More details are in Remark \ref{rm: disct equilibrium form}.
Compared to ADMM and augmented Lagrangian \cite{benamou2015augmented,benamou2017variational}, the subproblem in fictitious play can be solved with better accuracy in less time.

\section{Discretization and Implementation}
\label{sec: imple}

This section first discusses the finite difference approximation of \eqref{eq: mfg} on $[0,1]$. This discretization can be easily generalized to a cuboid in high dimensions. Based on the discretization, we propose a fictitious play accelerated by employing a hierarchical grid for solving the discrete system.

Let $n_t$ be a positive integer and $\dt=\frac{T}{n_t},t_n=n\dt, n=0,\cdots,n_t$.
Similarly, let $n_x$ be a positive integer and $\dx=\frac{1}{n_x},x_i=i\dx, i=0,\cdots,n_x$.
Denote $\calG$ as the collection of all grid points $(x_i,t_n),i=0,\cdots,n_x,n=0,\cdots,n_t$.
For any function $u$ defined on $\bbT^d\times[0,T]$, we denote the approximation of $u$ on $\calG$ as $u_{\calG}$ where $(u_{\calG})_\myij$ approximates $u(x_i,t_n)$.
When the context provides no ambiguity, we omit $\calG$.
In this section, we consider the homogeneous Neumann boundary condition for the HJB equation and the zero-flux boundary condition for the FP equation. To tackle this in the discretization, we introduce the following finite difference operators on the grid $\calG$.
Firstly, the elementary one-sided finite difference operators are defined as
\begin{equation*}
    (\dxp u)_\myij = \begin{cases}
        \frac{u_\ipj-u_\myij}{\dx}, & i=0,\cdots,n_x-1,\\
        0, & i=n_x,
    \end{cases}
    \quad
    (\dxm u)_\myij = \begin{cases}
        0, & i=0,\\
        \frac{u_\myij-u_\imj}{\dx}, & i=1,\cdots,n_x.
    \end{cases}
\end{equation*}
Then we define the central difference operator $(\dxc u)_\myij$ as the average of the one-sided finite difference, $[D_xu]_\myij$ as the collection of the one-sided finite differences at $(x_i,t_n)$
\begin{equation*}
    (\dxc u)_\myij = \frac{1}{2}(\dxp u)_\myij + \frac{1}{2}(\dxm u)_\myij,\quad
    [D_xu]_\myij = \left((\dxp u)_\myij,(\dxm u)_\myij\right),
\end{equation*}
and the discrete Laplace operator as
\begin{equation*}
    (\Delta_x u)_\myij =
    \begin{cases}
        \frac{u_\ipj-u_\myij}{\dx^2}, & i=0,\\
        \frac{u_\ipj-2u_\myij+u_\imj}{\dx^2}, & i=1,\cdots,n_x-1\\
        \frac{-u_\myij+u_\imj}{\dx^2}, & i=n_x\\
    \end{cases}
\end{equation*}
We call $u=(u_\myij)_{i=0,n=0}^{i=n_x,n=n_t}$ a function on $\calG$ and $[v]=(v^+,v^-)$ a velocity field on $\calG$ where $v^+,v^-$ are functions on $\calG$.
For any functions $u_1,u_2$ and vector fields $[v_1]=(v_1^+,v_1^-),[v_2]=(v_2^+,v_2^-)$ on $\calG$, we define the inner product on the grid $\calG$ as
\begin{equation*}
    \langle u_1,u_2 \rangle_{\calG}:=\Delta t\Delta x\sum_{i=0}^{n_x}\sum_{n=0}^{n_t} (u_1)_\myij(u_2)_\myij,\quad
    \left\langle [v_1],[v_2] \right\rangle_{\calG}:=
    \frac{1}{2}\left(
    \left\langle v_1^+, v_2^+ \right\rangle_{\calG} + 
    \left\langle v_1^-, v_2^- \right\rangle_{\calG}
    \right)
\end{equation*}
The discrete Laplacian operator is self-adjoint under the above inner product.
To preserve the adjoint relation between gradient and negative divergence, for a given velocity field $v$ on $\calG$, we define its divergence to be $-D_x^*[v]$ where $D_x^*$ is the adjoint operator of $D_x$,
\begin{equation*}
    (D_x^*[v])_\myij := 
    \begin{cases}
        -\frac{1}{2}\frac{(v^+)_{0,n} + (v^-)_{1,n}}{\dx}, & i=0,\\
        -\frac{1}{2}\left(\dxm(v^+) + \dxp(v^-)\right)_\myij, & i=1,\cdots,n_x-1,\\
        \frac{1}{2}\frac{(v^-)_{n_x,n} + (v^+)_{n_x-1,n}}{\dx}, & i=n_x,
    \end{cases}    
\end{equation*}
and then
\begin{equation*}
    \langle u,D_x^*[v] \rangle_{\calG} = \left\langle [D_x(u)],[v] \right\rangle_{\calG}.
\end{equation*}

We consider the Lax-Friedrichs Hamiltonian
\begin{equation*}
    \hamlf(x,\pxp,\pxm)
    =H\left( x, \left( \frac{\pxp+\pxm}{2} \right) \right)
    -\nu_n\frac{\pxp-\pxm}{2},
\end{equation*}
where $\nu_n$ is the numerical viscosity coefficient.

We approximate the interaction cost $f(x,\rho)$ and $f_T(x,\rho)$ at $x_i$ by $f_{\calG}(x_i,\rho_{\calG})$ and $f_{T,\calG}(x_i,\rho_{\calG})$, respectively.
For any function $\rho_{\calG}$ on grid $\calG$, we define its corresponding piece-wise constant function on $\bbT^d$ as $\bar{\rho}_{\calG}$ with
$\bar{\rho}_{\calG}(x) = \rho_i, x_i-\frac{\dx}{2}\leq x < x_i+\frac{\dx}{2}$.
For consistency, we assume that there exists a constant $C$ independent of $\dx$ and $\rho$ such that for any $x_i$ and $\rho_{\calG}$, 
$$|f_{\calG}(x_i,\rho_{\calG})-f(x_i,\bar{\rho}_{\calG})|<C\dx,
\quad |f_{T,\calG}(x_i,\rho_{\calG})-f_T(x_i,\bar{\rho}_{\calG})|<C\dx.$$ 
Below are some examples of cost functions and their consistent discretizations. 
\begin{enumerate}
    \item If $f(x,\rho)=f_0(\rho(x))$ is defined locally by some function $f_0:\bbR\to\bbR$, we can define $f_{\calG}(x_i,\rho_{\calG})=f_0(\rho_i)$.
    \item If $f(x,\rho)=\int_0^1 K(x,y)\rho(y)\dd y$ is a nonlocal cost defined by some kernel function $K$, it is usually approximated by $f_{\calG}(x_i,\rho_{\calG})=\sum_{j=0}^{n_x} K(x_i,x_j)\rho_j\dx$
    \item In section \ref{sec: num result}, we consider a cost of form $f(x,\rho)=f_0(x,\mu_1(\rho),\mu_2(\rho))$ where $\mu_1(\rho)=\int_0^1 x\rho(x)\dd x$ is the first moment and $\mu_2(\rho)=\int_0^1 x^2\rho(x)\dd x$ is the second moment of $\rho$.
    We first approximate the first and second moments with $\mu_{1,\calG}(\rho_{\calG})=\sum_{i=0}^{n_x}x_i\rho_i\dx$ and 
    $\mu_{2,\calG}(\rho_{\calG})=\sum_{i=0}^{n_x}x_i^2\rho_i\dx$ and then define $f_{\calG}(x_i,\rho_{\calG})=f_0(x_i,\mu_{1,\calG}(\rho_{\calG}),\mu_{2,\calG}(\rho_{\calG}))$. 
\end{enumerate}    

With the above notation, we denote the residues of the HJB equation and FP equation as follows,
\begin{equation}
    (\reshjb(\rho,\phi))_\myij:=
    \left\{
    \begin{aligned}
        &-\frac{\phi_\ijp-\phi_\myij}{\dt} - \nu(\Delta_x \phi)_\myij + \hamlf(x_i,[D_x\phi]_\myij)-f(x_i,\rho_{\cdot,n+1}),\\
        &\qquad\qquad\qquad\qquad\qquad\qquad i=0,1,\cdots,n_x,n=0,1,\cdots,n_t-1,\\
        &\phi_{i,n_t} - f_T(x_i,\rho_{\cdot,n_t}),\qquad i=0,1,\cdots,n_x,
    \end{aligned}
    \right.
\label{eq: res hjb}
\end{equation}
\begin{equation}
    (\resfp(\rho,\phi))_\myij:=
    \left\{
    \begin{aligned}
        &\frac{\rho_\ijp-\rho_\myij}{\dt} 
        - \nu(\Delta_x \rho)_\ijp + (D_x^*[\rho_{\cdot,n+1} v_{\cdot,n}])_i ,\\
        &\qquad \text{where }v^+_\myij = -\nabla_{p^+}\hamlf(x_i,[D_x\phi]_\myij), \\
        &\qquad\qquad\quad  v^-_\myij = -\nabla_{p^-}\hamlf(x_i,[D_x\phi]_\myij),\\
        &\qquad\qquad\qquad\qquad\qquad\qquad i=0,1,\cdots,n_x,n=0,1,\cdots,n_t-1,\\
        &\rho_{i,0} - \rho_0(x_i),\qquad i=0,1,\cdots,n_x.\\
    \end{aligned}
    \right.
\label{eq: res fp}
\end{equation}
And our goal is to find $\rho,\phi$ such that 
\begin{equation}
    \reshjb(\rho,\phi)=\mathbf{0},\resfp(\rho,\phi)=\mathbf{0}.
\label{eq: disct mfg}
\end{equation}

\begin{remark}[Existence and uniqueness of the solution to \eqref{eq: disct mfg}]
The Lax-Friedrichs Hamiltonian has the following nice properties
\begin{enumerate}
    \item $\hamlf$ is consistent with $H$, i.e. $\hamlf(x,p,p) = H(x,p)$.
    
    \item $\hamlf$ is nonincreasing w.r.t. $p^+$ and nondecreasing w.r.t $p^-$ given $\nu_n\geq \frac{1}{2}\left|\nabla_p H(x,\frac{p^+ + p^-}{2})\right|$. 
    
    \item $\hamlf$ is differentiable w.r.t. $p^+,p^-$ given $H$ being differentiable w.r.t. $p$.
    
    \item $\hamlf$ is jointly convex in $p^+, p^-$ given $H$ being convex in $p$.
        
\end{enumerate}
According to \cite{achdou2010mean,Achdou2013fdcvg}, 1-3 are needed to establish the existence of the solution to the discrete MFG system \eqref{eq: disct mfg} and 1-4 are needed to establish the uniqueness of the solution.
For the details, we refer readers to \cite{achdou2010mean,Achdou2013fdcvg}. It is also established therein that the solution to the discrete system \eqref{eq: disct mfg} converges to the solution of the continuous system \eqref{eq: mfg} as the mesh step approaches zero.

\end{remark}

\begin{remark}[The discretization of Fokker-Planck equation]
    Requiring $\resfp(\rho,\phi)=\mathbf{0}$ is equivalent to applying a semi-implicit Lax-Friedrichs scheme to the Fokker-Planck equation with zero-flux boundary condition, i.e.
    $$\frac{\rho_\ijp-\rho_\myij}{\dt} 
    - \left(\nu+\frac{\nu_n}{2}\dx\right)(\Delta_x \rho)_\ijp 
    + (D_x^*[\rho_{\cdot,n+1}{v}_{\cdot,n}])_i =0,$$
    where $ {v}^{\pm}_\myij = -\nabla_pH(x_i,(\dxc \phi)_\myij).$
    Under this scheme, the discrete mass is conserved, i.e.  
    $ \sum_{i=0}^{n_x} \rho_\myij = \sum_{i=0}^{n_x} \rho_\ijp $ holds for any $n=0,1,\cdots,n_t-1$.    
\end{remark}

\begin{remark}[Semi-implicit schemes]
    The discrete MFG system \eqref{eq: disct mfg} is semi-implicit.
    To be precise, the discrete HJB \eqref{eq: res hjb} is implicit in $\phi$ and explicit in $\rho$ while the discrete FP \eqref{eq: res fp} is implicit in $\rho$ and explicit in $\phi$.
    While one can also choose the implicit scheme, the semi-implicit scheme preserves the time-adjoint relation between $\rho$ and $\phi$, which is essential for preserving the variational formulation and equilibrium formulation in discretization. 
    Remark \ref{rm: disct pot mfg} and \ref{rm: disct equilibrium form} give more concrete explanations.
    In addition to the above advantage, according to \cite{achdou2010mean}, it is easier to show the uniqueness of the solution to the discrete MFG system with the semi-implicit scheme.
\end{remark}

\begin{remark}[Variational formulation for potential MFGs]
\label{rm: disct pot mfg}
Notice that if $L$ is the convex conjugate of $H$, then
\begin{equation*}
    \hamlf(x,\pxp,\pxm)=\sup_{q} \left\{ \frac{\pxp+\pxm}{2} q-L(x,q)\right\} -\nu_n\frac{\pxp-\pxm}{2}.
\end{equation*}
With this, when there exist $F,F_T$ mapping functions on $\calG$ to $\bbR$, such that 
\begin{equation*}
    \frac{\partial F}{\partial\rho_\myij}\left( \rho_{\cdot,n} \right) = \dx f(x_i,\rho_{\cdot,n}),\quad
    \frac{\partial F_T}{\partial\rho_\myij}\left( \rho_{\cdot,n} \right) = \dx f_T(x_i,\rho_{\cdot,n}),
\end{equation*}
the discrete system \eqref{eq: disct mfg} is the optimality condition of the optimization problem
\begin{equation}
\begin{aligned}
    \min_{(\rho,m)\in C_{\calG}}~ & \dt\dx\sum_{i=0}^{n_x}\sum_{n=0}^{n_t-1}\rho_\ijp L\left( x_i, \frac{m_\ijp}{\rho_\ijp} \right) 
     +\dt \sum_{n=0}^{n_t-1}F\left( \rho_{\cdot,n+1} \right) + F_T\left( \rho_{\cdot,n_t} \right) \\
\end{aligned}
\label{eq: potential mfg disct}
\end{equation}
with 
\begin{equation}
C_{\calG}:=\left\{(\rho,m):
\begin{aligned}
    &\frac{\rho_\ijp-\rho_\myij}{\dt} 
    - \left(\nu+\frac{\nu_n}{2}\dx\right)(\Delta_x \rho)_\ijp + (D_x^*(m,m))_\ijp =0,\\
    &\qquad \qquad i=0,\cdots,n_x,n=0,\cdot,n_t-1,\\
    & \rho_{i,0} = \rho_0(x_i), i=0,\cdots,n_x.     
\end{aligned}
\right\}
\label{eq: disct cst set}
\end{equation} being the discrete constraint set and
$\phi_\myij$ being the dual variable to the $(i,n)$-th constraint of $\rho,m$.
Such a variational formulation does not exist if we discretize HJB and FP with implicit schemes.
\end{remark}

\subsection{Fictitious play for the discrete system}
\label{subsec: disct ficplay}

While requiring $\reshjb(\rho,\phi)=\mathbf{0}$ or $\resfp(\rho,\phi)=\mathbf{0}$ can be addressed backward or forward with relative ease, satisfying both simultaneously is more complex. There is no straightforward time-marching method to achieve this; consequently, the entire system must be solved concurrently. This is challenging due to the presence of nonlinear Hamiltonian and cost terms.
To numerically solve the discrete system \eqref{eq: disct mfg}, we alternatively solve the HJB \eqref{eq: res hjb} and the FP \eqref{eq: res fp} equations. 
Algorithm \ref{alg: disct ficplay} summarizes the procedure. In practice, we solve the HJB backward with Newton's method.

\begin{algorithm}[htb]
\caption{Fictitious Play for Discrete System \eqref{eq: disct mfg}}
\begin{algorithmic}
    \STATE{Parameters} $\rho_0,f,f_T,\nu,\nu_n, 0<\delta_k\leq 1$
    \STATE{Initializations} $\phi^{(0)}$
    \STATE{Given ${\phi}^{(0)}$, solve the discretized FP $\resfp\left({\rho}^{(0)},{\phi}^{(0)}\right)=\mathbf{0}$ for ${\rho}^{(0)}$ forward in $t$.}
    \FOR{$k=1,2,\cdots,K$}
    
        \STATE{Given $\rho^{(k-1)}$, solve the discretized HJB $\reshjb\left(\rho^{(k-1)},\hat{\phi}^{(k)}\right)=\mathbf{0}$ for $\hat{\phi}^{(k)}$ backward in $t$.} 
               
        \STATE{Given $\hat{\phi}^{(k)}$, solve the discretized FP $\resfp\left(\hat{\rho}^{(k)},\hat{\phi}^{(k)}\right)=\mathbf{0}$ for $\hat{\rho}^{(k)}$ forward in $t$.}     
        
        \STATE{Execute pointwise density average \eqref{eq: ficplay disct avg} to obtain $\rho^{(k)},\phi^{(k)}$}
        \begin{equation}
            \rho^{(k)}=\left(1-\delta_k\right)\rho^{(k-1)}+\delta_k\hat{\rho}^{(k)}.
        \label{eq: ficplay disct avg}
        \end{equation}
    \ENDFOR
    \STATE{Output} $\rho^{(K)},\hat{\phi}^{(K)}$.
\end{algorithmic}
\label{alg: disct ficplay}
\end{algorithm}

\begin{remark}[Equilibrium form for the discrete system \eqref{eq: disct mfg} and reformulation of Algorithm \ref{alg: disct ficplay}]
\label{rm: disct equilibrium form}
    Denote the discrete counterpart of $J$ as
    \begin{equation*}
    \begin{aligned}
        J_{\calG}(\tilde{\rho},\tilde{m};\rho):=
        &\dt\dx\sum_{i=0}^{n_x}\sum_{n=0}^{n_t-1}\tilde{\rho}_\ijp L\left( x_i, \frac{\tilde{m}_\ijp}{\tilde{\rho}_\ijp} \right) \\
        & +\dt\dx\sum_{i=0}^{n_x}\sum_{n=0}^{n_t-1}f\left(x_i,\rho_{\cdot,n+1} \right) \tilde{\rho}_{\ijp} 
         +\dx\sum_{i=1}^{n_x} f_T\left(x_i, \rho_{\cdot,n_t} \right) \tilde{\rho}_{i,n_t}.
    \end{aligned}
    \end{equation*}
    In the same spirit of Remark \ref{rm: disct pot mfg}, one can show that with mild assumptions, $(\rho,\phi)$ solves the discrete MFG system \eqref{eq: disct mfg} if and only if 
    \begin{equation*}
        m_\ijp = -\rho_\ijp(\dxc\phi)_\myij
    \end{equation*}
    and $(\rho,m)$ solves the equilibrium problem
    \begin{equation*}
        (\rho,m)=\argmin_{(\tilde{\rho},\tilde{m})\in C_{\calG}}~  J_{\calG}(\tilde{\rho},\tilde{m};\rho). 
    \end{equation*}    
    Consequently, the discrete fictitious play (Algorithm \ref{alg: disct ficplay}) can be written as
    \begin{align}
        &(\hat{\rho}^{(k)},\hat{m}^{(k)})\in\argmin_{(\tilde{\rho},\tilde{m})\in C_{\calG}} J_{\calG}(\tilde{\rho},\tilde{m};\rho^{(k-1)}),
        \label{eq: ficplay lin disct}\\
        &(\rho^{(k)},m^{(k)})=(1-\delta_k)(\rho^{(k-1)},m^{(k-1)}) + \delta_k (\hat{\rho}^{(k)},\hat{m}^{(k)}).
        \label{eq: ficplay avg rhom disct}
    \end{align}
    This is the discrete counterpart of \eqref{eq: ficplay lin}-\eqref{eq: ficplay avg rhom} for the continuous fictitious play (Algorithm \ref{alg: ficplay}). 
    Therefore, the convergence analysis for Algorithm \ref{alg: ficplay} also holds for Algorithm \ref{alg: disct ficplay} under sufficient regularity assumptions and fine meshes.
\end{remark}

It is important to notice that the convergence rate depends only on the Lipschitz constant $L_f$ and convexity $\lambda_R$, which are $\mathcal{O}(1)$ with respect to $\dx,\dt$. Therefore, the convergence rate of the discrete fictitious play (Algorithm \ref{alg: disct ficplay}) does not depend on the mesh step. 
This is a big advantage compared to primal-dual since the convergence rate of primal-dual depends on the operator norm of constraint and in our case, the differential operator norm is $\mathcal{O}(\frac{1}{\min\{\Delta t,\Delta x\}})$.

\subsection{Acceleration: incorporating hierarchical grids}
\label{subsec: multigrid}

In this section, we describe how we incorporate a hierarchical grid to accelerate Algorithm \ref{alg: disct ficplay}.

We notice that the MFG system is essentially a two-point boundary value problem and fictitious play decomposes it into two initial value problems, which can be viewed as a shooting method. 
Thus when solving the discrete system \eqref{eq: disct mfg}, it is natural first to solve it on a coarse grid to let the boundary information flow into the interior and refine the grid to improve solution accuracy.
This can reduce the computational complexity at the beginning, provide a good initialization on fine grids, and reduce the number of iterations on fine grids.
With the Lax-Friedrichs scheme, which incorporates a numerical viscosity term proportional to the mesh step, our hierarchical grid approach can stabilize the algorithm for first-order MFGs where the system is dominated by advection.

Precisely, consider the same setting at the beginning of this section. 
Let $\dx,\dt$ be the mesh step on the coarsest grid $\calG_0$. 
We refine the grid by $\frac{1}{2}$ and obtain a hierarchy grids $\calG_l\subset\calG_{l+1},l=0,\cdots,L-1$ where the mesh step of $\calG_l$ is $2^{-l}(\dx,\dt)$.
We define the prolongation operator $\pro$ such that for any given function $u^{\calG_l}$ on $\calG_l$, $u^{\calG_{l+1}}:=\pro\left(u^{\calG_l}\right)$ is a function on $\calG_{l+1}$, and 
$$u^{\calG_{l+1}}_{i,n}=\begin{cases}
    u^{\calG_l}_{\frac{i}{2},\frac{n}{2}}, & \text{ $i$ and $n$ are even, }\\
    \frac{1}{2} u^{\calG_l}_{\frac{i-1}{2},\frac{n}{2}} + \frac{1}{2} u^{\calG_l}_{\frac{i+1}{2},\frac{n}{2}}, & \text{ $i$ is odd, $n$ is even, }\\
    \frac{1}{2} u^{\calG_l}_{\frac{i}{2},\frac{n-1}{2}} + \frac{1}{2} u^{\calG_l}_{\frac{i}{2},\frac{n+1}{2}}, & \text{ $i$ is even, $n$ is odd, }\\
    \frac{1}{4} u^{\calG_l}_{\frac{i-1}{2},\frac{n-1}{2}} + \frac{1}{4} u^{\calG_l}_{\frac{i+1}{2},\frac{n-1}{2}}
    + \frac{1}{4} u^{\calG_l}_{\frac{i-1}{2},\frac{n+1}{2}} + \frac{1}{4} u^{\calG_l}_{\frac{i+1}{2},\frac{n+1}{2}}, & \text{ $i$ and $n$ are odd. }\\    
\end{cases}$$
With the above notations and Algorithm \ref{alg: disct ficplay}, we summarize our hierarchical grid fictitious play in Algorithm \ref{alg: multigrid}.
\begin{algorithm}[ht]
\caption{Hierarchical Grid Fictitious Play}
\begin{algorithmic}
    \STATE{Parameters} $\epsilon,$ $\rho_0,f,f_T,\nu,\nu_n$
    \STATE{Initialization} ${\phi}^{\calG_0}$
    \FOR{$l=1,\cdots,L$}
    \STATE{Interpolate ${\phi}$: } ${\phi}^{\calG_l,(0)}=P_{l-1}^{l}(\hat{\phi}^{\calG_{l-1}}).$
    \STATE{Run Algorithm \ref{alg: disct ficplay} on $\calG_l$ with initialization $\phi^{\calG_l,(0)}$ to accuracy $10^{L-l}\epsilon$.}
    \STATE{Save the result:} $\rho^{\calG_l}=\rho^{\calG_l,(K_l)},
    ~{\phi}^{\calG_l}={\phi}^{\calG_l,(K_l)}$.
    \ENDFOR
    \STATE{Output} $\rho^{\calG_L},{\phi}^{\calG_L}$.
\end{algorithmic}
\label{alg: multigrid}
\end{algorithm}

\begin{remark}[Interpolation]
    While Algorithm \ref{alg: multigrid} interpolates $\phi$, we remark that one can also interpolate $\rho$ as the initialization of the fine grid.
    In our numerical experiment, because the true boundary condition of the HJB is not always consistent with the Neumann boundary that we enforce, interpolating $\rho$ causes the HJB residue to be large on the fine grid at the beginning.
    Interpolating $\phi$ on the one hand preserves the adjustment on the coarse grid. On the other hand, since the value of $\rho$ is usually very small on space boundaries, the velocity field derived from $\phi$ has little impact on the value of $\rho$. 
    In addition, initializing the algorithm with $\phi$ gives the exact velocity field that induces $\rho^{(0)}$ and helps to compute the gain function $g$ at the beginning and makes back-tracking line search accessible.
\end{remark}

\begin{remark}
    Our method uses a hierarchy of grids but is different from multigrid methods \cite{Trottenberg2000multigrid} used in \cite{achdou2012multigrid,andreev2017multigridadmm,briceno2018multigridadmm}, where they have both prolongation and restriction.
    Essentially, our hierarchical grid approach captures low-frequency modes on coarse grids first and gradually adds the high-frequency modes on fine grids. On the contrary, multigrid methods first capture the high-frequency modes on fine grids and reduce the low-frequency error on coarse grids. 
    The idea of our hierarchy grid acceleration is similar to the multilevel approach in \cite{yu2024pgdeucmfg}. But since their algorithm couples $\rho$ and $m$, their interpolation acts on both variables, while we only interpolate $\rho$ or $\phi$. 
\end{remark}

We conclude this section with a comparison of computational complexity in Table \ref{tab: complexity} when using the weight $\delta_k = \frac{1}{k+1}$ in the fictitious play. From Theorem \ref{thm: non-pot cvg rate}, it follows that asymptotically, $g(\rho^{(k)}, m^{(k)}) = \mathcal{O}\left(\frac{1}{k}\right)$.
Suppose our objective is to find $\rho^{\calG_L}, m^{\calG_L}$ such that $g(\rho^{\calG_L}, m^{\calG_L}) \leq \epsilon$ on grid $\calG_L$ with $n_x$ spatial points and $n_t$ time points, where $n_t = \mathcal{O}(n_x)$ (due to the use of implicit scheme in time) and $\frac{1}{\epsilon} = \mathcal{O}(n_x)$ (due to the first order scheme in space for HJ equation). 
For the hierarchical grid strategy (Algorithm \ref{alg: multigrid}), we choose $L = \mathcal{O}(\log_2 n_x)$. On each level $l=1, 2, \ldots, L$, one only needs to achieve an accuracy comparable with the grid size $\mathcal{O}(2^{-l})=\mathcal{O}(2^{L-l}\epsilon)$. Since we use the solution on grid level $l-1$ as the initial guess for grid level $l$, one only needs to reduce the error by $\frac{1}{2}$ from the initial error and hence requires $\mathcal{O}(1)$ iterations on each level. Each iteration on grid level $l$ has an operation count of $\mathcal{O}([2^{l-L}n_x]^{d+1})$, where $d$ is the spatial dimension.
As a result, the hierarchical grid strategy reduces the complexity of vanilla fictitious play from $\mathcal{O}(n_x^{d+2})$ to $\mathcal{O}(n_x^{d+1})$.

\begin{table}[h]
\caption{Complexity Comparison of Discrete Fictitious Play and the Hierarchical Grid Fictitious Play}
\label{tab: complexity}
\centering
\begin{tabular}{cc|cc}
\toprule
\multicolumn{2}{c|}{\makecell{Discrete fictitious play \\ (Algorithm \ref{alg: disct ficplay})}} & \multicolumn{2}{c}{\makecell{Hierarchical grid fictitious play \\ (Algorithm \ref{alg: multigrid})}}   \\ \midrule
per iteration      & iterations         & per iteration  & iterations  \\ \midrule
\multirow{4}{*}{$\mathcal{O}(n_x^{{d+1}})$} & \multirow{4}{*}{$\displaystyle \mathcal{O}(n_x)$} & $\mathcal{O}([2^{-L}n_x]^{d+1})$              & $\mathcal{O}(1)$               \\
                   &                    & $\mathcal{O}([2^{l-L}n_x]^{d+1})$              & $\mathcal{O}(1)$                    \\
                   &                    & $\vdots$              & $\vdots$                    \\
                   &                    & $\mathcal{O}(n_x^{d+1})$              & $\mathcal{O}(1)$           \\ \midrule
\multicolumn{2}{c|}{total: $\displaystyle \mathcal{O}(n_x^{d+2})$}                                         & \multicolumn{2}{c}{total: $\displaystyle \mathcal{O}(n_x^{d+1})$}                           \\          \bottomrule
\end{tabular}
\end{table}

\section{Numerical Results}
\label{sec: num result}

In this section, we conduct various experiments to show the effectiveness and efficiency of our hierarchical grid fictitious play algorithm. 
All of our numerical experiments are implemented in MATLAB on a PC with an Apple M3 Pro chip and 18 GB of memory.

For all of our experiments, we take $T=1$.
We denote the density function of a univariate Gaussian distribution as $\rho_G(x;\mu,\sigma):=\frac{1}{\sqrt{2\pi}\sigma}\exp{-\frac{1}{2}(\frac{x-\mu}{\sigma})^2}$.
Let the inner product and norm on the discrete grid $\calG$ be
$$\langle u,v \rangle_{\calG}:= \Delta x\Delta t \sum_{i,n}u_{i,n}v_{i,n},\quad
\left\| u \right\|_{\calG} :=  \sqrt{\langle u,u \rangle_{\calG} }.$$
We monitor three criteria that measure convergence.
The first is the gain of the best response
$$g^{(k)}:=J({\rho}^{(k)},{\phi}^{(k)};\rho^{(k)})
-J(\hat{\rho}^{(k+1)},\hat{\phi}^{(k+1)};\rho^{(k)}),$$
where
\begin{equation*}
\begin{aligned}
    J(\tilde{\rho},\tilde{\phi};\rho):=
    &\dt\dx\sum_{i=0}^{n_x}\sum_{n=0}^{n_t-1}\tilde{\rho}_\ijp L\left( x_i,-D_pH(x_i,(D_x^c\tilde{\phi})_\myij) \right)\\
    &+ \dt\dx\sum_{i=0}^{n_x}\sum_{n=0}^{n_t-1}f(x_i,\rho_{\cdot,n+1})\tilde{\rho}_{i,n+1} 
    + \dx\sum_{i=0}^{n_x}f_T(x_i,\rho_{\cdot,n_t})\tilde{\rho}_{i,n_t}.
\end{aligned}
\end{equation*}
While theoretically $g^k\geq 0$, numerically it may become negative when $g^{(k)}$ is close to 0. 
Therefore, we use $|g^{(k)}|\leq\epsilon$ as the stopping criterion and when presenting the semilogy plot of $g^{(k)}$, we use $|g^{(k)}|$.
The second is the consecutive residue $\|\hat{\rho}^{(k)}-\rho^{(k-1)}\|_{\calG}$. As discussed in Section \ref{sec: ficplay cvg}, under certain conditions, the distance between $\rho^{(k-1)}$ and the Nash Equilibrium $\rho^*$ is bounded by the consecutive residue.
The third measurement is FP residue $ \left\| \resfp\left( \rho^{(k)},\hat{\phi}^{(k+1)} \right) \right\|_{\calG}.$
Since $\hat{\phi}^{(k+1)}$ is obtained from solving HJB whose source term is given by $\rho^{(k)}$, the HJB residue $\left\| \reshjb\left( \rho^{(k)},\hat{\phi}^{(k+1)} \right) \right\|_{\calG}$ is only due to the discretization and HJB solver. The FP residue reflects how well $\rho^{(k)}$ approaches the fixed point of the best response mapping and should be small when $\rho^{(k)}$ is close to a Nash Equilibrium.

While the theoretical convergence result is established only for monotone games and density-independent terminal cost, the convergence analysis suggests that, provided that the terms due to non-monotone cost and dependent $f_T$ are adequately controlled, the fictitious play still produces a decreasing sequence of $g^{(k)}$ and thus converges.
Our numerical experiments also demonstrate success in non-monotone MFGs and more general forms of $f_T$.

\subsection{Divergence of fixed-point iteration}
\label{subsec: num fixpt}
\begin{figure}[htbp]
    \centering
    \includegraphics[width=5cm]{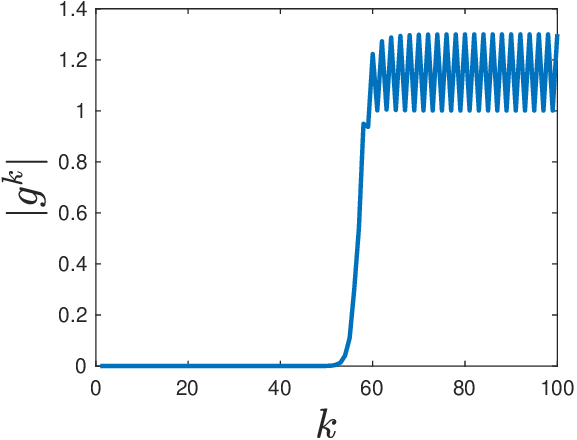}
    \includegraphics[width=5cm]{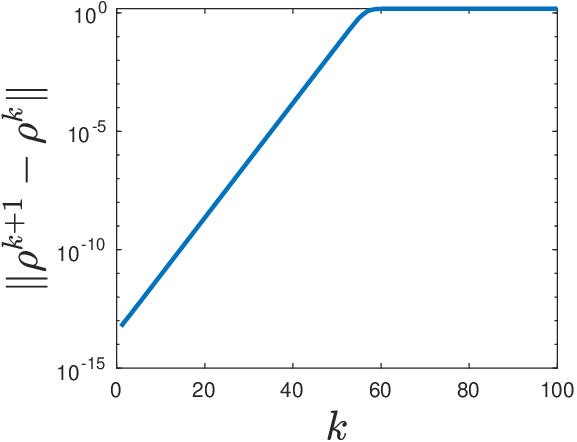}
    \includegraphics[width=5cm]{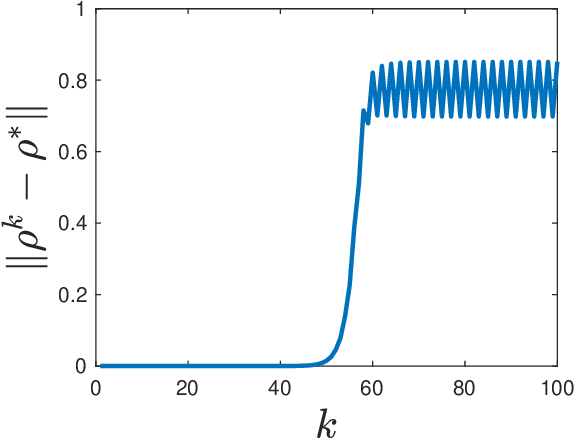}
    
\caption{Numerical results of Example \ref{subsec: num fixpt}. The fixed-point algorithm diverges even starting from the Nash Equilibrium.
Left: the semilog plot of $g^{(k)}$.
Center: consecutive residue $\|\hat{\rho}^{(k)}-\rho^{(k-1)}\|_{\calG}$.
Right: the error $\|\rho^{(k)}-\rho^*\|_{\calG}$.}
\label{fig: fixpt plots}
\end{figure}

We first show that fixed-point iteration may diverge under a very simple and ideal setting.
On $\bbT$, consider a MFG with $H(x,p)=\frac{1}{2}|p|^2,f(x,\rho)=\rho(x)+\exp(\sin(2\pi x))$ and $f_T(x,\rho)=0$. 
This is a potential MFG with a monotone interaction cost and therefore admits at most one Nash Equilibrium.
Let $\rho_0(x)=\rho_G(x;0,0.2)$. 
We first apply fictitious play to solve the problem, and the algorithm converges in 38 iterations, reaching a residue of $10^{-13}$. This gives us a numerical equilibrium $(\rho^*,v^*)$ up to machine accuracy. 
We then initialize the fixed-point iteration with the equilibrium and report the results in Figure \ref{fig: fixpt plots}. The plots show that the fixed-point iteration diverges, even starting from the Nash Equilibrium. In fact, it oscillates between two states in the end.

\subsection{Demonstration of convergence analysis.}
\label{subsec: num cvg}

In this example, we consider a problem on $[-5,5]$ to verify our convergence analysis of the fictitious play, specifically Proposition \ref{lem: stability} and the linear convergence.
We take non-zero viscosity $\nu=0.1$, quadratic Hamiltonian $H(x,p)=\frac{1}{2}|p|_2^2$, local interaction cost $f(x,\rho)=\rho(x)$, fixed terminal cost $f_T(x,\rho)=0$ and choose the initial density $\rho_0(x) = \rho_G(x;0,0.5)$.
For the discrete system, we take $n_x=10^4$ points in the space and $n_t=30$ in time. We choose the numerical viscosity to be $\nu_n=1$.
The algorithm is initialized with $\phi^{(0)}=\mathbf{0}$ which implies a zero velocity field.
In this example, we take the weight $\delta_k=0.5$ being a constant along iterations and set the stopping criteria $\epsilon=10^{-12}$.

\begin{figure}[htbp]
    \centering
    \includegraphics[width=5cm]{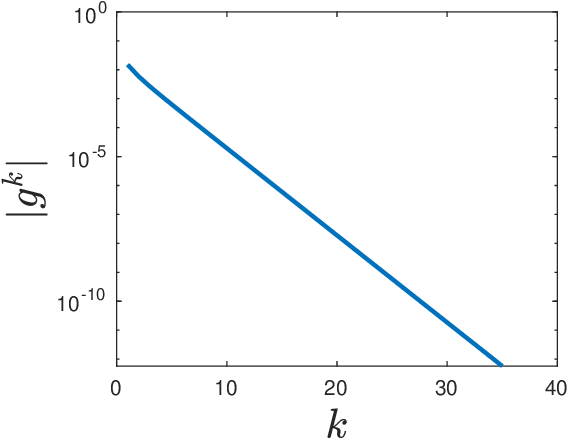}
    \includegraphics[width=5cm]{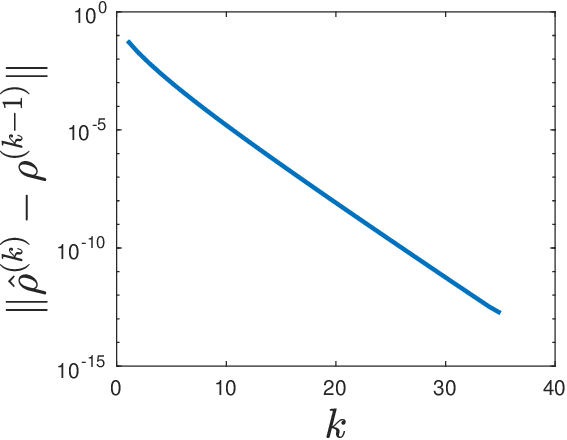}
    \includegraphics[width=5cm]{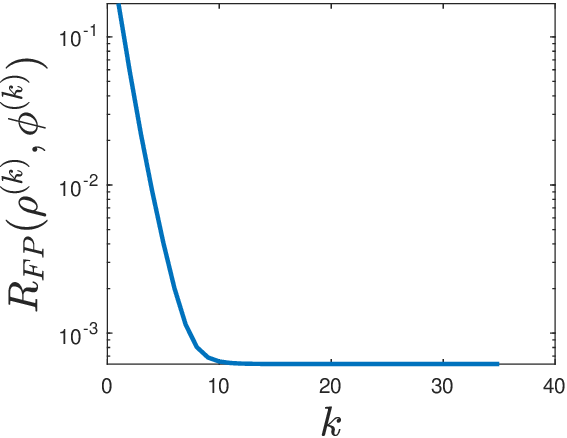}
    
    \hspace{.1cm}
    \includegraphics[width=5cm]{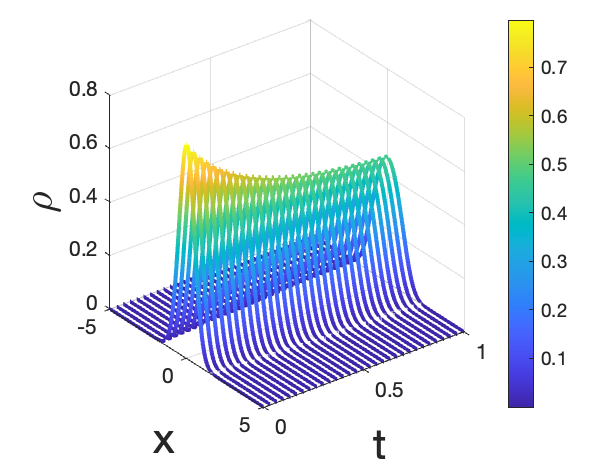}
    \includegraphics[width=5cm]{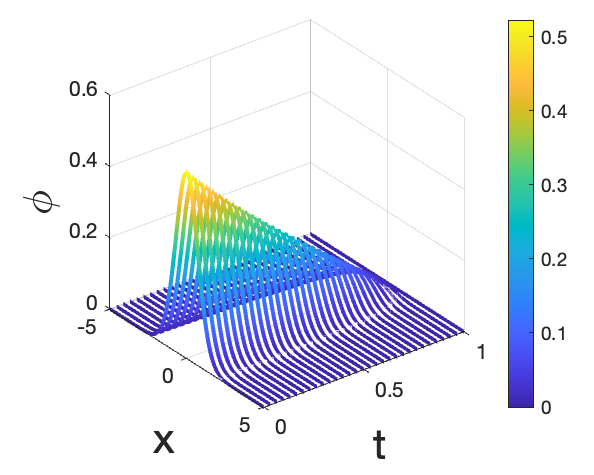}
    \hspace{5cm}
    
\caption{Numerical results of Example \ref{subsec: num cvg}. The algorithm converges in 35 iterations in 2.29s.
Top row: the semilog plot of $g^{(k)}$ (left), consecutive residue $\|\hat{\rho}^{(k)}-\rho^{(k-1)}\|_{\calG}$ (center) and Fokker-Planck residue $\left\| \resfp\left( \rho^{(k)},\hat{\phi}^{(k+1)} \right) \right\|_{\calG}$ (right) versus the number of iteration.
Bottom row: the illustrations of $\rho^{(K)}$ (left) and $\phi^{(K)}$ (right).}
\label{fig: fixphiT plots}
\end{figure}

\begin{figure}[htbp]
    \centering
    \includegraphics[width=7cm]{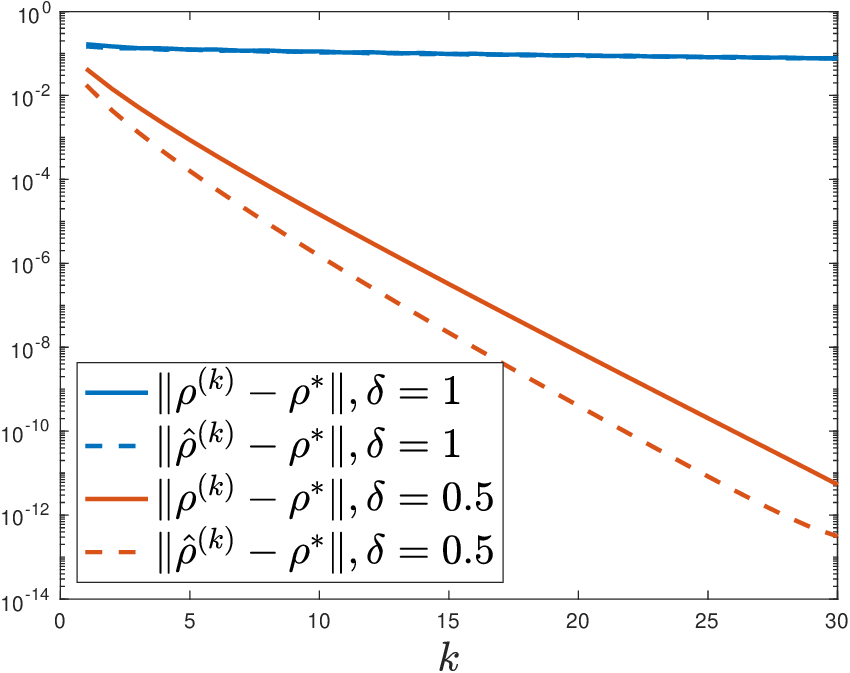}
    \includegraphics[width=7cm]{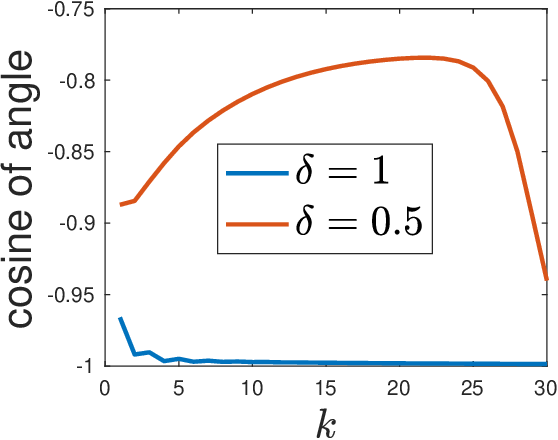}
\caption{Numerical results of Example \ref{subsec: num cvg}.
Left: Comparison of $\|\rho^{(k)}-\rho^*\|_{\calG}$ and $\|\hat{\rho}^{(k)}-\rho^*\|_{\calG}$. 
Right: the plot of $\frac{\langle \rho^{(k)}-\rho^*,\hat{\rho}^{(k)}-\rho^* \rangle_{\calG}}{\|\rho^{(k)}-\rho^*\|_{\calG}\|\hat{\rho}^{(k)}-\rho^*\|_{\calG}}$.}
\label{fig: fixphiT angles}
\end{figure}

Our algorithm converges in 35 iterations in 2.29s.
Figure \ref{fig: fixphiT plots} shows that the algorithm converges linearly.
We treat the output $\rho^{(K)}$ as the equilibrium $\rho^*$ and run the algorithm again with $\delta_k=0.5$ and $\delta_k=1$ for 30 iterations.
In this run, we measures $\|\rho^{(k)}-\rho^*\|_{\calG}$, $\|\hat{\rho}^{(k)}-\rho^*\|_{\calG}$ and $\frac{\langle \rho^{(k)}-\rho^*,\hat{\rho}^{(k)}-\rho^* \rangle_{\calG}}{\|\rho^{(k)}-\rho^*\|_{\calG}\|\hat{\rho}^{(k)}-\rho^*\|_{\calG}}$ and summarize the results in Figure \ref{fig: fixphiT angles}.

According to Proposition \ref{lem: stability}, $\langle \hat{\rho} - \rho^*, \rho - \rho^* \rangle \leq 0$, achieving equality only when $\hat{\rho} = \rho^*$. Actually the plot of $\frac{\langle \rho^{(k)}-\rho^*,\hat{\rho}^{(k)}-\rho^* \rangle_{\calG}}{\|\rho^{(k)}-\rho^*\|_{\calG}\|\hat{\rho}^{(k)}-\rho^*\|_{\calG}}$ in Figure \ref{fig: fixphiT angles} demonstrates that 
$\hat{\rho}^{(k)} - \rho^*$ and $\rho^{(k)} - \rho^*$ are almost in perfect opposite directions which makes fictitious play stable and efficient when using an appropriate convex weighted average of $\hat{\rho}^{(k)}$ and $\rho^{(k)}$ for an update. 
Additionally, the comparison between $\delta_k = 0.5$ and $\delta_k = 1$ suggests that more aggressive updates may not be a good choice. In particular, $\delta_k = 1$ corresponds to the fixed point iteration, which generates a sequence that jumps left and right without converging for this particular example. 

\subsection{Convergence rate independent of mesh step}
\label{subsec: num ind mesh}

This example aims to demonstrate that the convergence of fictitious play is independent of the mesh size and the discretization accuracy, which is a motivation for using hierarchical grids and the assumption for the complexity analysis in Table \ref{tab: complexity}. On the other hand, since the constraint involves differential operators, if an optimization method is used to solve a potential MFG, the condition number of the discretized system becomes larger when the mesh size becomes smaller, which causes slower convergence. 

\begin{figure}[htbp]
    \centering
        \includegraphics[width=5cm]{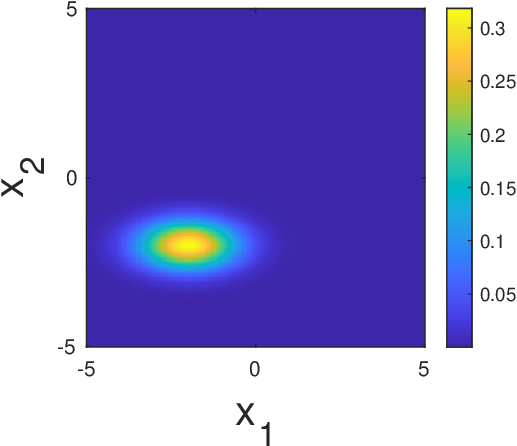}
        \includegraphics[width=5cm]{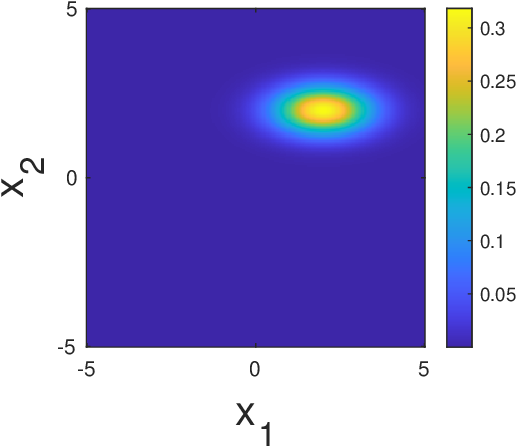}
        \includegraphics[width=5cm]{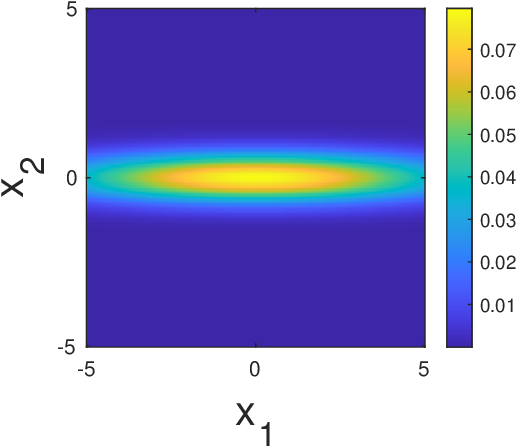}        
    \caption{Illustrations of Example \ref{subsec: num ind mesh}. Left: initial density $\rho_0$, center: desired terminal density $\rho_1$, right: convolutional kernel $K$}
    \label{fig: cvg indp mesh illus}
\end{figure}

We consider a non-potential MFG on $[-5,5]^2$ with standard Hamiltonian $H(x,p)=\frac{1}{2}|p|^2$ and viscosity $\nu=1$.
The interaction cost is a convolution $f(x,\rho)=10\int_{\bbT^d}K(x-y)\rho(y)\dd y$, where the kernel is a non-isometric Gaussian $K(x)=\rho_G(x_1;0,4)\rho_G(x_2;0,0.5)$ as illustrated in Figure \ref{fig: cvg indp mesh illus}. 
Due to this interaction cost, the population tends to avoid concentration along the $x_1$ direction.
We choose the initial and desired terminal density to be $\rho_0(x)=\rho_G(x_1;-2,1)\rho_G(x_2;-2,0.5)$, $\rho_1(x)=\rho_G(x_1;2,1)\rho_G(x_2;2,0.5)$.
As illustrated in Figure \ref{fig: cvg indp mesh illus}, they are non-isometric Gaussians as well.
Although the initial and terminal distributions have a larger variance along the $x_1$ direction, the interaction kernel induces a larger cost along the $x_1$ direction. To reduce the interaction cost, the players tend to avoid aggregating along the $x_1$ direction.
The terminal cost is $f_T(x,\rho)=150(\rho-\rho_1)$, which pushes the terminal density to the desired one.

\begin{figure}[htbp]
    \centering
        \includegraphics[width=5cm]{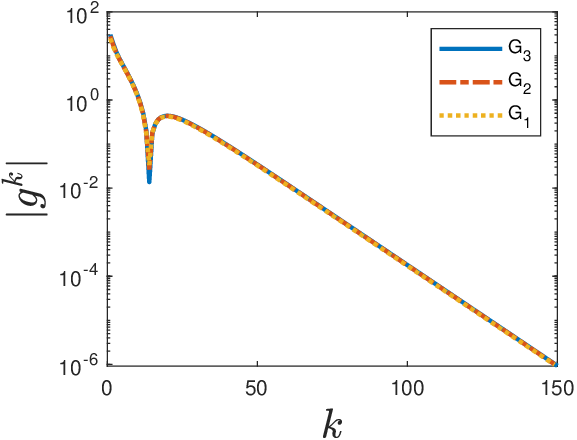}
        \includegraphics[width=5cm]{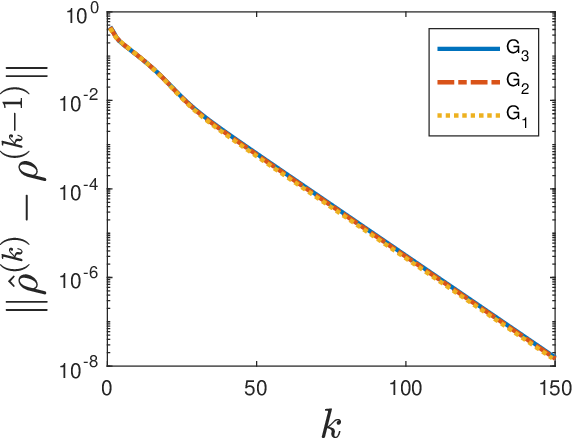}
        \includegraphics[width=5cm]{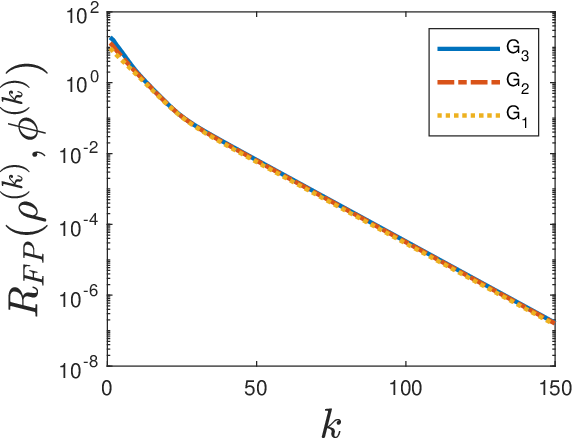}

    \caption{Numerical results of Example \ref{subsec: num ind mesh}: The algorithm converges on grids $\calG_1,\calG_2,\calG_3$ in 149, 150 and 150 iterations and 34.15s, 324.67s and 3592.56s, respectively. The figure shows the semilog plot of $g^{(k)}$ (left), consecutive residue $\|\hat{\rho}^{(k)}-\rho^{(k-1)}\|_{\calG}$ (center) and Fokker-Planck residue $\|\resfp(\rho^{(k)},\hat{\phi}^{(k+1)}\|_\calG$ (right) versus the number of iterations.}
    \label{fig: cvg indp mesh plots}
\end{figure}

\begin{figure}[htbp]
    \centering
    \includegraphics[width=4cm]{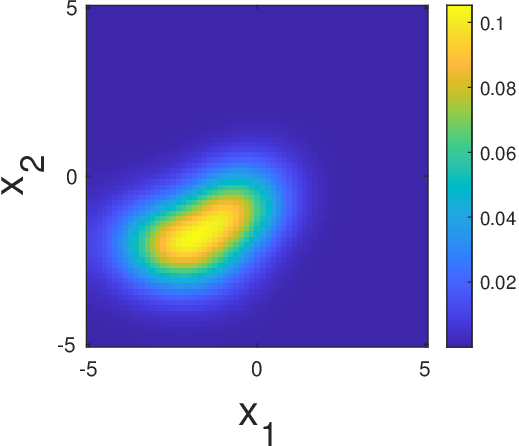}
    \includegraphics[width=4cm]{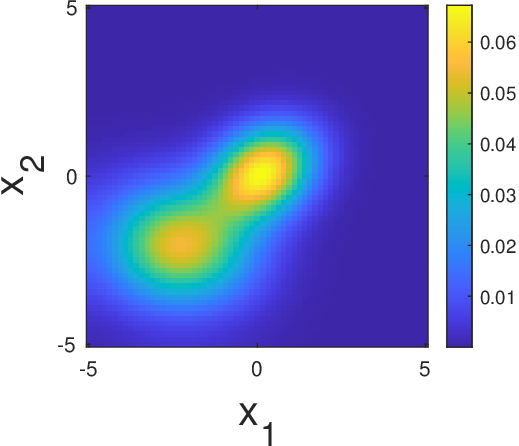}
    \includegraphics[width=4cm]{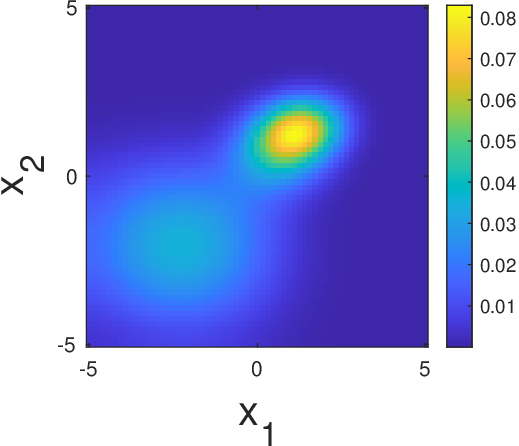}
    \includegraphics[width=4cm]{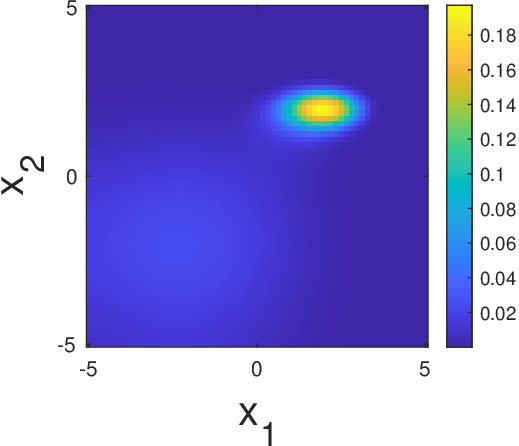}\\
    \includegraphics[width=4cm]{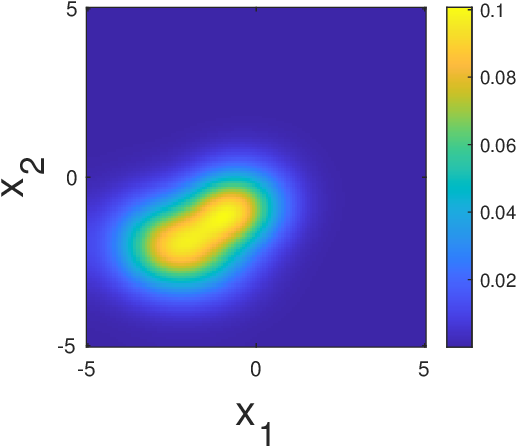}
    \includegraphics[width=4cm]{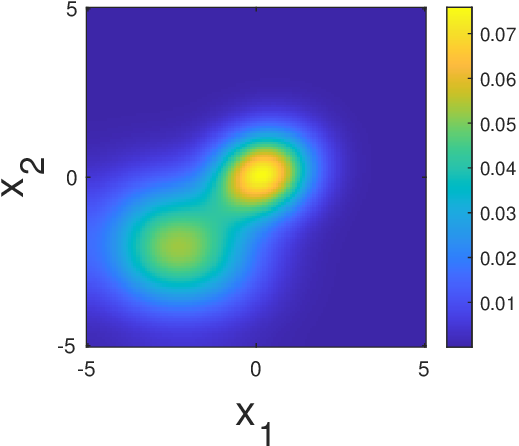}
    \includegraphics[width=4cm]{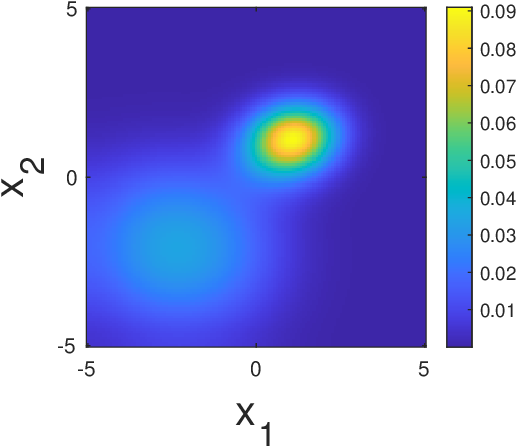}
    \includegraphics[width=4cm]{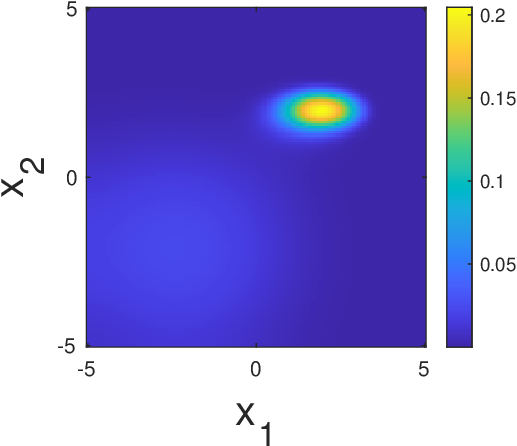}\\
    \includegraphics[width=4cm]{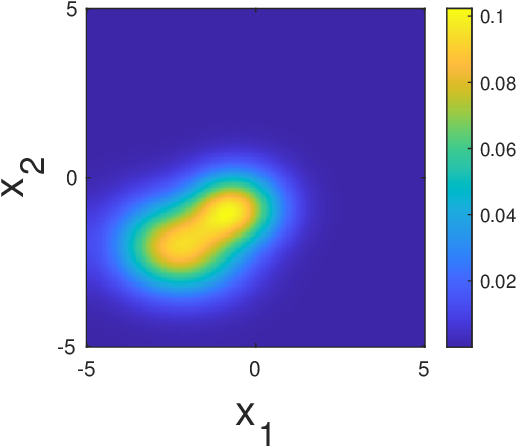}
    \includegraphics[width=4cm]{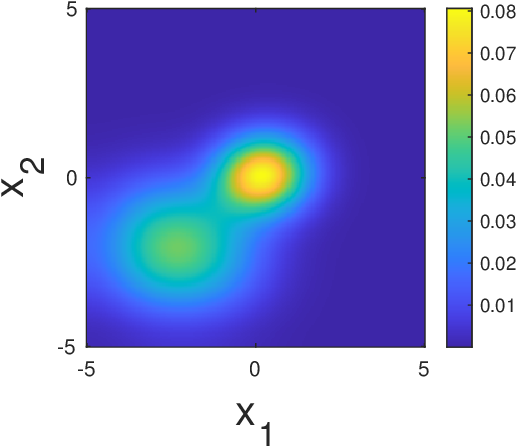}
    \includegraphics[width=4cm]{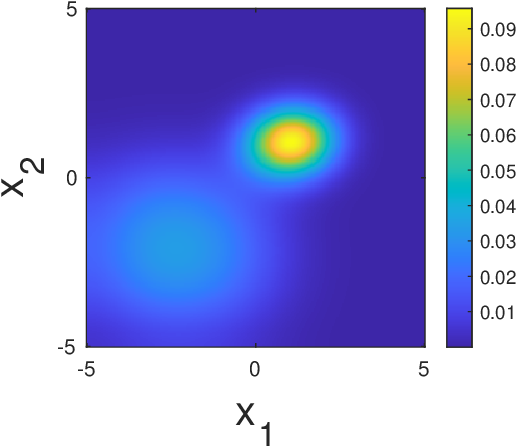}
    \includegraphics[width=4cm]{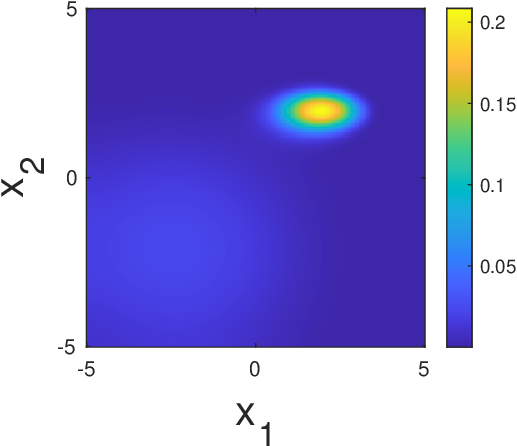}\\
\caption{Numerical results of Example \ref{subsec: num ind mesh}: illustration of $\rho^{(K)}$. Left to right: $t=0.25,0.5,0.75,1$, top to bottom: $l=1,2,3$.}
\label{fig: cvg indp mesh rho}
\end{figure}

In numerical simulation, we set numerical viscosity to be $\nu_n=0$, and weight to be a constant $\delta_k=0.1$.
We test our algorithm on different grids $\calG_l,l=1,2,3$ ($n_{x_1}=n_{x_2}=32\times 2^l$ and $n_t=4\times 2^l$) with the same initialization and set the tolerance to be $10^{-6}$.
Numerical results in Figure \ref{fig: cvg indp mesh plots} show that the convergence behavior of the algorithm on different grids is almost the same.
The algorithm converges on $\calG_1,\calG_2,\calG_3$ in 149, 150 and 150 iterations and 34.15s, 324.67s and 3592.56s, respectively.
And the numerical solution of the Nash Equilibrium $\rho$ in Figure \ref{fig: cvg indp mesh rho} also looks very similar and meets our expectations.

\subsection{Acceleration with back-tracking line-search}
\label{subsec: num btls}

Next, we illustrate the effectiveness of the line-search strategy through a MFG with a non-quadratic Hamiltonian on $[-1,1]$.
Let $H(x,p) = h(x) + \frac{1}{\gamma}|p|^\gamma$ where$\gamma=1.5$ and $h(x)=-\sin(2\pi x)$.
We take initial distribution to be Gaussian $\rho_0(x) = \rho_G(x;0,0.1)$, viscosity $\nu=0.1$, non-local interaction cost $f(x,\rho)=100(I-\Delta)^{-2}\rho$ and terminal cost $f_T(x,\rho)=\rho(x)$.

\begin{figure}[htbp]
    \centering
    \includegraphics[width=5cm]{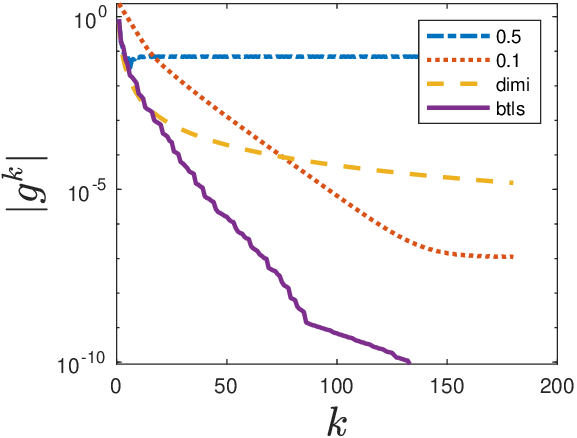}
    \includegraphics[width=5cm]{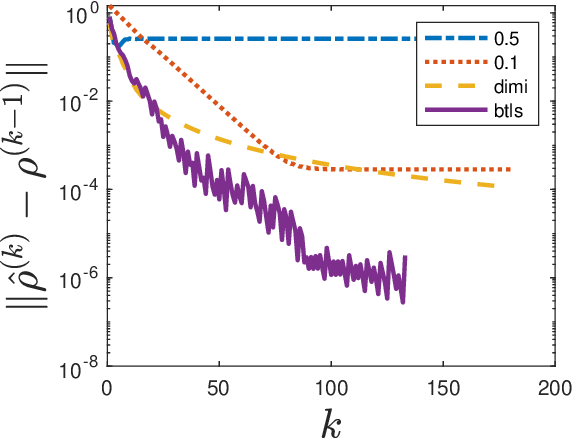}
    \includegraphics[width=5cm]{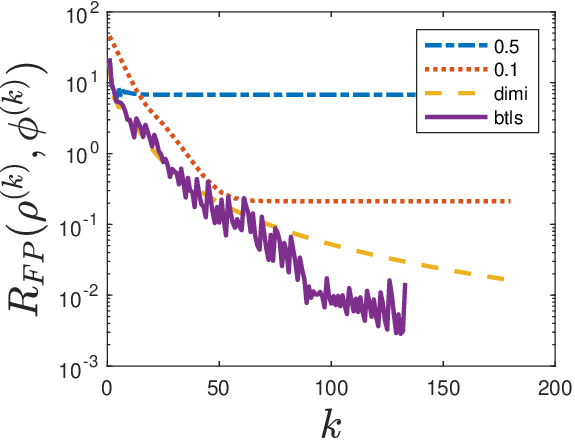}
    
    \includegraphics[width=5cm]{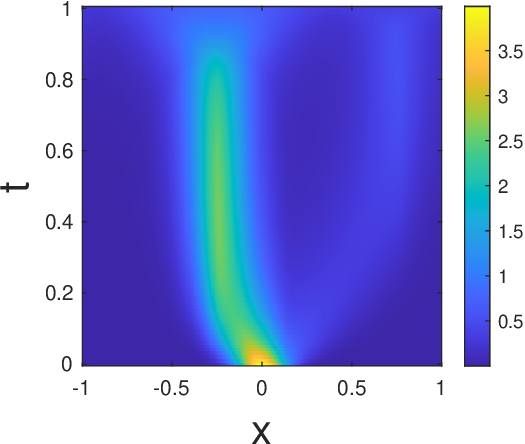}
    \includegraphics[width=5cm]{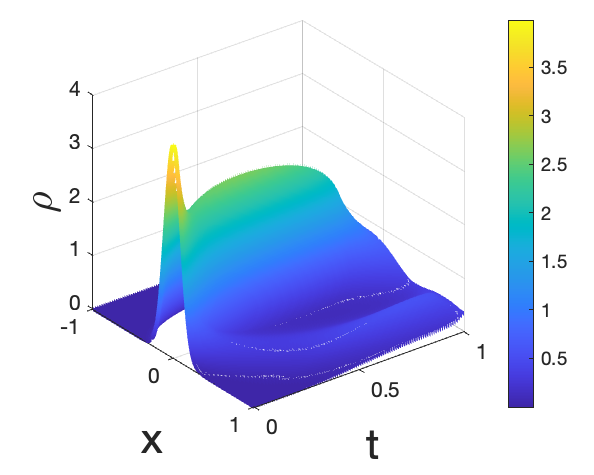}
    \includegraphics[width=5cm]{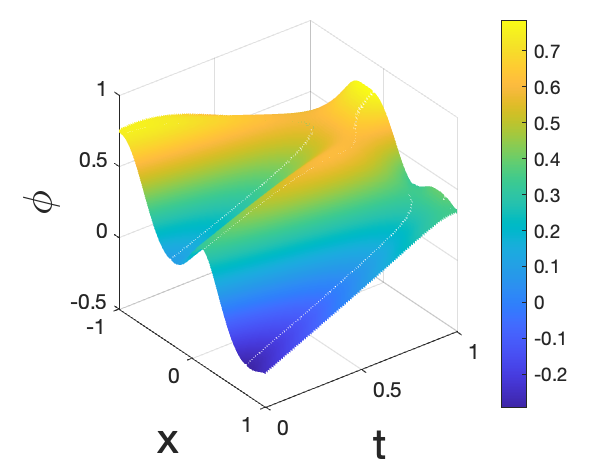}
\caption{Numerical results of Example \ref{subsec: num btls}.
The algorithm with back-tracking line-search weight $\delta_k$ converges in 172 iterations in 28.64s.
Top row: the semilog plot of $g^{(k)}$ (left), consecutive residue $\|\hat{\rho}^{(k)}-\rho^{(k-1)}\|_{\calG}$ (center) and FP residue $\|\resfp(\rho^{(k)},\hat{\phi}^{(k+1)})\|_{\calG}$ (right) versus the number of iteration.
Bottom row: the illustration of $\rho^{(K)}$ (left and center) and $\hat{\phi}^{(K)}$ (right) with back-tracking line-search weight $\delta_k$.
}
\label{fig: nonquad plots}
\end{figure}

We discretize the space with $n_x=2\times10^3$ segments, the time domain with $n_t=100$ segments, and take numerical viscosity $\nu_n=0$.
The algorithm initializes with a zero velocity field and the tolerance is $\epsilon=10^{-10}$. 
We select $\delta_k$ by the back-tracking line-search method in Section \ref{subsec: acc btls} and the parameters are $\delta=1,\beta=0.5$ and $\zeta=0.8$.
The algorithm with line-search weight converges in 133 iterations in 23.13s.
We also run the algorithm with diminishing weight $\delta_k=\frac{2}{k+2}$ and constant weight $\delta_k=0.1$, $\delta_k=0.5$ for 180 iterations and compare the convergence in Figure \ref{fig: nonquad plots}, where ``dimi'' stands for diminishing weight and ``btls'' stands for back-tracking line-search weight.
Figure \ref{fig: nonquad plots} shows that back-tracking line-search efficiently improves the convergence compared to diminishing weight or non-optimal pre-select constant weight (0.1), and it also avoids getting stuck at some point with an aggressive pre-select constant weight (0.5).
We remark that the oscillation in the consecutive residue and Fokker-Planck residue is reasonable since the back-tracking line search only guarantees the decay of $g^{(k)}$ and how the consecutive residue is bounded by $g^{(k)}$ depends on the local convexity property between $\rho^{(k-1)}$ and $\hat{\rho}^{(k)}$, which can be different for different $k$.
Since the weight sequence generated by the back-tracking line search is not smooth, oscillations occur.
This oscillation is not observed for diminishing weight or constant weight because the weight sequences are smoother there and therefore generate a smoother sequence of $\rho^{(k)}$.

\subsection{Stabilization with hierarchical grid strategy.}
\label{subsec: num multigrid stab}

In this example, we demonstrate that the hierarchical grid strategy helps to stabilize the algorithm through a first-order MFG ($\nu=0$).

Let $\mu_1(\rho),\mu_2(\rho)$ be first and second-order moments of the distribution $\rho$.
We consider a MFG on $\bbR$ with
initial distribution density $\rho_0(x) = \rho_G(x;c,\sigma_0)$,
Hamiltonian $H(x,p)=\frac{1}{2}|p|_2^2$, 
interaction cost $f(x,t,\rho)=\frac{1}{2}\left(\alpha^2+\frac{\nu^2}{(\mu_2(\rho)-\mu_1^2(\rho))^2}\right) \left( x-\mu_1(\rho) \right)^2 + 2\alpha x + \left( \frac{1}{2}(2at+b)^2 + \alpha\nu - \frac{\nu^2}{\mu_2(\rho)-\mu_1^2(\rho)} \right)$,
terminal cost $f_T(x,\rho)=-(2a+b)x - \frac{1}{2}\left( \alpha - \frac{\nu}{\mu_2(\rho)-\mu_1^2(\rho)} \right) \left( x-\mu_1(\rho) \right)^2$.
Here $a,b,c,\sigma_0,\alpha$ are pre-selected constants.
Let $\mu(t)=at^2+bt+c$ and $\sigma(t)=\sigma_0\exp{(\alpha t)}$.
A solution to this system is
$$\rho(x,t) = \rho_G(x;\mu(t),\sigma(t)),\quad
 \phi(x,t) = -(2at+b)x - \frac{1}{2}\left( \alpha - \frac{\nu}{\sigma(t)} \right) \left( x-\mu(t) \right)^2.$$

In numerical experiments, we take $a=0,b=0.5,c=-0.25$ and $\sigma_0=0.5,\alpha=-0.1$. Therefore we expect the solution $\rho(\cdot,t)$ to be a Gaussian distribution with mean $0.5t-0.25$ and standard deviation $0.5\times e^{-0.1t}$.
We set the space boundary as $[-5,5]$, take $n_x=2^{14}$ points in space, and $n_t=2^9$ points in the time domain.
The numerical viscosity is set to be 1 and the tolerance is set to be $\epsilon=10^{-6}$.
We initialize the algorithm with $\rho^{(0)}(x,t)=\rho_0(x)$ and run it with weight $\delta_k=0.1$. 

Algorithm \ref{alg: disct ficplay} fails when we impose the Dirichlet boundary condition for both the HJB and FP equations. This failure occurs because the incorrect initial guess is incompatible with the boundary condition, and in the absence of physical viscosity, the solution of the HJB equation becomes highly non-smooth near the space boundary (boundary layer). The lack of smoothness propagates through the fictitious play iterations, ultimately causing instability and divergence.
To resolve this issue, we switch to Algorithm \ref{alg: multigrid} with $L=5$. The relatively large numerical viscosity on the coarse grid allows the algorithm to first generate a relatively smooth solution on the coarse grid as well as a good initial guess for the next level finer grid, eventually leading to convergence to the true solution of the system.

\begin{figure}[htbp]
    \centering
    \includegraphics[width=5cm]{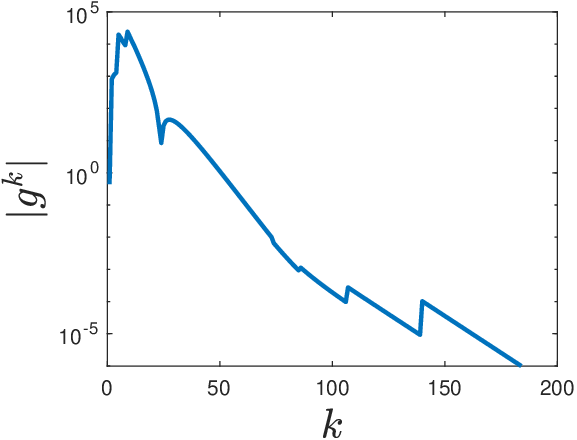}
    \includegraphics[width=5cm]{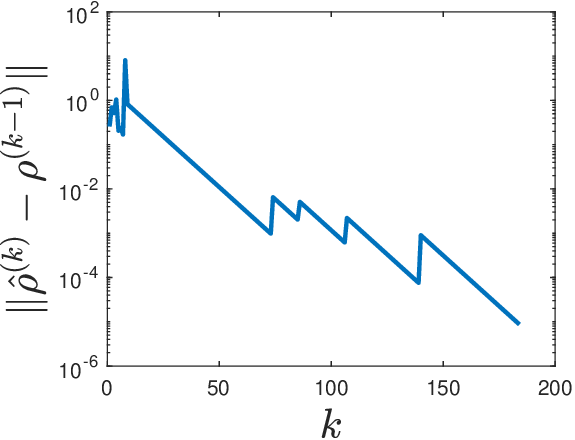}
    \includegraphics[width=5cm]{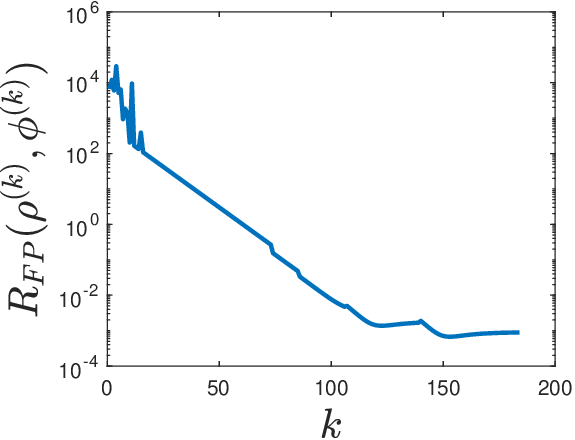}

    \includegraphics[width=5cm]{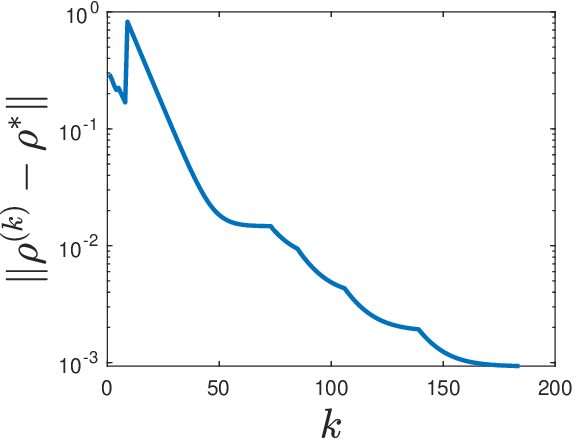}
    \includegraphics[width=5cm]{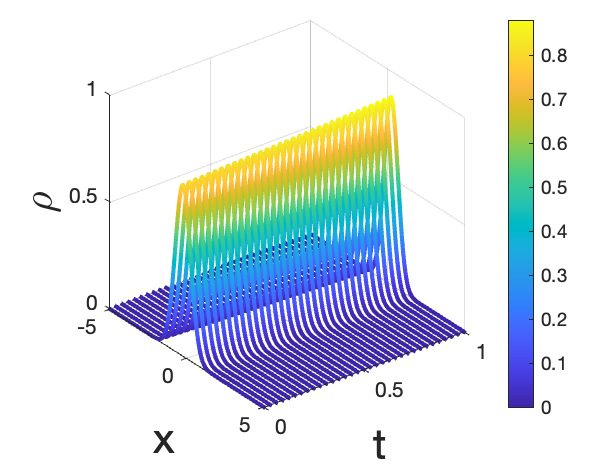}
    \includegraphics[width=5cm]{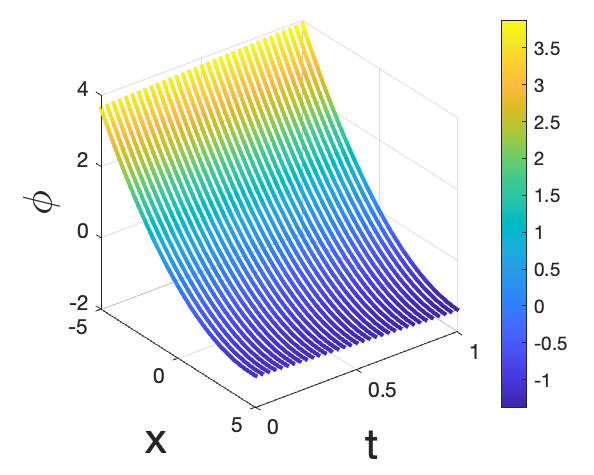}
\caption{Numerical results of Example \ref{subsec: num multigrid stab}. The algorithm finishes in 97.68s. 
Top row: the semilog plot of $g^{(k)}$ (left), consecutive residue $\|\hat{\rho}^{(k)}-\rho^{(k-1)}\|_{\calG}$ (center) and Fokker-Planck residue $\left\| \resfp\left( \rho^{(k)},\hat{\phi}^{(k+1)} \right) \right\|_{\calG}$ (right) versus the accumulated number of iteration.
Bottom row: the semilog plot of $\|\rho^{(k)}-\rho^*\|_{\calG}$ (left) versus the accumulated number of iterations, the illustrations of $\rho^{(K)}$ (center) and $\phi^{(K)}$ (right) on the fine grid.}
\label{fig: movegau1 ml plots}
\end{figure}

The first row of Figure \ref{fig: movegau1 ml plots} shows the gain function, the consecutive residue, and the FP residue against the accumulated number of iterations. 
The jumps are caused by the grid changes when the solution on a coarse grid is interpolated to a finer grid, which will generate a larger error initially due to a change of discretization on different mesh sizes.
Since the exact solution is available in this example, we plot the error against the accumulated number of iterations in the second row. It shows that the hierarchical grid fictitious play effectively converges to the Nash Equilibrium.

\subsection{Acceleration with hierarchical grid strategy}
\label{subsec: num multigrid acc}

In addition to stabilizing the algorithm, the hierarchical grid strategy also helps accelerate the algorithm. With the same cost functions in Example 4, we choose a different set of parameters and boundary conditions to illustrate the acceleration.

We take $a=6,b=-5,c=0$ and $\sigma_0=1,\alpha=-0.5$, viscosity $\nu=0.1$ and weight $\delta_k=0.25$. The boundary conditions are homogeneous Neumann boundary conditions. The numerical viscosity is set to be 0 since there is a physical viscosity.
The other settings are the same as in Example 4.

\begin{figure}[htbp]
    \centering
    \includegraphics[width=5cm]{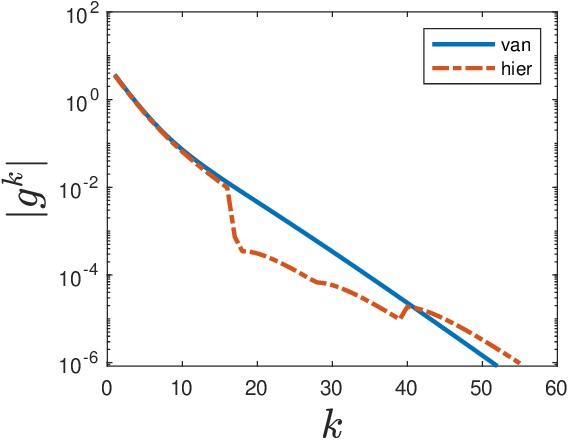}
    \includegraphics[width=5cm]{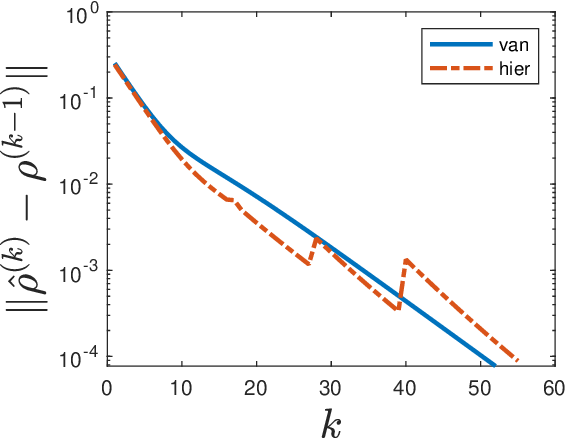}
    \includegraphics[width=5cm]{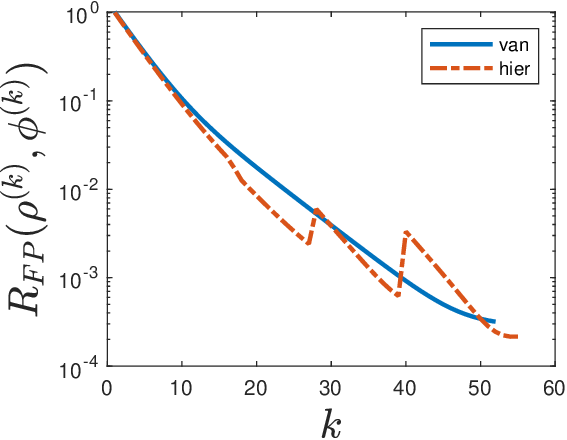}

    \includegraphics[width=5cm]{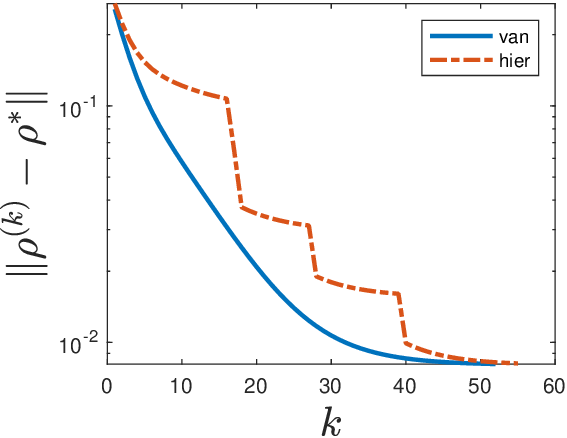}
    \includegraphics[width=5cm]{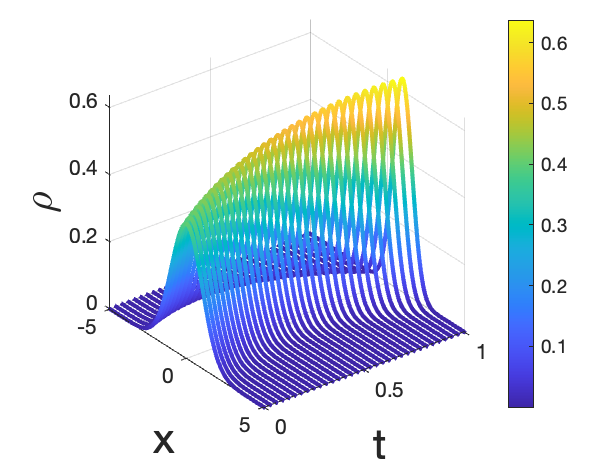}
    \includegraphics[width=5cm]{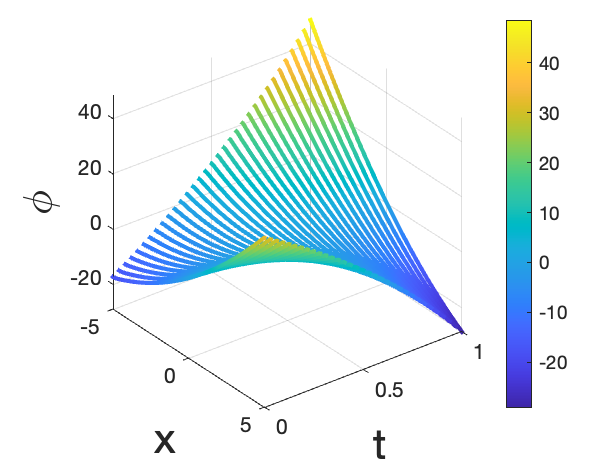}
\caption{Numerical results of Example \ref{subsec: num multigrid acc}. The vanilla algorithm finishes in 52 iterations (all on the finest grid) in 72.05s. The hierarchical grid algorithm finishes in 55 iterations (16 iterations on the finest grid) in 30.13s. 
Top row: the semilog plot of $g^{(k)}$ (left), consecutive residue $\|\hat{\rho}^{(k)}-\rho^{(k-1)}\|_{\calG}$ (center) and Fokker-Planck residue $\left\| \resfp\left( \rho^{(k)},\hat{\phi}^{(k+1)} \right) \right\|_{\calG}$ (right) versus the accumulated number of iteration.
Bottom row: the semilog plot of $\|\rho^{(k)}-\rho^*\|_{\calG}$ (left) versus the accumulated number of iterations, the illustrations of $\rho^{(K)}$ (center) and $\phi^{(K)}$ (right) on the fine grid.}
\label{fig: movegau2 ml plots}
\end{figure}

We compare the results of Algorithm \ref{alg: disct ficplay} and Algorithm \ref{alg: multigrid} in Figure \ref{fig: movegau2 ml plots}.
It takes 72.06s for Algorithm \ref{alg: disct ficplay} to converge in 52 iterations.
And Algorithm \ref{alg: multigrid} only takes 30.13s, which is less than half of Algorithm \ref{alg: disct ficplay}, to converge. This is because the outputs on the coarse grid give a good initial guess on the fine grid and therefore Algorithm \ref{alg: multigrid} only spends 16 iterations on the fine grid.

\subsection{Mean-field planning with hierarchical grid strategy}
\label{subsec: num multigrid mfp}

In this example, we solve a mean-field planning (MFP) problem \eqref{eq: mfp} in $[-5,5]^2\subset\bbR^2$. Different from \eqref{eq: mfg}, both boundary conditions of a MFP are on population density.
\begin{equation}
\left\{
\begin{aligned}
    -\partial_t\phi-\nu\Delta\phi+H(x,\nabla\phi(x,t))=f(x,\rho(t)),&   \\
    \partial_t\rho-\nu\Delta\rho-\nabla\cdot(\rho D_pH(x,\nabla\phi))=0,& 
    \quad \rho(x,0)=\rho_0(x),\rho(x,1)=\rho_1(x).   
\end{aligned}
\right.
\label{eq: mfp}
\end{equation}
As discussed in \cite{Achdou2012mfp}, one can solve the MFP \eqref{eq: mfp} by considering the penalized MFG \eqref{eq: mfgpp} with large positive $\eta$.
\begin{equation}
\left\{
\begin{aligned}
    -\partial_t\phi-\nu\Delta\phi+H(x,\nabla\phi(x,t))=f(x,\rho(t)),&\quad \phi(x,1) = \eta(\rho(x,1)-\rho_1(x)),   \\
    \partial_t\rho-\nu\Delta\rho-\nabla\cdot(\rho D_pH(x,\nabla\phi))=0,& 
    \quad \rho(x,0)=\rho_0(x).   
\end{aligned}
\right.
\label{eq: mfgpp}
\end{equation}
Let the viscosity $\nu=1$.
We would like to solve the problem with interaction cost  $f(x,\rho(t))=\begin{cases}
    +\infty,& |x|_2^2\leq 2,\\
    0,&\text{otherwise}
\end{cases}$, which is the indicator function of a disk and is independent of the population density. And we take the initial and terminal density to be two Gaussian distributions with means $[-3,-3]$ and $[3,3]$, respectively, and both with covariance matrix $\operatorname{diag}(0.25,0.25)$.
The solution $\rho$ should avoid the disk and transport from one Gaussian distribution to the other.

We solve \eqref{eq: mfgpp} on a hierarchy of grids with increasing penalty parameters. 
On the $l$-th level grid $\calG_l$, the domain is discretized with $n_t=2\times 2^l,n_{x_1}=n_{x_2}=8\times 2^l$. 
To simulate the indicator function and $\eta\to \infty$, on the $l$-th grid, we take $\eta=50l$ and solve \eqref{eq: mfgpp} with $f(x,\rho(t))=\begin{cases}
    2\eta,& |x|_2^2\leq 2,\\
    0,&\text{otherwise}.
\end{cases}$
At each level, we run Algorithm \ref{alg: disct ficplay} until tolerance $\epsilon=1\times10^{-6}$ is achieved or a maximum of 200 iterations is achieved.
The initialization of the $l$-th layer is the prolongation of the output on the $l-1$-th layer.

\begin{figure}[htbp]
    \centering
    \includegraphics[width=4cm]{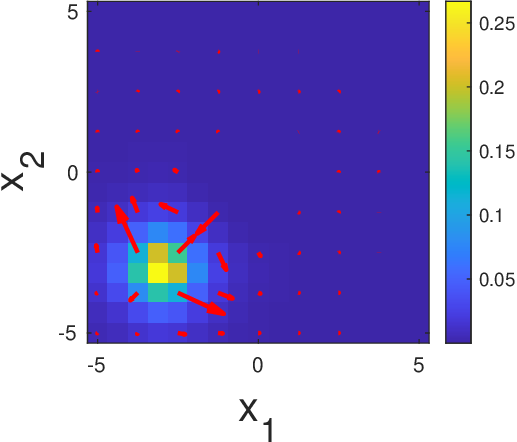}
    \includegraphics[width=4cm]{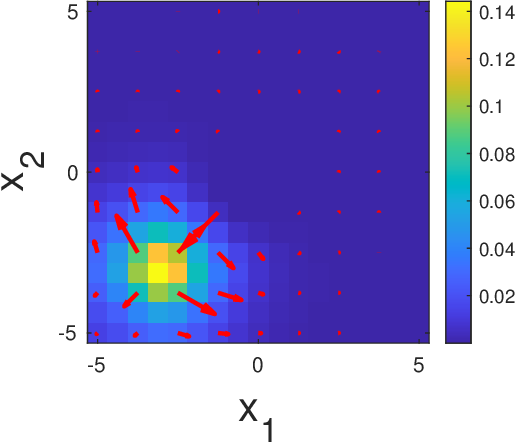}
    \includegraphics[width=4cm]{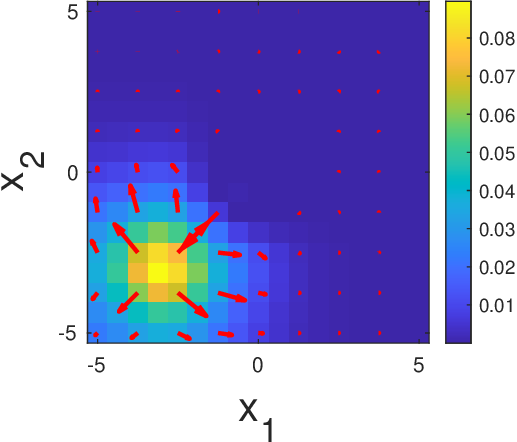}
    \includegraphics[width=4cm]{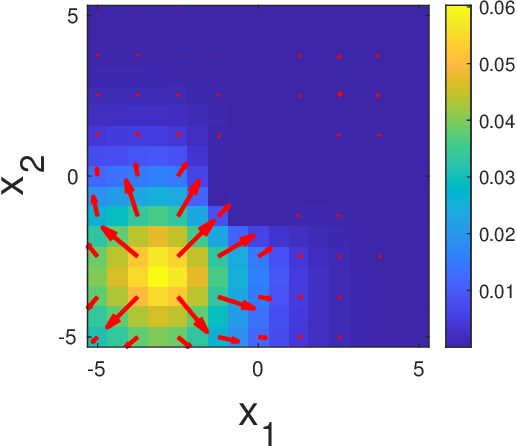}\\
    \includegraphics[width=4cm]{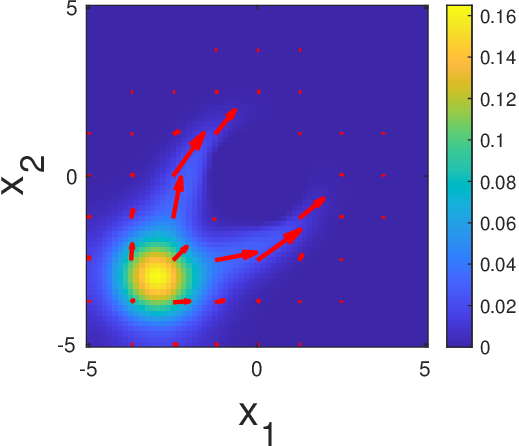}
    \includegraphics[width=4cm]{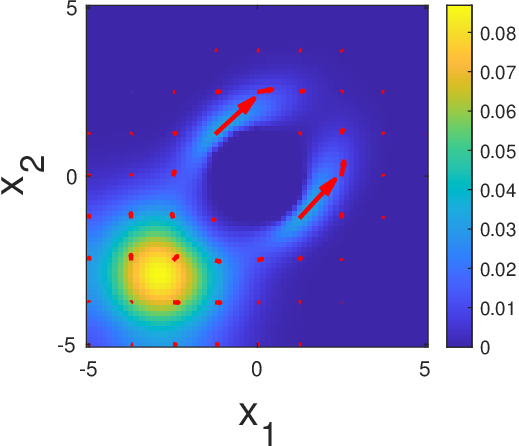}
    \includegraphics[width=4cm]{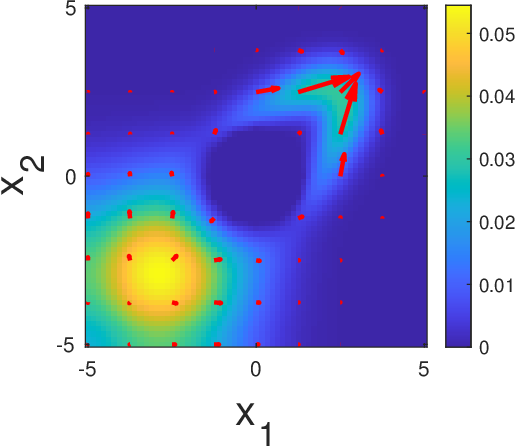}
    \includegraphics[width=4cm]{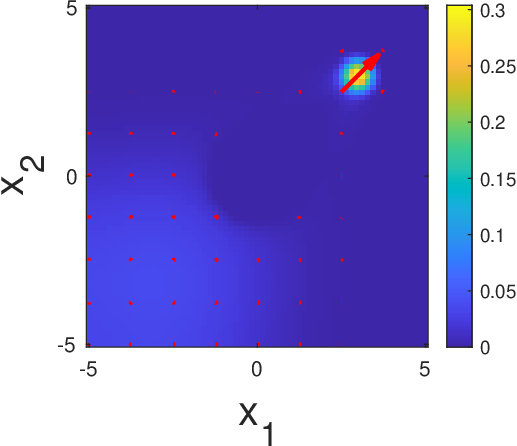}\\    
    \includegraphics[width=4cm]{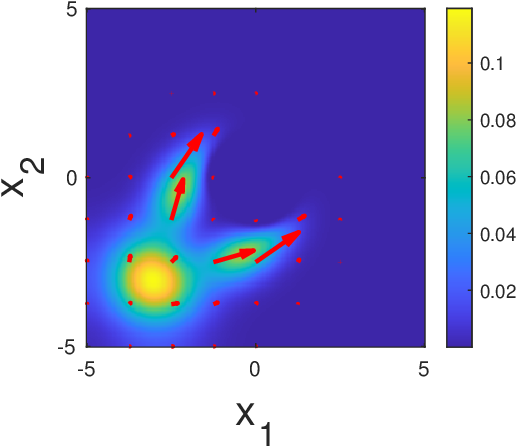}
    \includegraphics[width=4cm]{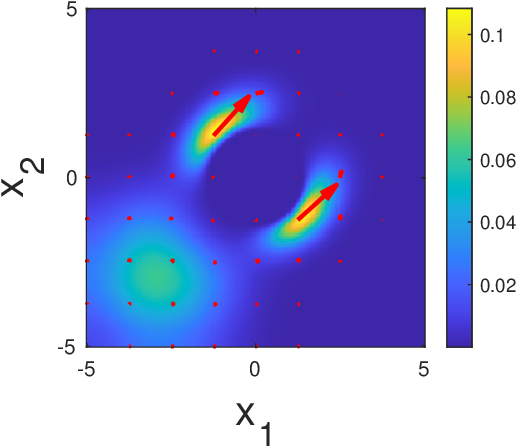}
    \includegraphics[width=4cm]{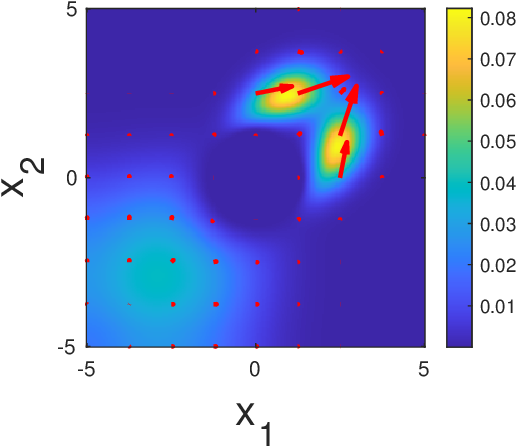}
    \includegraphics[width=4cm]{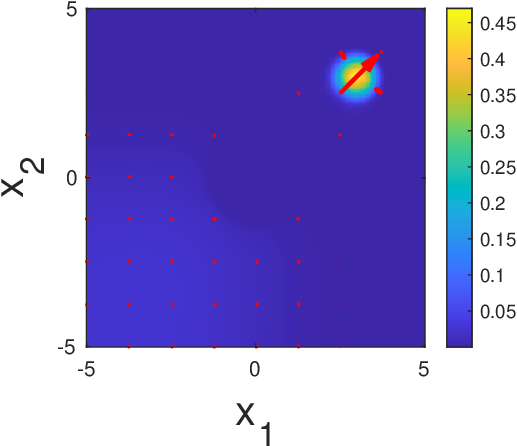}    
\caption{Numerical results of Example \ref{subsec: num multigrid mfp}: illustration of $\rho^{(K)}$. Left to right: $t=0.25,0.5,0.75,1$, top to bottom: $l=1,3,5$.}
\label{fig: mfp rho}
\end{figure}

Figure \ref{fig: mfp rho} shows the snapshots of $\rho$ for $l=1,3,5$.
As we can see, as we increase the penalty $\eta$, the mass avoids the disk and the terminal density is closer to the desired distribution.

\section{Conclusion}
\label{sec: conclu}

This paper studies the convergence and acceleration of fictitious play. 
We show that under mild assumptions, fictitious play converges linearly.
The convergence analysis does not rely on the existence of a potential and is therefore applicable to general non-potential MFGs.
Our convergence analysis shows that the weight is crucial to the convergence rate, and therefore, we propose a back-tracking line search strategy to choose the best weight in every iteration.
In addition, we speed up the algorithm by working with a hierarchy of grids to accelerate the information transfer in the time and space joint domain.
We are exploring the efficient implementation of fictitious play in high-dimensional space and the convergence analysis when the intermediate optimization cannot be solved exactly.

\section*{Acknowledgement}
J.Y. and X.C. were partially supported by NSF DMS-2237842 and Simons Foundation (grant ID: MPS-MODL-00814643).
J-G.L was supported by NSF DMS-2106988.
H. Zhao was partially supported by NSF DMS-2012860 and DMS-2309551.

\bibliography{reference}
\bibliographystyle{plain}

\end{document}

%% file: notations.tex
\newcounter{example}

\newcommand{\bbR}{\mathbb{R}}
\newcommand{\calP}{\mathcal{P}}

\newcommand{\bbT}{\mathbb{T}}

\newcommand{\dxp}{D_x^+}
\newcommand{\dxm}{D_x^-}
\newcommand{\dxc}{D_x^c}

\newcommand{\myij}{{i,n}}
\newcommand{\ipj}{{i+1,n}}
\newcommand{\imj}{{i-1,n}}
\newcommand{\ijp}{{i,n+1}}

\newcommand{\dx}{\Delta x}
\newcommand{\dt}{\Delta t}

\newcommand{\hamlf}{H^{\text{LF}}}
\newcommand{\pxp}{p^+}
\newcommand{\pxm}{p^-}

\newcommand{\reshjb}{R_{\text{HJB}}}
\newcommand{\resfp}{R_{\text{FP}}}

\newcommand{\calG}{\mathcal{G}}

\newcommand{\pro}{P_{l}^{l+1}}

\newcommand{\dd}{\mathrm{d}}
\newcommand{\half}{\frac{1}{2}}

\DeclareMathOperator*{\argmin}{argmin}

\newtheorem{theorem}{Theorem}[section]
\newtheorem{proposition}[theorem]{Proposition}
\newtheorem{lemma}[theorem]{Lemma}
\newtheorem{corollary}[theorem]{Corollary}

\theoremstyle{definition}
\newtheorem{definition}[theorem]{Definition}
\newtheorem{assumption}[theorem]{Assumption}

\theoremstyle{remark}
\newtheorem{remark}[theorem]{Remark}